\documentclass[11p,cmyk,usenames,dvipsnames]{siamart220329}
\usepackage{amsmath}
    \usepackage{amsfonts}
    \usepackage{amssymb,bbm}
   
    \usepackage{stix, wasysym}

\usepackage{amsbsy}

\usepackage{algorithm}
\usepackage{algpseudocode}
\usepackage{grffile}
\usepackage{hyperref}
\usepackage{color}
\usepackage{comment}
\usepackage{lipsum}
\usepackage{graphicx}%
\usepackage{caption}
\usepackage{subcaption}
\usepackage[normalem]{ulem}
\usepackage{enumerate}
\usepackage{epstopdf}
\usepackage[prependcaption,colorinlistoftodos]{todonotes}

\newtheorem{remark}[theorem]{Remark}	
\makeatletter 
\@mparswitchfalse%
\makeatother
\normalmarginpar

\newcommand{\Rd}{\ensuremath{{\mathbb{R}^d}}}

\newcommand{\ep}{\ensuremath{\varepsilon}}

\newcommand{\bu}{\ensuremath{\mathbf{u}}}

\newcommand{\by}{\ensuremath{\mathbf{y}}}
\newcommand{\bz}{\ensuremath{\mathbf{z}}}

\def\P{{\mathcal P}}
\def\M{{\mathcal M}}

\newcommand{\R}{{\mathord{\mathbb R}}}

\newcommand{\N}{{\mathord{\mathbb N}}}
\newcommand{\supp}{{\mathop{\rm supp\ }}}

\newcommand{\F}{\mathcal{F}}
\newcommand{\G}{\mathcal{G}}

\newcommand{\C}{\mathcal{C}}
\newcommand{\A}{\mathcal{A}}

\newcommand{\bt}{\mathbf{t}}
\newcommand{\bv}{\mathbf{v}}
\newcommand{\ba}{\mathbf{a}}
\newcommand{\bF}{\mathbf{F}}

\newcommand{\E}{\mathcal{E}}

\newcommand{\bEta}{\boldsymbol{\eta}}

\newcommand{\bnu}{\boldsymbol{\nu}}

\def\P{{\mathcal P}}

\def\epsilon{\varepsilon}

\def\F{\mathcal{F}}

\newcommand{\be}{\begin{equation}}
\newcommand{\ee}{\end{equation}}
\newcommand{\bes}{\begin{equation*}}
\newcommand{\ees}{\end{equation*}}

\newtheorem{assumption}[theorem]{Assumption}

    \makeatletter
    \newcommand{\l@chapter}{\@dottedtocline{1}{1.5em}{2.6em}}
    \newcommand{\l@section}{\@dottedtocline{1}{1.5em}{2.6em}}
    \newcommand{\l@paragraph}{\@dottedtocline{1}{1.5em}{2.6em}}
    \newcommand{\l@subparagraph}{\@dottedtocline{1}{1.5em}{2.6em}}
    \newcommand{\l@figure}{\@dottedtocline{1}{1.5em}{2.6em}}
    \newcommand{\l@table}{\@dottedtocline{1}{1.5em}{2.6em}}
    \newcommand{\l@subsection}{\@dottedtocline{2}{4.0em}{3.6em}}
    \newcommand{\l@subsubsection}{\@dottedtocline{3}{7.4em}{4.5em}}
    \makeatother

\title{A Blob Method for Mean Field Control with Terminal Constraints\thanks{Submitted to the editors
DATE. 
\funding{K. Craig's work was  supported by NSF DMS grant 2145900. K. Elamvazhuthi's work was supported by AFOSR MURI FA9550-18-1-0502. H. Lee's work was partially supported by grants NSF DMS-1952339 and TWCF0333 from the Templeton World Charity Foundation.}}}
\author{Katy Craig\thanks{Department of Mathematics, University of California, Santa Barbara
(\email{kcraig@ucsb.edu}). Corresponding author.}
\and Karthik Elamvazhuthi\thanks{Department of Mechanical Engineering, University of California, Riverside (\email{karthike@ucr.edu}).}
\and Harlin Lee\thanks{School of Data Science and Society, University of North Carolina at Chapel Hill
(\email{harlin@unc.edu}). This work was partially completed while the author was at Department of Mathematics, University of California, Los Angeles.}}
\headers{A Blob Method for Mean Field Control with Terminal Constraints}{Katy Craig, Karthik Elamvazhuthi, and Harlin Lee}

\begin{document}
\maketitle


\begin{abstract}
In the present work, we develop a novel particle method for  a general class of mean field control problems, with source and terminal constraints. Specific examples of the problems we consider include the dynamic formulation of the $p$-Wasserstein metric, optimal transport around an obstacle, and measure transport subject to acceleration controls. Unlike existing numerical approaches, our particle method is meshfree and does not require global knowledge of an underlying cost function or of the terminal constraint. A key feature of our approach is a novel way of enforcing the terminal constraint via a soft, nonlocal approximation, inspired by recent work on blob methods for diffusion equations. We prove convergence of our particle approximation to solutions of the continuum mean-field control problem in the sense of $\Gamma$-convergence. A byproduct of our result is an extension of existing discrete-to-continuum convergence results for mean field control problems to more general state and measure costs, as arise when modeling   transport around obstacles, and more general constraint sets, including controllable linear time invariant systems. Finally, we conclude by implementing our method numerically and using it to compute   solutions the example   problems discussed above. We conduct a detailed numerical investigation of the convergence properties of our method, as well as its behavior in sampling applications and for approximation of optimal transport maps.
\end{abstract}

 	\section{Introduction}   \label{sec:problem_formulation}
The goal of the present paper is to develop  a  particle method for solving the following mean-field control  problem:
\begin{align} \label{star} \tag{$\rm\star$}
&\min_{(\mu,\bu) \in \mathcal{A}(m_0,m_1)} \int_0^1 \int_{\mathbb{R}^d}\psi(\bu(y,t))d\mu_t(y)dt + \int_0^1 \int_\Rd L(y,\mu_t) d \mu_t(y) dt  \\
& \ \mathcal{A}(m_0,m_1) := \left\{ (\mu, \bu):  \mu \in AC([0,1]; \mathcal{P}_1(\Rd)) \ , \ \  \bu \in L^1_{d\mu_t \otimes dt} ( \Rd \times [0,1]; U))   \right. \nonumber \\
&\left.     \quad \qquad \qquad \qquad \partial_t \mu_t  + \nabla \cdot \left( (\bF(\cdot, \mu_t) \mu_t  + \bu(\cdot, t) \mu_t \right) = 0  , \ 
 \mu_0 = m_0 ~  \mu_1=m_1 \right\}  . \nonumber
\end{align}
In particular, we seek an evolving probability measure $\mu$ and a control $\bu$, constrained to take values in a vector space $U$, so that (i) $\mu$ evolves from an initial measure $m_0$ toward a terminal measure $m_1$, (ii)  $(\mu,\bu)$ satisfy the continuity equation constraint, subject to additional affects from a measure dependent vector field $\bF$, and (iii) $(\mu,\bu)$ accomplish this while minimizing the integrals of the control cost $\psi$ and the state and measure cost $L$. See section \ref{notationsec}  for precise definitions of the relevant function spaces and notion of solution.

Problems of the form (\ref{star}) arise in a variety of contexts, including statistical mechanics, biology, economics, and control theory; see  \cite{Seung,carmona2018probabilistic,carrillo2014derivation,carrillo2010particle,fornasier2014mean,orrieri2019variational,porretta2014planning,ruthotto2020machine, burger2021mean, bonnet2021necessary, cavagnari2022lagrangian, jimenez2020optimal, bonnet2022semiconcavity} and the references therein. More recently, problems of this form have also attracted interest in the machine learning community, in the context of sampling and generative modeling \cite{albergo2022building,kobyzev2020normalizing}.
The main contribution of the present work is a new   method for computing  optimal solutions of (\ref{star})  via a particle approximation: see equation (\ref{MFCN}) and its numerical implementation in section \ref{numericssec}. We prove convergence of this approximation, in the sense of  {$\Gamma$-convergence}, and then implement the method numerically  to compute several fundamental examples of (\ref{star}), including the dynamic formulation of the $p$-Wasserstein distance, optimal transport around obstacles, and  measure transport subject to acceleration controls.

 Traditionally, numerical methods for solving mean field control problems  with terminal constraints have relied on an \emph{Eulerian} grid-based approach, such as that proposed in Benamou and Brenier's original work on the dynamic formulation of the 2-Wasserstein distance \cite{BenamouBrenier}.  While this approach has been valuable, especially in lower-dimensional contexts \cite{achdou2012mean,elamvazhuthi2018optimal,elamvazhuthi2023dynamical}, it presents a notable limitation in high dimensions.  In contrast, our particle method is inherently meshfree, thus feasible   in high dimensions. 
 
 A second difference between our method and  classical methods is that  our   approach computes the optimal trajectory from $m_0$ to $m_1$ by updating the location of each individual particle using only \emph{local information} from $m_1$ near that particle. This is in contrast to classical grid-based approaches, for which the terminal constraint $m_1$ must be  known globally at all grid points. As a consequence of this, we believe our approach has promise as an optimal transport based sampling method,  interpreting the final particle locations at time $t=1$ as the samples of $m_1$. Furthermore, as is evident in our numerical experiments (see section \ref{gaussiansection}),  the final particle locations obtained by our method exhibit more structure than traditional iid samples, which can   be advantageous when using the samples to approximate integrals of smooth functions against $m_1$. 

Our approach is also different from existing methods for solving  (\ref{star}) based on a \emph{Lagrangian} transport map approach, in which one reformulates the problem as a Monge problem and the solution is given by an {optimal transport map}  \cite{agrachev2009optimal,hindawi2011mass,chen2016optimal}. Classical methods for computing such a map demand an explicit, closed-form expression for the cost function, which is unavailable in many scenarios, such as  optimal transport around obstacles. Furthermore, even when an optimal transport map is found, further computation is required to infer the optimal trajectories. In contrast,  our particle approach does not require global knowledge of an underlying cost function, and it directly outputs approximate optimal trajectories. Moreover, in cases where the exact optimal transport map possesses strong monotonicity properties, such as the dynamic formulation of the p-Wasserstein distance, interpolating between these particle trajectories can then give a numerical approximation of the optimal transport map. (See section \ref{erroralonggraddescent} for an error analysis of this approximation in the case $p=2$.)


Our approach is strongly inspired by recent work by Fornasier,   Lisini,   Orrieri,  and Savar{\'e} \cite{fornasier2019mean}, which studied the convergence of finite agent controls to (\ref{star}), in the mean field limit. In this work, the authors analyze when solutions of a mean field control problem \emph{without} a terminal constraint, $\mu_1 = m_1$,  may be approximated by particle solutions of a spatially discrete optimization problem. The heuristic idea behind this approach can be seen as follows: given an approximation of $m_0$ by an empirical measure, $m_0 \approx (1/N)\sum_{i=1}^N \delta_{y_i} ,$ $\{ y_i\}_{i=1}^N \subseteq \Rd$,
if the velocity field $\bv = \bF + \bu$ in the continuity equation constraint is sufficiently regular, the continuity equation constraint reduces to a system of ordinary differential equations for the trajectories of the particles, $\dot{y}_i(t) = \bv(y_i(t),t)$, where $\mu_t= (1/N) \sum_{i=1} \delta_{y_i(t)}$, and the spatial integrals in the objective function reduce to finite sums.

The main strategy of the present paper is to extend the  approach developed by Fornasier et. al. \cite{fornasier2019mean} to the setting of mean field problems with terminal constraints (\ref{star}), while at the same time generalizing the hypotheses on the state and measure cost $L$ and the constraint set $U$. Indeed, a byproduct of our result is that we prove convergence of the discrete to the continuum problems under these weaker assumptions, which is new even in the absence of a terminal constraint. On the other hand, our extension  to  control problems with a terminal constraint requires a novel approach, via a soft, nonlocal approximation of the  constraint. This is inspired by classical vortex blob methods for the Euler and Navier Stokes equations  \cite{BealeMajda1,BealeMajda2} and the aggregation equation \cite{CraigBertozzi}, which have more recently been extended to the case of diffusion equations \cite{carrillo2019blob,craig2023blob, craig2023nonlocal, carrillo2023nonlocal,lions2001methode, oelschlager1990large, burger2023porous}.
  
With this strategy in mind, we consider solutions of (\ref{star}) under the following assumptions on the source and target distributions $m_0$ and $m_1$, the range of the control map $U$, the state and measure cost $L$,    the control cost $\psi$, and the measure dependent vector field $\bF$: 
  \begin{assumption}
\label{asmp1}
We suppose that the following hold:
\begin{enumerate}[(i)]
\item $m_0, m_1 \in \P_1(\Rd)$ with $d m_1(x) = m_1(x) dx$ for $m_1 \in L^2(\Rd)$. \label{mas}
\item $U \subseteq \mathbb{R}^d$ is a subspace.
\item $L :\R^d \times \P_1(\R^d) \rightarrow [0,+\infty]$ satisfies one of the following: \label{Las}
\begin{enumerate}[(a)]
\item $L$ is jointly uniformly continuous and $\{L = +\infty\} = \emptyset$; \label{Lcts}
\item $L$ is independent of the second variable, i.e., $L(y,\mu) = L(y,\nu)$, for all $y \in \Rd$, $\mu,\nu \in \P_1(\Rd)$,   $\{ L = +\infty \}$ is open, and $L$ is continuous on $\{L < +\infty \}$.   \label{Llsc} 
\end{enumerate}
\item $\psi :U \rightarrow [0, +\infty) : u \mapsto \frac{1}{p}|u|^p$, for $p >1$, \\
or, more generally, $\psi$ may be any function satisfying Assumption \ref{asmp1prime} below;
\item   $\bF: \Rd \times \P_1(\Rd) \rightarrow \R^d$ is jointly uniformly continuous and there exist constants $C_{F},C_F' \geq 0$ so that 
\[ |\bF(y,\mu)| \leq C_F + C_{F}' \left( |y|  + \int_\Rd |z| d \mu(z) \right) , \quad \forall y \in \mathbb{R}^d , \ \mu \in \P_1(\Rd). \]
\end{enumerate}
\end{assumption}
The main distinction between the above hypotheses and previous work by Fornasier, et. al. \cite{fornasier2019mean} is item  (\ref{Llsc}), which is of   interest when the state and measure cost $L$ is used to enforce an obstacle; see section \ref{motivatingexamplessec} below. We also note that we commit a mild notational abuse in part (\ref{mas}) above: if a measure $m$ is absolutely continuous with respect to Lebesgue measure, we use $m$ to denote both the measure and its density with respect to Lebesgue.  See  also Remark \ref{notL2} below, for an approach to relaxing the assumption $m_1 \in L^2(\Rd)$ in the context of numerical simulation.

Of particular interest in the present paper, and the cases for which we obtain the strongest convergence results, are the cases of an unconstrained control $U = \Rd$ and the case of a controllable linear time invariant system, $L = 0$, $F(y, \mu) = Ay$, for which the system is controllable between any two states and $\supp m_0 \subseteq N(A)$; see Theorem \ref{NtoinftyminimizersRdLTI}.  While the former was   considered by Fornasier, et. al. \cite{fornasier2019mean}, to our knowledge the latter is new in context of particle approximation of mean field control, with or without a terminal constraint. A key motivation underlying both of these hypotheses is that they ensure the discrete optimization problem at the particle level is \emph{feasible}, in the sense that there exists an element of the constraint set for which the objective function is finite.

The remainder of the introduction is organized as follows. In Section \ref{mainresultssec}, we describe   our particle approximation and state our main results. In Section \ref{motivatingexamplessec}, we describe three key examples of mean field control problems that motivate our Assumption \ref{asmp1} above.  In Section \ref{futurework}, we outline the strategy of our approach.

\subsection{Main results} \label{mainresultssec}
In order to state our main results, we begin by introducing an equivalent formulation of our mean-field control problem (\ref{star}) in momentum coordinates. This is motivated by the fact that, in  the original formulation, $\bu$ belongs to a  function space  that depends on $\mu$, and we seek to remove this dependence.

Define the control cost functional
\begin{align*}
&\Psi : \M(  \Rd \times [0,1];U) \times C([0,1];\P_1(\Rd)) \to [0,+\infty] , \\ 
 &  \Psi(\bnu|\mu) := 
    \begin{cases}
    \int_0^1 \int_\Rd \psi(\bu(y,t))d{\mu}_t(y)dt &{\rm if}~ d\bnu(y,t) =  \bu(y,t) d\mu_t(y)   dt , \\
  +  \infty  &{\rm otherwise.}
\end{cases}
\end{align*}
Likewise, define the penalization   for the terminal constraint, $\F: \P_1(\Rd) \to [0, +\infty]$,
\begin{equation}
\F(\mu) := 
\begin{cases}
0 & \text{if}~\mu=m_1,\\
+ \infty &\text{otherwise}.
\end{cases}
\end{equation}
The original problem (\ref{star}) is equivalent to   the following mean field control problem in momentum coordinates
\begin{align} \label{mainproblem} \tag{{\rm MFC}}
&\min_{(\mu,\bnu) \in \mathcal{C}(m_0)} \mathcal{E}(\mu, \bnu) ,  \nonumber 
 \quad  \mathcal{E}(\mu,\bnu):= \G (\mu,\bnu) + \F(\mu_1) , \nonumber \\
&\ \G (\mu,\bnu) := \Psi(\bnu|\mu) + \int_0^1 \int_\Rd L(y,\mu_t) d \mu_t(y) dt  ,\nonumber \\
&\ \C(m_0):= \{ (\mu,\bnu) \in AC([0,1]; \P_1(\Rd)) \times \mathcal{M}( \Rd \times [0,1]; U) :   \nonumber \\
&\qquad \qquad  \qquad \qquad \qquad \qquad\qquad     \partial_t \mu_t  + \nabla \cdot \left( (\bF(\cdot, \mu_t) \mu_t  + \bnu_t \right) = 0 , \ \mu_0 = m_0 \}.\nonumber
\end{align}
See Lemma \ref{momequiv} below for a proof of this equivalence.

With this formulation of the   problem in hand, we now describe our approach to softening the terminal constraint and constructing the particle approximation. The first step is to replace the characteristic function $\F $ with the soft $L^2$ penalty:
\begin{align} \label{Fepsdef}
\F_\epsilon(\mu) : = \begin{cases} \frac{1}{\epsilon} \int_\Rd | \mu(y) - m_1(y)|^2 dy &\text{ if } d \mu(y) = \mu(y) dy \text{ and } \mu \in L^2(\Rd) , \\ +\infty, &\text{ otherwise.}  \end{cases}
\end{align}
This leads to the optimization problem
\begin{align} \tag{$\rm MFC_{\epsilon}$} \label{eq:main_optimrrbb}
&\min_{(\mu, \bnu) \in \C(m_0)} \E_\epsilon(\mu,\bnu) , \quad \E_\epsilon(\mu,\bnu):= \G(\mu,\bnu)+   \F_\epsilon(\mu_1) . \\ \nonumber
\end{align}

We have the following convergence result for solutions of (\ref{eq:main_optimrrbb}) to (\ref{mainproblem}) as $\epsilon \to 0$. 
\begin{proposition}[Convergence as $\epsilon \to 0$]
\label{thm:convep}
Suppose Assumption \ref{asmp1} holds and (\ref{mainproblem}) is feasible.
Suppose that $(\mu_\epsilon,\bnu_\epsilon)$ is a sequence of solutions to (\ref{eq:main_optimrrbb}). Then, there exists $(\mu,\bnu) \in \C(m_0)$ so that, up to a subsequence,
\begin{align} \label{desiredconvergencepstozero} (\mu_\ep, \bnu_\ep) \to    (\mu, \bnu) \text{ in } C([0,1]; \P_1(\Rd)) \times  \M( \Rd \times [0,1];U) 
\end{align}
and $(\mu,\bnu)$ solves (\ref{mainproblem}).  
\end{proposition}
\begin{remark}[Convergence to unique minimizer] \label{uniquenessremark}
The fact that the preceeding convergence result  only holds ``up to a subsequence'' is due the generality of our assumptions, which  do not necessarily ensure that solutions of (\ref{mainproblem}) are unique. In particular, it is possible that different subsequences of $(\mu_\epsilon, \bnu_\epsilon)$   approximate different minimizers of (\ref{mainproblem}).

  On the other hand, if one imposes additional hypotheses to ensure that the minimizer $(\mu,\bnu)$ of (\ref{mainproblem}) is unique, as is true in many important examples (see section \ref{motivatingexamplessec}), then Proposition \ref{thm:convep} implies that every sequence $(\mu_\ep, \bnu_\ep)$ has a further subsequence that converges to $(\mu,\bnu)$. Therefore, $(\mu_\ep, \bnu_\ep)$ itself must converge to $(\mu,\bnu)$.
\end{remark}

While smaller values of $\epsilon$ lead to minimizers $\mu_\epsilon$ that are closer to the desired target measure $m_1$ at time $t=1$, it is   clear from the definition of $\F_\epsilon$ that minimizers are forced to satisfy $\mu_1 \in L^2$, which will always fail in the case of particle measures, $\mu_t = (1/N) \sum_{i=1} \delta_{y_i(t)}$.
We navigate this issue by incorporating an additional regularization into the terminal constraint. 

Consider a mollifier $k_\delta$ satisfying the following assumption:
\begin{assumption} \label{psias}
Suppose $k \in C_b(\Rd)$ is nonnegative, even, and $\int |x| k(x) dx <+\infty$. For any $\delta >0$, let $k_\delta(y) = k(y/\delta)/\delta^d$.
\end{assumption}
 For a probability measure $\nu$, the convolution of $\nu$ with $k_\delta$ is a bounded, continuous function, given by
\[ k_\delta*\nu(y) = \int_\Rd k_\delta (y-x) d \nu(x) . \]
In this way, we can define the regularized functional
\begin{align} \label{Fepsdeldef}
\F_{\epsilon,\delta}(\mu) : =   \frac{1}{\epsilon} \int_\Rd \left| k_\delta*\mu(y) - k_\delta* m_1(y) \right|^2 dy .
\end{align}
When $\int k = 1$ and $\mu \in L^2(\Rd)$, $\F_{\epsilon,\delta}(\mu)$   converges to $\F_\epsilon(\mu)$ as $\delta \to 0$. However, unlike $\F_\epsilon(\mu)$, $\F_{\epsilon,\delta}(\mu)$ is always finite --- in fact,  bounded by $2\|k\|_\infty/(\epsilon \delta^d)$ --- for all $\mu \in \P(\Rd)$. This leads to the optimization problem
\begin{align}
\tag{$\rm MFC_{\epsilon,\delta}$} \label{eq:main_optimsrrbb2}
&\min_{(\mu, \bnu) \in \C(m_0)} \E_{\epsilon,\delta}(\mu,\bnu) , \quad \E_{\epsilon,\delta}(\mu,\bnu):= \G(\mu,\bnu)+ \F_{\epsilon,\delta}(\mu_1) .   \\ \nonumber
\end{align}
We will often use that, expanding the square, using associativity of convolution, and abbreviating $K_\delta := k_\delta*k_\delta$, we may express $\F_{\epsilon,\delta}$ as follows:
\begin{align}
 \F_{\epsilon,\delta}(\mu) \label{altFepsdelta}  &= \frac{1}{\epsilon} \left[ \int (k_\delta*k_\delta* \mu) d \mu - 2 \int (k_\delta*k_\delta*m_1) d \mu+ \int  (k_\delta*k_\delta*m_1) d m_1\right]    \\
& = \frac{1}{ \epsilon} \left[ \int \left( (K_\delta*\mu)   - 2   (K_\delta* m_1)  \right) d \mu + C_{\delta,m_1} \right], \nonumber
\end{align}
where   $C_{\delta,m_1} : = \int (K_\delta* m_1) d m_1$ is a constant independent of $\mu$.

Under mild hypotheses, we show that, 
 as long as (\ref{eq:main_optimrrbb}) feasible,  solutions of  (\ref{eq:main_optimsrrbb2}) converge to a solution of (\ref{eq:main_optimrrbb}) as $\delta \to 0$, in the sense of the following theorem. (Note that, if (\ref{mainproblem}) is feasible, then both (\ref{eq:main_optimrrbb}) and (\ref{eq:main_optimsrrbb2}) are feasible for all $\epsilon, \delta >0$.)

\begin{theorem}[Convergence as $\delta \to 0$] \label{deltatozerothm}Suppose Assumptions \ref{asmp1} and \ref{psias} hold and $\int k = 1$. Fix $\epsilon >0$, and suppose that (\ref{eq:main_optimrrbb}) is feasible.
Then, for any  sequence  $(\mu_{\delta},\bnu_{\delta})$ of solutions to (\ref{eq:main_optimsrrbb2}), there exists $(\mu,\bnu)$ so that, up to a subsequence,
\begin{align} \label{desiredconvergencepstozeroagain}
(\mu_\delta, \bnu_\delta) \to (\mu, \bnu) \text{ in }  C([0,1]; \P_1(\Rd)) \times \M(\Rd \times [0,1];U) 
\end{align}
and $(\mu ,\bnu )$ solves (\ref{eq:main_optimrrbb}). 
\end{theorem}
As described in Remark \ref{uniquenessremark}, if one additionally assumes that the solution of (\ref{eq:main_optimrrbb}) is unique, then any sequence of minimizers of (\ref{eq:main_optimsrrbb2}) converges to this unique solution, without passing to a subsequence.

As indicated above, the key advantage of the regularized problem \eqref{eq:main_optimsrrbb2}, compared to \eqref{mainproblem}, is that the regularized problem admits a natural particle discretization.
%
Replacing $\mu_0$ with its empirical approximation, 
\[  \mu_0 \approx \frac{1}{N}\sum_{i=1}^N \delta_{y_{i,0}} , \qquad \by_0 :=  [ y_{i,0}]_{i=1}^N \in (\Rd)^N , \]
 we obtain the following finite dimensional, ODE constrained optimization problem: 
    \begin{align} \label{MFCN} \tag{$\rm MFC_{\epsilon,\delta,N}$}
&\min_{(\by, \bu) \in \mathcal{A}_N(\by_0)} \E_{\ep,\delta,N}(\by,\bu)   \\ \nonumber
\end{align}
    where
    \begin{align}
\E_{\ep,\delta,N}(\by,\bu) &:=  \frac{1}{N}\sum_{i=1}^N  \int_0^1 \psi(u_i(t))dt + \frac{1}{N} \sum_{i=1}^N \int_0^1  L_N(y_i(t),\mathbf{y}(t)) dt  \label{EepdelN} \\
    & \quad             + \frac{1}{\epsilon} \left[ \frac{1}{N^2}\sum_{i,j=1}^N  K_\delta(y_i(1)-y_j(1))- \frac{2}{N}\sum_{i=1}^N (K_\delta*m_1)(y_i(1)) + C_{\delta,m_1} \right]  \nonumber  \\
 \mathcal{A}_N(\by_0) &:= \left\{ (\by,\bu) : \by \in AC \left([0,1]; (\mathbb{R}^d)^N \right) \ , \ \ \bu \in L^1([0,1];U^N)  \right. \nonumber \\
& \qquad    \qquad  \ \    \left. \dot{y}_i(t) = \bF_N(y_i(t), \by(t)) + u_i(t)  \ , \  \  \forall i = 1, \dots, N  \ , \ \ \by(0) = \by_0 \right \}  . \nonumber
\end{align}
where $C_{\delta,m_1}$ is the constant from equation (\ref{altFepsdelta}), $U^N$ denotes the cartesian product of $U$, and the   system of ordinary differential equations holds in the Carath\'eodory sense.

In order to describe our hypotheses on $L_N$ and $\bF_N$, which follow Fornasier et. al., we begin by recalling the notions of symmetry and convergence introduced in their work \cite{fornasier2019mean}. A map $G_N: \Rd \times (\Rd)^N \to \R^k$ is \emph{symmetric} if
  \[ G_N(x,\by) = G_N(x, \sigma(\by)) \text{ for all permutations } \sigma: (\Rd)^N \to (\Rd)^N . \]
  Given a symmetric and continuous map $G_N$, we will consider the following notion of convergence to $G: \Rd \times \P_1(\Rd) \to \R^k$:
  \begin{align} \label{GNtoGconvergence} 
&\text{ for any sequence of natural numbers $N \to +\infty$ } \\
&\text{ and } \by_{N} := \left[y_{i,N} \right]_{i=1}^{N} \subseteq (\Rd)^{N} \text{ satisfying }  \mu_N := \frac{1}{N} \sum_{i=1}^{N} \delta_{y_{i,N}} \xrightarrow{N \to +\infty } \mu \text{ in } W_1, \nonumber\\ &\text{ we have } 
\lim_{N \to +\infty} \sup_{z \in C} \left| G_{N}(z, \by_{N}) - G(z,\mu) \right| = 0 , \quad \text{ for all } C \subset \subset \mathbb{R}^d . \nonumber
\end{align}

With this in hand, we turn to the hypotheses on $L_N$ and $\bF_N$:
  \begin{assumption} We assume $L_N$ and $F_N$ satisfy either of the following hypotheses: \label{MFCNas}\ 
    \begin{enumerate}[(i)]
\item
\begin{enumerate}[(a)]
\item  If $L$ satisfies Assumption \ref{asmp1}(\ref{Lcts}), we assume   \label{LNcts}
$L_N :\R^d \times (\Rd)^N \rightarrow [0,+\infty)$ is symmetric and continuous  and $L_N$ converges to $L$ as in equation (\ref{GNtoGconvergence}).
\item  If $L$ satisfies Assumption \ref{asmp1}(\ref{Llsc}),   we assume \label{LNlsc}   $L_N :\R^d \times (\Rd)^N \rightarrow [0,+\infty]$ is symmetric, lower semicontinuous, and  independent of the second variable, i.e. $L_N(x ,\by) = L_N(x, \bz)$ for all $x \in \Rd$, $\by, \bz \in (\Rd)^N$. Furthermore, we assume $L_N \nearrow L$ pointwise.
\end{enumerate}
\item $\bF_N: \Rd \times (\Rd)^N \rightarrow \R^d$ is continuous, symmetric, and there exist   $C_{F_N},C_{F_N}' \geq 0$ so that 
\[ |\bF_N(z,\by)| \leq C_{F_N} + C_{F_N}'  ( |z|  + \| \by \| ) , \quad \forall z \in \mathbb{R}^d , \ \by \in (\mathbb{R}^d)^N .  \]
Furthermore, we assume the compatibility condition
\[ \bF_N(z,\by) - \bF(z,\mu) \in U \text{ for all } z \in \Rd, \ \by \in (\Rd)^N , \ \mu \in \P_1(\Rd) . \]
Lastly, we assume $\bF_N$ converges to $\bF$ in as in equation (\ref{GNtoGconvergence}).
\end{enumerate}
\end{assumption}

We now consider existence of solutions to  (\ref{MFCN}). We show that if either (i) the control is unconstrained and the initial locations lie in the domain of the state and measure cost $L_N$ or if (ii) the state and measure cost $L_N$ is   continuous, then   (\ref{MFCN})   feasible. It then follows quickly that, whenever (\ref{MFCN}) is feasible, a solution exists. Our proof is a   mild adaptation of \cite[Proposition 4.2]{fornasier2019mean}, extending this result to the case when $L_N$ satisfies Assumption \ref{MFCNas}(\ref{LNlsc}).

\begin{proposition} \label{discretesolutionsexist}
Suppose Assumptions \ref{asmp1}, \ref{psias} and  \ref{MFCNas}    hold. Fix $\epsilon >0$, $\delta >0$, $N \in \mathbb{N}$ and $\by_0 \in (\mathbb{R}^d)^N$. Suppose at least one of the following structural assumptions holds:
\begin{enumerate}[(a)]
\item the control is unconstrained, $U =\Rd$ and  $y_{i,0} \in \{L_N <+\infty\}$ for all $i = 1, \dots, N$;
\item  the state and measure cost $L_N$ satisfies Assumption \ref{MFCNas}(\ref{LNcts}) .
\end{enumerate} 
 Then   (\ref{MFCN}) is feasible.
 
 Furthermore, whenever (\ref{MFCN}) is feasible, then a solution of (\ref{MFCN}) exists.
\end{proposition}

As a corollary of this and the preceding convergence results, we obtain sufficient conditions for existence of minimizers to the continuum optimization problems we consider.
\begin{corollary} \label{existenceofsolutions}
Suppose Assumptions \ref{asmp1}, \ref{psias}, and \ref{MFCNas}   hold and  
 $\int \phi(|x|)  d m_0 <+\infty$, where $\phi$ is as in Assumption \ref{asmp1prime} below.
 Furthermore, suppose that (\ref{mainproblem}) is feasible and that at least one of the following holds:
\begin{enumerate}[(a)]
\item the control is unconstrained, $U =\Rd$, and  $\supp m_0 \subseteq \{L<+\infty\}$;
\item  the state and measure costs $L$ and $L_N$ satisfy Assumptions  \ref{asmp1}(\ref{Lcts}) and  \ref{MFCNas}(\ref{LNcts}) .
\end{enumerate} 
Then, for all $\epsilon, \delta >0$, solutions of  (\ref{eq:main_optimsrrbb2}), (\ref{eq:main_optimrrbb}), and (\ref{mainproblem}) exist.
\end{corollary}

Finally, we consider the convergence of the discrete particle approximation (\ref{MFCN}) to the continuum. In particular, we show that, if  the problem either has an unconstrained control or is a controllable linear time invariant system,  then, for fixed $\epsilon, \delta >0$, solutions of (\ref{MFCN}) converge to a solution of (\ref{eq:main_optimsrrbb2}) as $N \to +\infty$. 
\begin{theorem}[Convergence of minimizers as $N \to +\infty$] \label{NtoinftyminimizersRdLTI} Suppose Assumptions \ref{asmp1}, \ref{psias}, and \ref{MFCNas}   hold, and $\phi$ is as in Assumption \ref{asmp1prime} below. Suppose there exists a convex set $\Omega$ so that 
\[ \supp m_0 \subset \subset \Omega \subseteq \{L<+\infty\} . \]
  Fix $\epsilon, \delta >0 $, and suppose (\ref{eq:main_optimsrrbb2}) is feasible.
Furthermore, assume at least one of the following structural assumptions holds: \begin{enumerate}[(a)]
\item Unconstrained control: $U = \Rd$ \label{URdcaseN}
\item Controllable linear time invariant system: \label{LTIcaseN}
\begin{enumerate}[(i)]
\item $L = L_N = 0$;
\item There exists $A \in M_{d\times d}(\Rd)$ so that  $\bF_N(y,\cdot)= \bF(y, \cdot) = Ay$;
 \item There exists a  full rank matrix $B \in M_{d \times k}(\Rd)$ with $R(B) = U$ for which  the system is  {controllable}: that is, the {controllability grammian},
$\Gamma(T) = \int_0^T e^{-\tau A} B B^T e^{-\tau A^T} d \tau $ is nonsingular for all $T>0$.
\item For any sequence $\eta  \to 0$, 
\[ \eta \phi( | \Gamma^{-1}(\eta)(I - e^{-\eta A}) z |) \to 0 \text{ uniformly for } z \in \supp m_0 .\]
 In particular, it is sufficient that $\supp m_0 \subseteq N(A)$. \label{chocolate}
\end{enumerate}
\end{enumerate}

Fix  $\by_{N,0} \in \Omega^N$  with $\sup_{i,N} | y_{i,N,0} | <+\infty$, satisfying
 \begin{align} \label{birdies}
 \frac{1}{N} \sum_{i=1}^N \delta_{y_{i,N,0}} \xrightarrow{N \to +\infty} m_0 \text{ in }\P_1(\Rd) .
 \end{align} 
Then, for any sequence $(\by_N, \bu_N) \in \mathcal{A}(\by_{N,0})$ of minimizers of (\ref{MFCN}),   there exists $(\mu, \bnu) \in \mathcal{C}(m_0)$ so that, up to a subsequence,
\begin{align} \label{conv1N}
&\frac{1}{N}\sum_{i=1}^N \delta_{y_{i,N}} \to \mu \text{ in }C([0,1]; \P_1(\Rd)) ,\\ 
&\frac{1}{N} \sum_{i=1}^N u_{i,N}(t) \delta_{y_{i,N}(t)} dt \to \bnu \text{ in } \M(\Rd \times [0,1]; U) . \label{conv2N}
\end{align}
Furthermore, any such limit point $(\mu, \bnu)$ is a minimizer of (\ref{eq:main_optimsrrbb2}). \end{theorem}

\begin{remark}[Discrete to continuum  in the absence of a terminal constraint] \label{extensionremark}
While the main focus on the present paper is mean field control in the presence of a terminal constraint,   we note that the previous theorem   includes the case of no terminal constraint, which arises when $k = k_\delta = 0$, so $\F_{\epsilon,\delta} = 0$. In this way,   the previous theorem extends the  gamma convergence result of Fornasier, Lisini, Orrieri, and Savar\'e \cite[Theorem 3.3]{fornasier2019mean} to the following two cases:
\begin{itemize}
\item  unbounded state and measure costs $L$ and $L_N$, as arise when modeling an obstacle,
\item  controllable linear time invariant systems, with constrained controls, $U \neq \R^d$.
\end{itemize}
\end{remark}

We now state our  main result, which shows that solutions of (\ref{MFCN}) converge to a solution of (\ref{mainproblem}) in the limit $\epsilon, \delta \to 0$ and $N \to +\infty$.  

\begin{theorem}[Convergence of Minimizers as $\epsilon, \delta \rightarrow 0$ and $N \rightarrow \infty$]  \label{maincoffeetheorem1} Suppose Assumptions \ref{asmp1}, \ref{psias}, and  \ref{MFCNas}    hold  and $\int k = 1$. Suppose (\ref{mainproblem}) is feasible and  there exists a convex set $\Omega$  so that
\[ \supp m_0 \subset \subset \Omega \subseteq \{ L <+\infty\} . \]
Furthermore, assume   either the unconstrained control hypothesis (a) or the controllable linear time invariant system   hypothesis (b) from  Theorem \ref{NtoinftyminimizersRdLTI}.  

Fix $\by_{N,0} \in \Omega^{N}$   with $\sup_{i,N} | y_{i,N,0} | <+\infty$ satisfying
\[ \frac{1}{N} \sum_{i=1}^{N} \delta_{y_{i,N,0}} \xrightarrow{N\to +\infty} m_0 \text{ in }\P_1(\Rd) . \]
Then, for any sequence
 $(\by_{\epsilon,\delta,N}, \bu_{\epsilon,\delta,N}) \in \mathcal{A}(\by_{N,0})$ of minimizers of (\ref{MFCN}),    there exists subsequences $\epsilon_m , \delta_m ,$ and $N_m $ so that, defining  
 \begin{align}  \label{cookies1}
 \by_m &:= \by_{\epsilon_m, \delta_m,N_m} , \qquad \bu_m = \bu_{\epsilon_m,\delta_m, N_m} ,  \\
\mu_m &:= \frac{1}{N}\sum_{i=1}^N \delta_{y_{i,m}} , \qquad \bnu_m := \frac{1}{N} \sum_{i=1}^N u_{i,m}(t) \delta_{y_{i,m}}  , \label{cookies2}
\end{align}
we have
\begin{align}
\mu_m \to \mu  \text{ in }C([0,1]; \P_1(\Rd)) \label{coffee1} , \\
 \bnu_m \to \bnu \text{ in } \M(\Rd \times [0,1]; U) ,  \label{coffee2}
\end{align}
where $(\mu,\bnu)$ is a minimizer of (\ref{mainproblem}).
\end{theorem}

Note that, in the preceding theorem, even if (\ref{mainproblem}) has a unique minimizer, our convergence result will only hold up to subsequences. This is due to the facts that the discretization parameter $N$ must grow sufficiently quickly with respect to $\epsilon$ and $\delta$ and  that $\delta$ must decay sufficiently quickly with respect to $\epsilon$. We leave the question of quantifying the relationship between $\epsilon, \delta,$ and $N$ for which convergence holds, without passing to a subsequence, to future work.

\subsection{Motivating Examples} \label{motivatingexamplessec}

	We now describe three important examples of mean field control problems of the form (\ref{star}), which motivate our numerical study.
	
A first special case is   the dynamic formulation of the p-Wasserstein distance on $\P_p(\Rd)$, the space of probability measures with finite $p$th moment, $M_p(\mu) = \int |x|^p d \mu(x)$, $p >1$,
	\begin{align}
	\label{eq:BBprb} \tag{$W_p^p$}
	&\min_{(\mu, \bu) \in {\mathcal{A}_{W}(m_0,m_1)}} \int_0^1 \int_{\Rd} |\bu(y,t)|^pd\mu_t(y) dt \\
	& \ \mathcal{A}_{W}(m_0,m_1) := \left\{ (\mu, \bu):  \mu \in AC([0,1], \mathcal{P}_1(\Rd)) , \   \bu \in L^1_{d \mu_t \otimes dt}( \Rd \times [0,1]; \Rd))   \right. \nonumber \\
&\left.      \ \ \quad     \qquad \qquad \qquad     \qquad \qquad \qquad \qquad    \partial_t \mu_t  + \nabla \cdot \left(  \bu_t \mu_t \right) = 0 \  , \ 
 \mu_0=m_0, \mu_1=m_1  \right\}  , \nonumber
\end{align}
which arises from (\ref{star}) by taking the choices $\psi = |\cdot|^p$, $L = 0$, $U= \Rd$, and $\bF = 0$, as introduced in the $p=2$ case by Benamou and Brenier \cite{BenamouBrenier} and generalized to $p>1$ by Ambrosio, Gigli, and Savar\'e \cite{ambrosiogiglisavare}.  To discretize this problem  via a particle approximation, (\ref{MFCN}), we take $L_N = 0$ and $\bF_N = 0$. When $m_0, m_1 \in \P_p(\Rd)$ and $m_1$ is absolutely continuous with respect to Lebesgue measure, there exists a unique minimizer of (\ref{eq:BBprb}), and the optimal value  of the objective function coincides with the $p$th power of the $p$-Wasserstein distance
\[ W_p^p(m_0, m_1) = \eqref{eq:BBprb}. \]

A second special case of the mean field control problem is the dynamic formulation of the $p$-Wasserstein distance around an  obstacle $\Omega$, represented by a open subset of $\Rd$ on which the interpolating measure $\mu_t$ is forbidden from placing mass, 
	\begin{align}
	\label{BBobstacle} \tag{$W_p^p(\Omega^c)$}
	&\min_{(\mu, \bu) \in \mathcal{A}_{\Omega^c}(m_0,m_1) } \int_0^1 \int_{\Rd} |\bu(y,t)|^pd\mu_t(y) dt \\
	 \ \mathcal{A}_{\Omega^c}(m_0,m_1) &:= \left\{ (\mu, \bu):  \mu \in AC([0,1], \mathcal{P}_1(\Rd))  ,  \bu \in L^1_{d \mu_t \otimes dt}( \Rd \times [0,1]; \Rd))   \right. \nonumber \\
&\left.      \    \ \  \ \           \partial_t \mu_t  + \nabla \cdot  (  \bu_t \mu_t  ) = 0   , 
 \mu_0=m_0, \mu_1=m_1 ,   \    \supp \mu_t \subseteq \Omega^c \text{ for a.e. } t \in [0,1] \right\} . \nonumber
\end{align}
This arises from (\ref{star}) under the same choices as the preceding example, with the exception that 
\begin{align*}
L(y,\mu) = \chi_{\Omega^c}(y) , \qquad \chi_{\Omega^c}(y) := \begin{cases} 0 &\text{ if } y \in \Omega^c ,\\ +\infty &\text{ if } y \in \Omega. \end{cases}
\end{align*}
To discretize this problem with the nonlocal terminal constraint, (\ref{MFCN}), we may take  $\bF_N = 0$  and $L_N$ to be any lower semicontinuous function $L_N(y, \mu_t) = L_N(y)$ that vanishes on $\Omega^c$ and converges up to $+\infty$ on $\Omega$. This variant of the optimal transport problem can also be interpreted, more geometrically, as optimal transport on the manifold with holes, $\Rd \setminus \Omega$. Unlike in the case without obstacles, minimizers in general are not unique, unless $\Omega^c$ is convex.
   
A third important example of the mean field control problem arises when an acceleration control is imposed on the evolution of the measure $\mu_t$,
\begin{align} \tag{$\text{MFC}_a$} \label{MFCa}
			& \min_{(\mu, \ba) \in \mathcal{A}_{a}(m_0,m_1)} \int_0^1\int_{\R^{2d}} |\ba(x,v,t)|^2d\mu_t(x,v)  dt \\
			&  \mathcal{A}_a(m_0,m_1) = \left\{ (\mu, \ba):  \mu \in AC([0,1], \mathcal{P}_1(\R^{2d})) ,   \ba \in L^1_{d \mu_t \otimes dt}( \R^{2d} \times [0,1];   \Rd))   \right. \nonumber \\
&\left.      \ \      \qquad \qquad \qquad   \partial_t \mu_t  + \nabla \cdot ( [v ~0]^t \mu_t) + \nabla \cdot  ( [0 ~\ba ]^t \mu_t  ) = 0   , 
 \mu_0=m_0, \mu_1=m_1  \right\}  , \nonumber
\end{align}
which arises from (\ref{star}) by   the choices $y = (x,v) \in \R^{2d}$, $\psi = |\cdot|^2$, $L = 0$, $\bu = (0,\ba) \in U= \{0 \} \times \Rd \subseteq \R^{2d}$,  and $\bF((x,v),\mu) = \bF((x,v)) = [v~0]^t $. To discretize this problem with the nonlocal terminal constraint, (\ref{MFCN}), we take $L_N = 0$ and $\bF_N((x,v),\bz)   = [v ~0 ]^t $. In this case, uniqueness of solutions is known \cite{hindawi2011mass,chen2016optimal}, when the probability measures are absolutely continuous and have compact support.

\subsection{Strategy of approach}\label{futurework}
The remainder of the paper is organized as follows. In section \ref{preliminariessection}, we specify the precise function spaces and notions of solution that we consider and prove the equivalence of the original formulation of the mean field control problem and the formulation in momentum coordinates. In section \ref{convergencesection}, we prove our main convergence results, relating minimizers of the optimization problems 
(\ref{MFCN}), (\ref{eq:main_optimsrrbb2}), (\ref{eq:main_optimrrbb}), and (\ref{mainproblem}) in the limits as $\epsilon, \delta \to 0$ and $N \to +\infty$. Finally, in section \ref{numericssec}, we implement our particle method numerically and use it to compute solutions of the dynamical formulation of the p-Wasserstein metric, optimal transport around an obstacle, and measure transport with acceleration constraints, as well as numerically analyzing the convergence properties of the method.

\section{Preliminaries} \label{preliminariessection}
\subsection{Notation} \label{notationsec}
Given $X \subseteq \Rd$, let  $\mathcal{M}(X)$ denote the space of finite signed Borel measures on $X$, let $\mathcal{P}(X)$ denote the space of Borel probability measures, and, for any $p \geq 1$, let $\mathcal{P}_p(X)$ denote  the subset of $\mathcal{P}(X)$ with finite $p$th moments, $M_p(\mu):=\int_X |x|^p d \mu(x) < +\infty$. We consider the space $\P_p(\Rd)$ to be endowed with the $p$-Wasserstein metric, $W_p$. We will often use the fact that convergence in the $W_p$ metric for any $p \geq 1$ implies narrow convergence, which is to say, convergence in the duality with bounded continuous functions $f:\Rd \to \R$.  For further details on the Wasserstein metrics and optimal transport,  we refer the reader to one of the many excellent textbooks on the subject \cite{ambrosiogiglisavare, villani2003topics, santambrogio2015optimal, figalli2021invitation,ambrosio2021lectures}.

Given $U$, a subspace of $\R^d$, let $\M(X;U)$ denote the space of  vector-valued finite signed Borel measures on $X$ with range in $U$. We will consider $\M(X;U)$ to be endowed with the narrow topology, which is to say, convergence in  $\M(X;U)$ is given by convergence in the duality with bounded continuous functions $f:X \to U$. For $X \subseteq \Rd$, $Y \subseteq \R^d$, and a Borel map $\bt :X\rightarrow Y$, define $\bt_{\#} : \mathcal{P}(X) \rightarrow \mathcal{P}(Y)$  by $(\bt_{\#}\mu)(B) = \mu(\bt^{-1}(B))$ for all Borel subsets $B \subseteq Y$.  For any Borel set $B$, we let $1_B$ denote the indicator function on $B$, that is, $1_B(x) = 1$ if $x \in B$ and $1_B(x) = 0$ if $x \not \in B$.

For any product space $\Pi_{i=1}^N X_i $, $X_i \subseteq \Rd$, let $\pi^i$ denote the projection onto the $i$th coordinate. When elements of the product space are denoted by a distinguished variable, e.g. $(x,y,z) \in X \times Y \times Z$, we   write $\pi^x, \pi^y, \pi^z$ for   projections onto the corresponding components.

 The set of   continuous curves in $\P_1(\Rd)$ will be denoted by $C([0,1];\P_1(\Rd))$  and the set of absolutely continuous curves will be denoted by $AC([0,1]; \P_1(\Rd))$. We will often identify curves $\mu  \in C([0,1];\P_1(\Rd))$ with the element $\tilde{\mu} \in \P_1(  \mathbb{R}^d \times [0,1])$ satisfying
\begin{align} \label{curvesmeasuresrem}
\int_{ \Rd \times [0,1]} f(y,t) d\tilde{\mu}(y,t) = \int_0^1 \int_\Rd f(y,t)d\mu_t(y)dt 
\end{align}
for all $f \in $ $C_b(  \Rd \times [0,1])$. We will abbreviate $\tilde{\mu}$ by $ d \mu_t \otimes dt$.
 Under the hypotheses of Assumption \ref{asmp1}, $\mu \in AC([0,1]; \P_1(\Rd))$ and $ \bu \in L^1_{d\mu_t \otimes dt} ( \Rd \times [0,1]; U))  $ is a distributional solution of the equation $\partial_t \mu_t  + \nabla \cdot \left( (\bF(\cdot, \mu_t) \mu_t  + \bu(\cdot, t) \mu_t \right) = 0$ if, for all $\varphi \in C^\infty_c(\Rd \times [0,1])$,
 \begin{align} \label{firstdistributionalsolution}
\int_0^1 \int_\Rd \partial_t \varphi(y,t) + \nabla \varphi(y,t) \cdot \left( (\bF(y, \mu_t)   + \bu(y,t)   \right) d \mu_t(y) \ dt = 0 .
\end{align}

Given $\bnu \in \M(  \R^d \times [0,1];U)$ and $ {\mu} \in C([0,1];\P_1(\Rd))$, if $\bnu$ is absolutely continuous with respect to the measure $ \tilde{\mu} = d \mu_t \otimes dt$, we let $\bu :  \Rd \times[0,1] \rightarrow U$ denote the Radon Nikodym derivative: $\bu(y,t) \ d \tilde{\mu}(y,t) = \bu(y,t) d \mu_t(y) dt= d\bnu(y,t)  $. Furthermore, for any  $\bnu \in \M(  \R^d \times [0,1];U)$,  we let $\bnu_t$ denote the disintegration of $\bnu$ with respect to $\pi^t_\# \bnu \in \M([0,1])$, where $\pi^t$ denotes the projection onto the second, temporal component of $  \Rd \times [0,1]$. In this way, for all $f \in C_b(  \Rd \times [0,1])$,
\[ \int_{\Rd \times [0,1]} f(y,t) d \bnu(y,t) = \int_0^1 \left( \int_\Rd f(y,t) d \bnu_t(y) \right) d (\pi^t_\# \bnu)(t). \] 
We say that $ (\mu,\bnu) \in AC([0,1]; \P_1(\Rd)) \times \mathcal{M}(\Rd \times [0,1] ; U)$ is a distributional solution of $ \partial_t \mu_t  + \nabla \cdot \left( (\bF(\cdot, \mu_t) \mu_t  + \bnu_t \right) = 0 $ if, for all $\varphi \in C^\infty_c(\Rd \times [0,1])$,
 \begin{align} \label{seconddistributionalsolution}
\int_0^1   \int_\Rd \left( \partial_t \varphi(y,t) + \nabla \varphi(y,t) \cdot   \bF(y, \mu_t)  \right)d \mu_t(y) dt    + \int_0^1 \int_\Rd \nabla \varphi(y,t) \cdot d \bnu  (y,t) = 0 .
\end{align}

\subsection{Control cost and momentum coordinates}
We now state the general hypothesis on our control cost $\psi$, following previous work by Fornasier, et. al. \cite[Section 2.2]{fornasier2019mean}. 

\begin{assumption} \label{asmp1prime}
Assume  the control cost $\psi: U \to [0,+\infty)$ is convex, lower semicontinuous, and satisfies $\psi(0) = 0$. Furthermore, assume there exists a \emph{moderating function} $\phi: [0,+\infty) \to [0, +\infty)$   that is strictly convex, continuously differentiable, satisfies $\phi(0) = 0$ and $\phi'(0)=0$, is superlinear at $+\infty$, and for which there exists $K >0$ so that 
\begin{align} \label{doubledouble} \phi(2r) \leq K (1+\phi(r)) \quad \text{ for all } r \in [0, +\infty) .
\end{align}
Finally, assume  that the control cost and the moderating function are related as follows: there exists $C >0$ so that 
\begin{align} \label{phipsirelation} \phi(|x|) -1 \leq \psi(x) \leq C(1+\phi(|x|)) \ , \quad \forall x \in U .\end{align}
\end{assumption}

Examples of control costs satisfying this assumption include
\begin{itemize}
\item $\psi(x)= \frac{1}{p}|x|^p$ for $p>1$;
\item $\psi(x) = \begin{cases} \frac{1}{p} |x| &\text{ for } |x|\leq 1 \\ \frac{1}{p}|x|^p &\text{ for } |x|>1 \end{cases}$ \quad for $p>1$.
\end{itemize}
In particular, the role of the moderating function $\phi$ is that it allows one to decouple the   smoothness and monotonicity properties  from the local behavior of $\psi$. In the following lemma, we gather two elementary observations about the relationship between $\psi$ and $\phi$ that we will use in what follows.

\begin{lemma} \label{fancydoubling}
Suppose $\phi$ and $\psi$  satisfy Assumption \ref{asmp1prime}.
\begin{enumerate}[(i)]
\item There exists $R_\phi>0$ so that  \label{psisuperlinearity} $|x| \leq \phi(|x|) + R_\phi \leq \psi(x) + R_\phi + 1 \ , \ \forall x \in U$.
\item There exists $K' >0$ so that, for all $D \geq 1$, $
\phi(Dr) + Dr \leq K' D^{K'} (\phi(r)+r) , $  $\forall r \in [0,+\infty).$\end{enumerate}
\end{lemma}

We now make precise the sense in which (\ref{star}) and its the momentum formulation (\ref{mainproblem}) are equivalent.  See Appendix \ref{basicpropertiesappendix} for the proof.
\begin{lemma} \label{momequiv}
If $(\mu, \bu)$ solves  (\ref{star}), then $(\mu, \bu \ d \mu_t \otimes dt)$ solves (\ref{mainproblem}). On the other hand, if $(\mu, \bnu)$ solves (\ref{mainproblem}), then $d \bnu \ll d\mu_t \otimes dt $ and $(\mu, d \bnu/(d\mu_t \otimes dt))$ solves (\ref{star}).
\end{lemma}

\section{Convergence of minimizers as $\epsilon \to 0$, $\delta \to 0$, and $\N \to +\infty$}  \label{convergencesection}
We now turn to the proofs of our major results, relating minimizers of the problems (\ref{MFCN}), (\ref{eq:main_optimsrrbb2}), (\ref{eq:main_optimrrbb}), and (\ref{mainproblem}) in the limits as $\epsilon \to 0$, $\delta \to 0$, and $N \to +\infty$. We begin, in section \ref{lcscpt}, by establishing some fundamental lower semicontinuity and compactness properties of the continuum mean field control problem. In section \ref{epstozerosec}, we consider the $\epsilon \to 0$ limit. 
In section \ref{deltatozerosec}, we consider the behavior as $\delta \to 0$, focusing on the proof of $\Gamma$-convergence of the objective functionals, after which convergence of minimizers follows quickly. Finally, in section  \ref{discretesec}, we consider the $N \to +\infty$ limit. Our arguments in   section \ref{discretesec} are strongly inspired by previous work by Fornasier, et. al. \cite{fornasier2019mean}. Furthermore, as explained in Remark \ref{extensionremark},  we succeed in extending these arguments to our more general hypotheses on the state and measure cost $L$ and to controllable linear time invariant systems, with $U \neq \Rd$.

\subsection{Lower semicontinuity and compactness} \label{lcscpt}
In this section, we collect some fundamental lower semicontinuity and compactness properties for the continuum mean field control problem. See   appendix \ref{basicpropertiesappendix} for the proofs.

First, we note that, under the hypotheses of Assumption \ref{asmp1}, the functional $\mathcal{G}$
is lower semicontinuous. 
\begin{lemma}
\label{lem:epcont}
Suppose Assumption \ref{asmp1} holds. Then the functional $\mathcal{G}$,  defined in  (\ref{mainproblem}),
is lower semicontinuous on $C([0,1],\mathcal{P}_1(\mathbb{R}^d)) \times \mathcal{M}( \mathbb{R}^d \times [0,1]; U) $. 
\end{lemma}

Next, we observe that sublevels of $\G$ in the constraint set $\C(m_0)$ are sequentially compact.

\begin{lemma} \label{compactnesseps}
Suppose Assumption \ref{asmp1} holds. If $(\mu_n, \bnu_n) \in \C(m_0)$ and $\sup_{n  } \G(\mu_n,\bnu_n) <+\infty$, then there exists $(\mu,\bnu) \in \C(m_0)$ so that, up to a subsequence,
\begin{align} \label{convergencecompactness} (\mu_n, \bnu_n) \text{ converges to } (\mu, \bnu) \text{ in }  C([0,1]; \P_1(\Rd)) \times \M(\Rd \times[0,1];U)
\end{align}
 and $\G(\mu,\bnu) \leq \sup_{n  } \G(\mu_n,\bnu_n) $.
\end{lemma}

\subsection{Convergence as $\epsilon \to 0$: soft to hard terminal constraint}  \label{epstozerosec} 
We now collect our results on  the convergence of (\ref{eq:main_optimrrbb}) to (\ref{mainproblem}) as $\epsilon \to 0$. The proofs follow standard $\Gamma$-convergence arguments, which we defer to Appendix \ref{epstozeroappendix}.

We first observe that the sequence of functionals $\E_\epsilon$, defined in (\ref{eq:main_optimrrbb}), $\Gamma-$converges to the functional $\E$, defined in   (\ref{mainproblem}), as $\epsilon \rightarrow 0$. 

\begin{proposition}[{$\Gamma$-convergence as $\epsilon \to 0$}] \label{propgammaeptozero} Suppose Assumption \ref{asmp1} holds.
\begin{enumerate}[(i)]
\item \label{firstpartep} For every sequence $(\mu_\ep,\bnu_\ep)$ converging in $C([0,1];\P_1(\Rd)) \times \M( \Rd \times [0,1]; U)$ to $(\mu,\bnu)$, we have $\E(\mu,\bnu) \leq \liminf_{\ep\rightarrow 0}\E_{\epsilon}(\mu_\epsilon,\bnu_\epsilon) .$
\item For every $(\mu,\bnu) \in C([0,1];\P_1(\Rd)) \times \M( \Rd \times [0,1] ; U)$, 
$\E(\mu,\bnu) \geq \limsup_{\ep\rightarrow 0}\E_{\epsilon}(\mu,\bnu) .$ \label{secondpartep}
\end{enumerate}
\label{prop:epgammacon}
\end{proposition}

It is then    an immediate consequence of this result that   minimizers of (\ref{eq:main_optimrrbb}) converge to a minimizer of (\ref{mainproblem}), up to a subsequence, as stated in Proposition \ref{thm:convep}; see appendix \ref{epstozeroappendix}.

\subsection{Convergence as $\delta \to 0$: nonlocal to local penalization on terminal measure} \label{deltatozerosec}
In the present section, we show that, for fixed $\epsilon >0$, the sequence of functionals $\E_{\epsilon,\delta}$, defined in (\ref{eq:main_optimsrrbb2}), $\Gamma-$converges to the functional $\E_\epsilon$, defined in  equation \ref{eq:main_optimrrbb}, as $\delta \rightarrow 0$. 

\begin{proposition}[$\Gamma$-convergence as $\delta \to 0$] 
Suppose Assumptions \ref{asmp1} and \ref{psias} hold and $\int k =1$.\label{gammadelta}
Fix $\epsilon>0$.
\begin{enumerate}[(i)]
\item For every sequence $(\mu_{\delta},\bnu_{\delta})$ converging in $ C([0,1];\P_1(\Rd)) \times \M(\Rd \times [0,1] ; U)$ to $(\mu, \bnu)$, we have $\E_\ep(\mu,\bnu) \leq \liminf_{\delta \rightarrow 0}\E_{\ep,\delta}(\mu_{\delta},\bnu_{\delta}) . $ \label{firstpartdeltatozero}
\item For all $(\mu,\bnu) \in C([0,1];\P_1(\Rd)) \times \M(\Rd \times [0,1] ; U)$,
$\E_\ep(\mu,\bnu) \geq \limsup_{\delta\rightarrow 0}\E_{\ep,\delta}(\mu,\bnu) . $ \label{secondpartdeltatozero}
\end{enumerate}
\end{proposition}

\begin{proof}
First, we consider part (\ref{firstpartdeltatozero}). Without loss of generality, we may pass to a subsequence so   that  
\begin{align} \lim_{\delta \rightarrow 0}\E_{\epsilon,\delta}(\mu_\delta,\bnu_\delta)   = \liminf_{\delta \rightarrow 0}\E_{\epsilon,\delta}(\mu_\delta,\bnu_\delta)  < +\infty.
\label{Eepdelbound}
\end{align}
In particular, this ensures that
\[ \sup_\delta \ep^{-1} \|k_\delta*\mu_{\delta,1} - k_\delta*m_1\|_{L^2(\Rd)}^2  = \sup_{\delta} \F_{\ep,\delta}(\mu_{\delta,1}) <+\infty . \]

Since $m_1 \in L^2(\Rd)$, $k_\delta*m_1 \to m_1$ in $L^2(\Rd)$, and in particular, $k_\delta*m_1$  is uniformly bounded in $L^2(\Rd)$. Thus $k_\delta*\mu_{\delta,1}$ must be uniformly bounded in $L^2(\Rd)$, and, up to another subsequence, $k_\delta*\mu_{\delta,1}$ converges weakly in $L^2(\Rd)$.
By the convergence of $\mu_\delta $   to $\mu$, we have that $\mu_{\delta,1}$ converges to $\mu_1$  in $\P_1(\Rd)$. Furthermore,  \cite[Lemma 2.3]{carrillo2019blob} ensures $k_\delta*\mu_{\delta,1} \to \mu_1$ in distribution. Thus,  by uniqueness of limits, $k_\delta*\mu_{\delta,1}$ converges to $\mu_1$ weakly in $L^2(\Rd)$ and, in particular, $\mu_1 \in L^2(\Rd)$.

Now we consider the limits of the functionals along these sequences. Since $\mu_{\delta,1}$ converges to $\mu_1$  in $\P_1(\Rd)$, by  \cite[Theorem 4.1]{carrillo2019blob},
\begin{align} \label{firstpartfepdel} \liminf_{\delta \to 0} \int K_\delta*\mu_{\delta,1} d \mu_{\delta,1} \geq  \int \mu_1^2 . 
\end{align}
Since $k_\delta*\mu_{\delta,1} \to \mu_1$ weakly in $L^2(\Rd)$ and $k_\delta*m_1 \to m_1$ strongly in $L^2(\Rd)$,  for  $C_{\delta,m_1}$ as in equation (\ref{altFepsdelta}),  we obtain
\begin{align} \label{secondpartfepdel}
&   \liminf_{\delta\to 0} - 2 \int    (K_\delta* m_1)    d \mu_{\delta,1} + C_{\delta,m_1}  \\
&\quad =   \liminf_{\delta\to 0}   - 2 \int  (k_\delta* \mu_{\delta,1}) (k_\delta*m_1)  +  \int (K_\delta* m_1) d m_1   \nonumber    \geq   -2 \int m_1 \mu_1 +  \int m_1^2 . \nonumber
   \end{align}
   Thus, combining equation (\ref{firstpartfepdel}) and (\ref{secondpartfepdel}) with the expression for $\F_{\ep,\delta}$ in equation (\ref{altFepsdelta}), we obtain
\[\F_{\ep}(\mu_1)  = \frac{1}{\ep} \| \mu_1 - m_1\|_{L^2(\Rd)}^2 \leq \liminf_{\delta \rightarrow 0} \F_{\ep,\delta}(\mu_{\delta,1})\]

Due to the lower semicontinuity of the functional $\mathcal{G}$, shown in Lemma \ref{lem:epcont}, we can therefore conclude that 
\[\E_{\epsilon}(\mu,\bnu) \leq \liminf_{\delta \rightarrow 0}\E_{\epsilon,\delta}(\mu_\delta,\bnu_\delta). \]

Now, we turn to part (\ref{secondpartdeltatozero}). We may assume that $\E_\ep(\mu,\bnu)<\infty$, otherwise the inequality is trivial. Thus, $\F_\ep(\mu_{1})<+\infty$, so $\mu_{1} \in L^2(\Rd)$ and $k_\delta*\mu_1 \to \mu_1$ in $L^2(\Rd)$. Thus, $\lim_{\delta \rightarrow 0} \F_{\ep,\delta}(\mu_1) = \F_{\ep}(\mu_1)  .$
Therefore, we   conclude that, 
\[\E_{\epsilon}(\mu,\bnu) \geq \limsup_{\delta \rightarrow 0}\E_{\epsilon,\delta}(\mu^\delta,\bnu^\delta). \]
\end{proof}

Theorem \ref{deltatozerothm}, on the convergence of minimizers of (\ref{eq:main_optimsrrbb2}) to a minimizer of (\ref{eq:main_optimrrbb})  now follows from  this $\Gamma$-convergence result via a standard argument; see Appendix \ref{epstozeroappendix}.

\subsection{Convergence as $N \to +\infty$: discrete to continuum} \label{discretesec}
We now turn to the convergence of the spatially discrete problem, in the continuum limit $N \to +\infty$. 
Note that, by definition of $\E_{\ep,\delta,N}$ in equation (\ref{EepdelN}), for any $(\by,\bu) \in (\Rd)^N \times U^N$,
\begin{align*} \E_{\ep,\delta,N}(\by,\bu) &=  \G_N(\by,\bu) +   \F_{\ep,\delta}\left( \frac{1}{N} \sum_{i=1}^N \delta_{y_{i}(1)} \right)  \\
    \G_N(\by,\bu)&:=   \frac{1}{N}\sum_{i=1}^N  \int_0^1 \psi(u_i(t))dt + \frac{1}{N} \sum_{i=1}^N \int_0^1  L_N(y_i(t),\mathbf{y}(t)) dt .
    \end{align*}
    We begin with the following lemma.

\begin{lemma} \ 
 \label{lemma}  Suppose Assumptions \ref{asmp1} and \ref{psias} hold, and 
 suppose  $\by_N  \in  (\Rd)^N  $ satisfies $ \frac{1}{N}\sum_{i=1}^N \delta_{y_i} \to \mu$ narrowly for some $\mu \in  \P(\Rd)$. Then, for all $\epsilon, \delta>0$,
\[  \lim_{N \to +\infty} \F_{\ep,\delta}\left( \frac{1}{N} \sum_{i=1}^N \delta_{y_i} \right) = \F_{\ep,\delta}(\mu). \]
\end{lemma}

\begin{proof}
Let $\mu_N = \frac{1}{N}\sum_{i=1}^N \delta_{y_i}$.
 Note that, using the expression (\ref{altFepsdelta}) and neglecting the constant $C_{\delta,m_1}$,  we may rewrite the left hand side of the expression in the following way:
\begin{align*}
 \int \left[ K_\delta * \mu_N - 2(K_\delta* m_1) \right] d \mu_N =  \iint  K_\delta (x-y) d \mu_N(x) d \mu_N(y) - 2 \int (K_\delta* m_1)  d \mu_N
\end{align*}
By Assumption \ref{psias}, $(x,y) \mapsto K_\delta(x-y)$ and $x \mapsto k_\delta*m_1(x)$ are bounded and continuous, so the result is an immediate consequence of the definition of narrow convergence.
\end{proof}

Now, combining the preceding lemma with the $\Gamma$-convergence result of Fornasier et. al.\cite[Theorem 3.2]{fornasier2019mean}, making appropriate adaptations when the state and measure costs satisfy our alternative assumptions   \ref{asmp1}(\ref{Llsc}) and   \ref{MFCNas}(\ref{LNlsc}), we  obtain the  following  result.
\begin{proposition}[{$\Gamma$-convergence as $N \to +\infty$}] \label{gammathm} 
Suppose Assumptions \ref{asmp1}, \ref{psias}, and \ref{MFCNas} hold.
Fix $\epsilon>0$ and $\delta >0$.
\begin{enumerate}[(i)]
\item  \label{fullliminf}  Consider $(\by_N, \bu_N) \in C([0,1]; (\Rd)^N \times L^1([0,1]; U^N)$, and suppose there exists  $(\mu, \bnu) $ for which
\begin{align}
&\frac{1}{N}\sum_{i=1}^N \delta_{y_{i,N}} \to \mu \text{ in }C([0,1]; \P_1(\Rd) , \\
 & \frac{1}{N} \sum_{i=1}^N u_{i,N}(t) \delta_{y_{i,N}(t)} dt \to \bnu \text{ in } \M(  \Rd \times [0,1] ; U) . \end{align}
Then, $\liminf_{N \to +\infty} \E_{\ep,\delta,N}(\by_N, \bu_N) \geq \E_{\epsilon, \delta}(\mu ,\bnu)$. 
\item \label{fulllimsup} Suppose $\int \phi(|x|) d m_0(x) < +\infty$. For every $(\mu,\bnu) \in C(m_0)$ such that 
  $\G(\mu,\bnu)<+\infty $ and $N \in \mathbb{N}$, 
 there exists  $\by_{N,0} \in (\Rd)^N$ and $  (\by_N, \bu_N) \in  \mathcal{A}_N(\by_{N,0})$ with $y_{i,N,t} \subseteq \supp \mu_t$ for every $i = 1, \dots, N$  and $t \in [0,1]$, so that
\begin{align} 
&\frac{1}{N}\sum_{i=1}^N \delta_{y_{i,N}} \to \mu \text{ in }C([0,1]; \P_1(\Rd) , \label{firstN} \\
 & \frac{1}{N} \sum_{i=1}^N u_{i,N}(t) \delta_{y_{i,N}(t)} dt \to \bnu \text{ in } \M(  \Rd \times [0,1] ; U) ,\label{secondN} \\
    &\limsup_{N \to +\infty} \E_{\ep,\delta,N}(\by_N, \bu_N) \leq \E_{\epsilon,\delta}(\mu, \bnu) .
\label{thirdN}
\end{align}
\end{enumerate}
\end{proposition}

\begin{proof}

First, we show part (\ref{fullliminf}). Up to passing to a subsequence, we may assume without loss of generality that
\begin{align} \label{ENfinitenesswlog} \lim_{N \to +\infty} \E_{\ep,\delta,N}(\by_N, \bu_N) = \liminf_{N \to +\infty} \E_{\ep,\delta,N}(\by_N, \bu_N) <+\infty .
\end{align}
By Lemma \ref{lemma}, 
\begin{align} \label{Fstuff} \F_{\ep,\delta}\left( \frac{1}{N} \sum_{i=1}^N \delta_{y_{i,N}(1)} \right)  \to \F_{\ep,\delta}(\mu_1) .
\end{align}

Suppose that $L$ and $L_N$ satisfy   Assumption \ref{asmp1}(\ref{Lcts}) and Assumption \ref{MFCNas}(\ref{LNcts}). By \cite[Theorem 3.2(i)]{fornasier2019mean}, we have
 \[ \liminf_{N \to +\infty} \G_N(\by_N,\bu_N) \geq \G(\mu,\bnu) . \]
 Since $\E_{\ep, \delta}(\mu,\bnu) = \G(\mu,\bnu) +   \F_{\ep,\delta}(\mu_1)$, combining the preceding inequality with (\ref{Fstuff}) gives the result.
 
 On the other hand, suppose that  $L$ and $L_N$ satisfy   Assumption \ref{asmp1}(\ref{Llsc}) and Assumption \ref{MFCNas}(\ref{LNlsc}).  By \cite[Theorem 3.2(i)]{fornasier2019mean}, we have 
\begin{align} \label{justPsi} \liminf_{N \to +\infty}   \frac{1}{N}\sum_{i=1}^N  \int_0^1 \psi(u_i(t))dt \geq \Psi(\bnu |\mu) .
\end{align}
Furthermore, for $\mu_N := \frac{1}{N} \sum_{i=1}^N \delta_{y_{i,N}}$, Fatou's lemma ensures
 \begin{align*}
 \liminf_{N \to +\infty}  \frac{1}{N} \sum_{i=1}^N \int_0^1  L_N(y_{i,N}(t)) dt  &= \liminf_{N \to +\infty} \int_0^1 \int_\Rd L_N d \mu_{N,t} dt \\
& \geq \int_0^1  \liminf_{N \to +\infty} \int_\Rd L_N d \mu_{N,t} dt .
 \end{align*}
 For any $M \in \mathbb{N}$, the fact that $L_N \nearrow L$ and $L_M$ is lower semicontinuous and nonnegative ensures
 \begin{align*}
 \liminf_{N \to +\infty} \int_\Rd L_N d \mu_{N,t} \geq \liminf_{N \to +\infty} \int_\Rd L_M d \mu_{N,t} \geq \int_\Rd L_M d \mu_{t} .
 \end{align*}
 Sending $M \to +\infty$ on the right hand side, the monotone convergence theorem ensures
 \begin{align} \label{LNliminf}
  \liminf_{N \to +\infty} \int_\Rd L_N d \mu_{N,t} \geq \int_\Rd L d \mu_{t}  . 
  \end{align}
 Combining this with (\ref{justPsi}), we obtain  $ \liminf_{N \to +\infty} \G_N(\by_N, u_N) \geq \G(\mu,\bnu) .$
Finally, combining this with (\ref{Fstuff}) gives the result.

We now show part (\ref{fulllimsup}). 
It is an immediate consequence of \cite[Theorem 3.2(i)]{fornasier2019mean} that (\ref{firstN}) and (\ref{secondN}) hold and 
\begin{align} \label{justPsi2} \limsup_{N \to +\infty}   \frac{1}{N}\sum_{i=1}^N  \int_0^1 \psi(u_i(t))dt \leq \Psi(\bnu |\mu) .
\end{align}
(The fact that we may choose $y_{i,N}(t) \in \supp(\mu_t)$ for every $i = 1, \dots, N$ and $t \in [0,1]$ can be seen by inspection of the proof: in the equation following \cite[equation (6.21)]{fornasier2019mean}, we may assume $\gamma_{i,k,m} \in \supp \tilde{\pi}^k \subseteq \supp \pi$.)
As before, by Lemma \ref{lemma}, 
\begin{align} \label{puppies}
\F_{\ep,\delta}\left( \frac{1}{N} \sum_{i=1}^N \delta_{y_{i,N}(1)} \right)  \to \F_{\ep,\delta}(\mu_1) .
\end{align} 

If $L$ and $L_N$ satisfy   Assumption \ref{asmp1}(\ref{Lcts}) and Assumption \ref{MFCNas}(\ref{LNcts}), \cite[Theorem 3.2(i)]{fornasier2019mean} also gives $\limsup_{N \to +\infty} \G_N(\by_N, \bu_N) \leq \G(\mu,\bnu). $
Since $\E_{\ep, \delta}(\mu,\bnu) = \G(\mu,\bnu) + \F_{\ep,\delta}(\mu_1)$, combining the preceding inequality with (\ref{puppies})  gives the result in this case.

On the other hand, suppose $L$ and $L_N$ satisfy  Assumption \ref{asmp1}(\ref{Llsc}) and Assumption \ref{MFCNas}(\ref{LNlsc}). Without loss of generality, we may assume $\E(\mu,\bnu)<+\infty$, so that 
\begin{align*}
\int_0^1 \int_\Rd L d \mu_t dt < +\infty &\implies \int_\Rd L d \mu_t < +\infty \text{ for a.e. } t \in [0,1] \\
&\implies \{ y_{i,N}(t)  \}_{i=1}^N \subseteq \supp \mu_t \subseteq \{L<+\infty\} \text{ for a.e. } t \in [0,1] .
\end{align*}
Thus, the fact that $L_N \nearrow L$ implies
\begin{align*}
&\limsup_{N\to +\infty} \int_0^1 \int_\Rd L_N d \mu_{N,t} dt \leq \limsup_{N\to +\infty}  \int_0^1 \int_\Rd L d \mu_{N,t} dt \\
&\quad =  \limsup_{N\to +\infty} \int_0^1 \int_\Rd L 1_{\{L<+\infty\}} d \mu_{N,t} dt  \leq  \int_0^1 \int_\Rd L 1_{\{L<+\infty\}} d \mu_t dt \leq  \int_0^1 \int_\Rd L d \mu_t dt,
\end{align*}
where the second to last inequality follows since $L$ is continuous on the closed set $\{L<+\infty\}$, so  $L 1_{\{L<+\infty\}}$ is upper semicontinuous. Combining this with    (\ref{justPsi2})  and (\ref{puppies}) gives the result.   \\  

\end{proof}

The preceding convergence result can now be used to show  that, for fixed $\epsilon, \delta >0$, as $N \to +\infty$, minimizers of the spatially discrete problem (\ref{MFCN}) converge to a solution of (\ref{eq:main_optimsrrbb2}), up to a subsequence. The result follows from a standard argument, which we defer to Appendix \ref{epstozeroappendix}.

\begin{theorem}[Convergence of minimizers as $N \to +\infty$] \label{Ntoinftyminimizers} Suppose Assumptions \ref{asmp1}, \ref{psias}, and \ref{MFCNas}   hold, and suppose  $\int \phi(|x|)  d m_0(x) ~<~+\infty$. Fix $\epsilon, \delta >0 $, and suppose (\ref{eq:main_optimsrrbb2}) is feasible.
 
 Then, there exists  $\by_{N,0} \in (\Rd)^N$ satisfying 
 \begin{align}
 \label{whynot1}
 y_{i,N,0} \in \supp m_0 \text{ for all } i =1, \dots, N, \\
  \label{whynot} 
 \frac{1}{N} \sum_{i=1}^N \delta_{y_{i,N,0}} \to m_0 \text{ in }\P_1(\Rd) ,
 \end{align} so that, if $(\by_N, \bu_N) \in \mathcal{A}(\by_{N,0})$ is a minimizer of (\ref{MFCN}), then 
 \begin{align} \label{goodinitialdata}
\lim_{N \to +\infty} \E_{\ep,\delta,N}(\by_N, \bu_N)  =\inf_{(\mu, \bnu) \in \C(m_0)} \E_{\epsilon,\delta}(\mu,\bnu).\end{align}
and  there exists $(\mu, \bnu) \in \mathcal{C}(m_0)$ so that, up to a subsequence, 
\begin{align} \label{conv1}
&\frac{1}{N}\sum_{i=1}^N \delta_{y_{i,N}} \to \mu \text{ in }C([0,1]; \P_1(\Rd)), \\ 
&\frac{1}{N} \sum_{i=1}^N u_{i,N}(t) \delta_{y_{i,N}(t)} dt \to \bnu \text{ in } \M(\Rd \times [0,1]; U) . \label{conv2}
\end{align}
Furthermore, any such limit point $(\mu, \bnu)$ is a minimizer of (\ref{eq:main_optimsrrbb2}).  
\end{theorem}
As described in Remark \ref{uniquenessremark}, if one additionally knows that the solution of (\ref{eq:main_optimsrrbb2}) is unique, then any sequence of minimizers $(\by_N, \bu_N)$  converges to the unique solution, without passing to a subsequence. However, note that while the above theorem ensures \emph{there exists} a sequence of initial conditions $\by_{N,0}$ for which minimizers of (\ref{MFCN}) converge to a minimizer of (\ref{eq:main_optimsrrbb2}), at this level of generality, it is not known if the result is true for \emph{all} choices of initial conditions.

In the special cases of unconstrained controls and controllable linear time invariant systems,  Theorem \ref{NtoinftyminimizersRdLTI} ensures that this result indeed continues to hold for all well-prepared initial conditions $\by_{N,0}$.  Our argument is strongly inspired by  the proof of   \cite[Theorem 3.3(iii)]{fornasier2019mean}, which we adapt to our more general hypotheses on the state and measure cost $L$ and the setting of controllable linear time invariant systems.

\begin{proof}[Proof of Theorem \ref{NtoinftyminimizersRdLTI}]

 Let $ {\by}_{N,0} \in \Omega^N$ be an arbitrary  sequence s.t.  $\sup_{i,N} | {\by}_{i,N,0} |<+\infty$ and (\ref{birdies}) holds. Suppose   we can show that, for any sequence of minimizers $( {\by}_N,  {\bu}_N) \in \mathcal{A}_N( {\by}_{N,0})$ of (\ref{MFCN}) we have
\begin{align} \label{geckos} \limsup_{N \to +\infty} \E_{\epsilon,\delta,N}( {\by}_N,  {\bu}_N) \leq  e_{m_0} <+\infty ,\end{align}
with $e_{m_0}$ as in equation (\ref{em0def}).
Then, up to a subsequence, we must have
\[ \sup_{N} \frac{1}{N}\sum_{i=1}^N  \int_0^1 \psi( {u}_{i,N}(t))dt  <+\infty . \]
Thus \cite[Theorem 3.1]{fornasier2019mean} ensures there exists $(\mu, \bnu) \in \mathcal{C}(m_0)$ so that, up to a subsequence, (\ref{conv1N})-(\ref{conv2N}) hold. Furthermore, for any such limit point $(\mu,\bnu)$,  Proposition \ref{gammathm}(\ref{fullliminf}) ensures
\begin{align} \label{whynotineq3}
e_{m_0} \leq    \E_{\epsilon, \delta}(\mu ,\bnu) \leq \liminf_{N \to +\infty} \E_{\ep,\delta,N}( {\by}_N,  {\bu}_N) \leq e_{m_0} .
\end{align}
Thus, $(\mu,\bnu)$ is a minimizer of (\ref{eq:main_optimsrrbb2}).

In this way, it suffices to show (\ref{geckos}). Let $\tilde{\by}_{N,0}$ be the sequence of initial conditions from Theorem \ref{Ntoinftyminimizers}. Since $m_0$ is compactly supported, we have $\sup_{i,N} |\tilde{y}_{i,N,0} |<+\infty$. Let
\begin{align*}
R = \max \{ \sup_{i,N} |y_{i,N,0} | , \sup_{i,N}|\tilde{y}_{i,N,0} | \} .
\end{align*}
Fix $(\tilde{\by}_N, \tilde{\bu}_N) \in 
\A(\tilde{\by}_{N,0})$ that minimize (\ref{MFCN}). Theorem \ref{Ntoinftyminimizers} ensures that
\begin{align} \label{bookmark}
\lim_{N \to +\infty} \E_{\ep,\delta,N}(\tilde{\by}_N, \tilde{\bu}_N)  =e_{m_0} ,
\end{align}
and there exists $(\tilde{\mu},\tilde{\nu}) \in \mathcal{C}(m_0)$ so that, up to a subsequence, (\ref{conv1}-\ref{conv2}) hold.
For the remainder of the proof, we work with this subsequence of $N$.

Define 
\begin{align} \label{zetadef}
\zeta(r) := \phi(r) + r .
\end{align}
By Lemma \ref{fancydoubling}, $\zeta(r) \leq 2\phi(r) +R_\phi$. Thus, by  \cite[Lemma 2.5]{fornasier2019mean},  
\[ \lim_{n \to +\infty} C_\zeta \left( \frac{1}{N} \sum_{i=1}^N \delta_{y_{i,N,0}}
, \frac{1}{N} \sum_{i=1}^N \delta_{\tilde{y}_{i,N,0}} \right) = 0 ,\]
where $C_\zeta(\mu,\nu)$ denotes the optimal value of the optimal transport problem between $\mu$ and $\nu$ with cost matrix $\zeta(x-y)$.  By Birkhoff's theorem, up to a permutation of the indices $i \mapsto y_{i,N}$, we may assume that
\[ c_N := C_\zeta \left( \frac{1}{N} \sum_{i=1}^N \delta_{y_{i,N,0}} , \frac{1}{N} \sum_{i=1}^N \delta_{\tilde{y}_{i,N,0}} \right) = \frac{1}{N} \sum_{i=1}^N \zeta(|y_{i,N,0} - \tilde{y}_{i,N,0}|) \text{ for all }N. \]

We now seek to make a small modification to $(\tilde{\by}_N, \tilde{\bu}_N)$ so that the initial condition agrees with $\by_{N,0}$ and, as $N \to +\infty$, the value of the discrete energy along the modified sequence converges to $e_{m_0}$. Our modification takes different forms, depending on the structural assumptions we consider. Fix $\eta \in (0,1)$, and recall that $\phi$ is increasing, convex, and satisfies the doubling condition (\ref{doubledouble}). Note that, for any $r,C \geq 0$,
\begin{align*}
\phi(r+C) \quad \leq \quad \begin{cases} \phi(2r) &\text{ if } r \geq C \\ \phi(2C) &\text{ if } r \leq C \end{cases}  \quad \leq  \quad \phi(2r) + \phi(2C).
\end{align*}

In case (a), in which $U = \Rd$,  define
\begin{align*}
y^{(a)}_{i,N,\eta}(t) &= \begin{cases} (1-t/\eta) y_{i,N,0} + (t/\eta) \tilde{y}_{i,N,0} &\text{ if } t \in [0,\eta) , \\
\tilde{y}_{i,N}(t-\eta) &\text{ if } t \in [\eta, 1] , \end{cases} \\
u^{(a)}_{i,N,\eta}(t) &= \begin{cases} (y_{i,N,0} - \tilde{y}_{i,N,0})/\eta - \bF_N(y^{(a)}_{i,N,\eta}(t), \by^{(a)}_{N,\eta}(t)) &\text{ if } t \in [0,\eta) , \\
\tilde{u}_{i,N}(t-\eta) &\text{ if } t \in [\eta, 1] .\end{cases}
\end{align*}
By construction, we have $(\by^{(a)}_{N,\eta}, \bu^{(a)}_{N,\eta} ) \in A(\by_{N,0})$. Furthermore, there exists $R>0$ so that,  for all $t \in [0,\eta)$,  $\eta \in (0,1)$,
\begin{align} \label{dontrunaway}
\sup_{i,N} |y^{(a)}_{i,N,\eta}(t)|   \leq R <+\infty . 
\end{align}
Thus, the convergence of $\bF_N(x, \by_{N,\eta})$ to $\bF(x, m_0)$ as $N \to +\infty$ locally uniformly in $x$ and $t$, implies  there exists $C_F >0$ so that
\[ \sup_{t \in [0,\eta),i,N  } |\bF_N(y^{(a)}_{i,N,\eta}(t), \by^{(a)}_{N,\eta}(t)  | < C_F. \]

Now, we estimate the state and measure cost in case (a) for $t\in[0,\eta)$. If $L$ and $L_N$ satisfy Assumptions \ref{asmp1}(\ref{Lcts}) and \ref{MFCNas} (\ref{LNcts}), then
 the convergence of $L_N(x, \by_{N,\eta})$ to $L(x, m_0)$ as $N \to +\infty$ locally uniformly in $x$ and $t$, and the estimate (\ref{dontrunaway}) implies 
\[ \sup_{t \in [0,\eta),i, N} | L_N(y^{(a)}_{i,N,\eta}(t), \by^{(a)}_{N,\eta}(t)  | <+\infty. \]
On the other hand,  if $L$ and $L_N$ satisfy Assumptions \ref{asmp1}(\ref{Llsc}) and \ref{MFCNas}(\ref{LNlsc}),
\begin{align*} 
\sup_{ t \in [0,\eta), i, N}   L_N(y^{(a)}_{i,N,\eta}(t), \by^{(a)}_{N,\eta}(t)) \leq \sup_{ t \in [0,\eta],i,N }   L(y^{(a)}_{i,N,\eta}(t))  <+\infty,
\end{align*}
where we use that $\{y^{(a)}_{i,N,\eta}(t)\}_{ t \in [0,\eta), i, N}$  is a bounded subset of $\Omega     \subseteq \{ L<+\infty\}$, so the image of this set under $L$ is also bounded.  
Thus, under either assumption on $L$, there exists $C_L >0$ so that
\begin{align*}
\frac{1}{N} \sum_{i=1}^N \int_0^\eta  L_N (y^{(a)}_{i,N,\eta}(t), \by^{(a)}_{N,\eta}(t)) dt \leq \eta C_L .
\end{align*}

 Finally, we estimate the control cost in case (a), estimating first in terms of the admissible function $\phi$. For $t \in [0, \eta)$,
\begin{align*}
 \phi\left ( \left|u^{(a)}_{i,N,\eta}(t)) \right| \right)  &= \phi \left( \left| (y^{(a)}_{i,N,0} - \tilde{y}^{(a)}_{i,N,0})/\eta - \bF_N(y^{(a)}_{i,N,\eta}(t), \by^{(a)}_{N,\eta}(t)) \right|  \right)\\
 & \leq \phi \left( \left| (y_{i,N,0}^{(a)} - \tilde{y}^{(a)}_{i,N,0})/\eta \right|  + C_F \right)  \\
  & \leq \phi \left( 2\left| (y^{(a)}_{i,N,0} - \tilde{y}^{(a)}_{i,N,0})/\eta \right| \right) + \phi(2C_F)
 \end{align*}

Now, we turn to  case (b). In this case,  define
\begin{align*}
y^{(b)}_{i,N,\eta}(t) &= \begin{cases} e^{tA}  y_{i,N,0} + \int_0^t e^{-\tau A}  {u}^{(b)}_{i,N,\eta}(\tau) d \tau  &\ \ \quad \qquad \text{ if } t \in [0,\eta) , \\
\tilde{y}_{i,N}(t-\eta)  &\ \ \quad \qquad \text{ if } t \in [\eta, 1] , \end{cases} \\
u^{(b)}_{i,N,\eta}(t) &= \begin{cases} -B B^T e^{-tA^T} \Gamma^{-1}(\eta) [y_{i,N,0} - e^{-\eta A} \tilde{y}_{i,N,0} ] &\text{ if } t \in [0,\eta) , \\
\tilde{u}_{i,N}(t-\eta) &\text{ if } t \in [\eta, 1] .\end{cases}
\end{align*}
Again, by construction, we have $(\by^{(b)}_{N,\eta}, \bu^{(b)}_{N,\eta} ) \in A(\by_{N,0})$. (See, for example, \cite[Theorem 5.2]{bacciotti2019stability}.)
Likewise, in case (b), for $t \in [0,\eta)$,
\begin{align*}
& \phi\left ( \left|u^{(b)}_{i,N,\eta}(t)) \right| \right)  \\
 &  = \phi \left( \left|  -B B^T e^{-tA^T} \Gamma^{-1}(\eta) [y_{i,N,0} - e^{-\eta A} \tilde{y}_{i,N,0} ] \right|  \right)\\
 & \leq \phi \left( \left|  B B^T e^{-tA^T} \Gamma^{-1}(\eta) [y_{i,N,0} - \tilde{y}_{i,N,0} ] \right| +  \left|  B B^T  e^{-tA^T} \Gamma^{-1}(\eta) [I - e^{-\eta A} ] \tilde{y}_{i,N,0} ] \right|  \right) \\
  & \leq \phi \left(2 \left|  B B^T e^{-tA^T} \Gamma^{-1}(\eta) [y_{i,N,0} - \tilde{y}_{i,N,0} ] \right| \right) + \phi \left( 2  \left|  B B^T e^{-tA^T} \Gamma^{-1}(\eta) [I - e^{-\eta A} ] \tilde{y}_{i,N,0} ] \right|  \right) \\
    &  \leq    \phi \left(C_{A,B}  \|\Gamma^{-1}(\eta)\| \left|  y_{i,N,0} - \tilde{y}_{i,N,0}   \right|  \right) + \phi \left( C_{A,B} | \Gamma^{-1}(\eta)(I-e^{-\eta A}) \tilde{y}_{i,N,0} |  \right) 
 \end{align*}
 for some $C_{A,B}>0$, depending on $A$ and $B$. 
 
   Since $\|\Gamma^{-1}(\eta)\|$ is a continuous function of $\eta >0$ that diverges to $+\infty$ as $\eta \to 0$, there exists a decreasing function $\sigma:(0,1) \to \R$ so that   $\sigma(\eta) \geq \frac{1}{\eta}$ and $\sigma(\eta) \geq \|\Gamma^{-1}(\eta)\|$ for all $\eta \in (0,1)$. Define
 \begin{align*}
 \omega(\eta, \tilde{y}_{i,N,0}) &= \begin{cases} 0 &\text{ in case (a),} \\
 C_{A,B} | \Gamma^{-1}(\eta)(I-e^{-\eta A}) \tilde{y}_{i,N,0} |  &\text{ in case (b).} \end{cases}\end{align*}
 
Therefore, considering both case (a) and case (b) simultaneously, there exists $C' \geq1 $ so that, for $K'>0$ as in Lemma \ref{fancydoubling} and $\zeta(|x|) = \phi(|x|) +|x|$,
\begin{align*}
  \phi \left ( \left|u_{i,N,\eta}(t)) \right| \right) &\leq \phi(C' \sigma(\eta) \left|  y_{i,N,0} - \tilde{y}_{i,N,0}   \right| ) + \phi(\omega(\eta,\tilde{y}_{i,N,0})) +  C' \\
  & \leq K' (C'\sigma(\eta))^{K'}   \zeta(  \left|  y_{i,N,0} - \tilde{y}_{i,N,0}   \right| )    + \phi(\omega(\eta,\tilde{y}_{i,N,0})) + C' 
  \end{align*}
Combining  this estimate with inequality (\ref{phipsirelation}),  there exists $C'' >0$ depending on $K'$, $C'$, and $C$ so that 
\begin{align*}
&\frac{1}{N} \sum_{i=1}^N \int_0^\eta \psi( u_{i,N,\eta}(t) ) dt \\
 &\quad \leq \eta C+ \frac{C}{N} \sum_{i=1}^N \int_0^\eta \phi( |u_{i,N,\eta}(t)| ) dt \\
&\quad \leq  \eta C''+ C'' \eta (\sigma(\eta))^{K'}   \left( \frac{1 }{N} \sum_{i=1}^N   \zeta(| (y_{i,N,0} - \tilde{y}_{i,N,0}|)  \right) + \eta C'' \frac{1}{N} \sum_{i=1}^N \phi(\omega(\eta,\tilde{y}_{i,N,0}))  \\
&\quad =   C'' \eta (\sigma(\eta))^{K'} c_N + o(1)  ,\text{ as } \eta \to 0,
\end{align*}
where we use the fact that, in case (a), $\omega \equiv 0$, and in case (b), hypothesis (\ref{chocolate}) ensures that 
\begin{align*} \eta \phi (\omega (\eta, \tilde{y}_{i, N,0} )) &\leq \eta K' (C_{A,B})^{K'} \zeta(| \Gamma^{-1}(\eta)(I-e^{-\eta A}) \tilde{y}_{i,N,0} | ) \\
& \leq \eta K' (C_{A,B})^{K'}  (2 \phi(| \Gamma^{-1}(\eta)(I-e^{-\eta A}) \tilde{y}_{i,N,0} | )+R_\phi) \to 0  ,
\end{align*}
uniformly in $i, N$.

Finally note that, by definition of $y_{i,N,\eta}$ in both case (a) and case (b), 
\begin{align*}
&\F_{\ep,\delta}\left( \frac{1}{N} \sum_{i=1}^N \delta_{{y_{i,N,\eta}}(1)} \right)  - \F_{\ep,\delta}\left( \frac{1}{N} \sum_{i=1}^N \delta_{\tilde{y}_{i,N}(1)} \right) \\
&\quad = \F_{\ep,\delta}\left( \frac{1}{N} \sum_{i=1}^N \delta_{\tilde{y}_{i,N}}(1-\eta) \right)  - \F_{\ep,\delta}\left( \frac{1}{N} \sum_{i=1}^N \delta_{\tilde{y}_{i,N}(1)} \right) .
\end{align*}
By Lemma \ref{lemma}, this  vanishes as $N \to +\infty$ and $\eta \to 0$, since the empirical measures in both arguments converge to $\tilde{\mu}_1$.

Combining the above estimates, we have shown that, in both case (a) and case (b),  
\begin{align*}
&  \E_{\epsilon,\delta,N}(\by_{N,\eta},\bu_{N,\eta}) \\
&\quad \leq \frac{1}{N} \sum_{i=1}^N \int_0^\eta \psi( \bu_{N,\eta}(t) ) dt + \frac{1}{N}\sum_{i=1}^N  \int_0^\eta L_N (y_{i,N,\eta}(t), \by_{N,\eta}(t)) dt  + \E_{\epsilon,\delta, N}(\tilde{\by}_{N},\tilde{\bu}_{N}) \\
&\qquad \quad +   \F_{\ep,\delta}\left( \frac{1}{N} \sum_{i=1}^N \delta_{{y_{i,N,\eta}}(1)} \right)  -   \F_{\ep,\delta}\left( \frac{1}{N} \sum_{i=1}^N \delta_{\tilde{y}_{i,N}(1)} \right)  \\
&\quad \leq  C'' \eta (\sigma(\eta))^{K'} c_N  + \E_{\epsilon,\delta, N}(\tilde{\by}_{N},\tilde{\bu}_{N}) + o(1) , \text{ as } \eta \to 0 ,  \ N \to +\infty. 
\end{align*}
Taking  $\eta = \sigma^{-1}(c_N^{-1/K'}) $ implies $\sigma(\eta) = c_N^{-1/K'}$ and $\eta (\sigma(\eta))^{K'} c_N = \eta \xrightarrow{N \to +\infty} 0$.
Thus,
\begin{align*}
\limsup_{N \to +\infty} \E_{\epsilon,\delta,N}(\by_{N,\eta},\bu_{N,\eta})  \leq \limsup_{N \to +\infty} \E_{\epsilon,\delta,N}(\tilde{\by}_{N},\tilde{\bu}_{N}) = e_{m_0}
\end{align*}
Finally, since $(\by_{N,\eta},\bu_{N,\eta}) \in \mathcal{A}(\by_{N,0})$ for all $\eta \in (0,1)$, for any sequence of minimizers $(\by_{N},\bu_N)  \in \mathcal{A}(\by_{N,0})$ of (\ref{MFCN}), we have
\[ \limsup_{N \to +\infty} \E_{\epsilon,\delta,N}(\by_{N},\bu_{N})  \leq\limsup_{N \to +\infty} \E_{\epsilon,\delta,N}(\by_{N,\eta},\bu_{N,\eta}) \leq e_{m_0} .\]
This shows (\ref{geckos}), which completes the proof. 
\end{proof}

%
%
%

Our main theorem, Theorem \ref{maincoffeetheorem1}, is now an immediate consequence of Theorem \ref{thm:convep}, Theorem \ref{deltatozerothm}, and Theorem \ref{NtoinftyminimizersRdLTI}, via a standard diagonal argument. See appendix \ref{epstozeroappendix}.

\section{Numerics} \label{numericssec}

\subsection{Numerical Implementation} \label{numericalimplementation}

We now apply the particle approximation developed in the previous sections to compute three key examples of the mean field optimal control problem (\ref{star}): the dynamic formulation of the 2-Wasserstein distance, 2-Wasserstein optimal transport around obstacles, and measure transport with acceleration constraints; see section \ref{motivatingexamplessec} for more details. For simplicity of exposition, we describe our approach in the context of the first two 2-Wasserstein examples. See section \ref{accsection} for the case of acceleration controls.

For classical $p=2$ optimal transport, with and without obstacles, the particle discretization of the mean field optimal control problem (\ref{MFCN}) becomes
\begin{align} \label{numericalobjective0}
		\inf_{y_i,u_i}& \left\{ \frac{1}{N}\sum_{i=1}^N\int_0^1|u_i(t)|^2dt  +   \frac{1}{N} \sum_{i=1}^N \int_0^1 L_N(y_i(t)) dt   \right. \\
		& \qquad   \qquad \qquad  \left. + \frac{1}{\epsilon}\left[ \frac{1}{N^2}\sum_{i,k=1}^N  K_\delta(y_i(1)-y_k(1))- \frac{2}{N}\sum_{i=1}^N (K_\delta*m_1)(y_i(1)) \right] \right\}  \nonumber
			\end{align}
			subject to the differential equations
			\begin{align}
			\dot{y_i}(t) = u_i(t) , \quad
			y_i(0)=y_{i,0}. \label{icconstraint}
\end{align}
(Note that we neglect the constant $C_{\delta,m_1}$ in the objective functional, since it does not affect minimizers.)
In order to numerically approximate solutions of  this constrained optimization problem, we must discretize time, incorporate the ODE constrains, and develop a method for approximating the minimizer.

 First, we discretize time $t \in [0,1]$ on a uniform grid with $M$ grid points and time step $h = 1/(M-1)$, approximating $y_i(t)$ by the vector
\begin{align}
    [ {y}_{i,j}]_{j=1}^{M}  &\approx \left[y_{i,0}, y_i(h),  y_i(2h), \ldots,  y_i(1-h), y_i(1) \right], \label{eq:xtilde2}
\end{align}
which, by definition, incorporates initial condition  constraint  (\ref{icconstraint}). Next, we approximate the velocity by a first order finite difference,
\begin{align}
u_i(t) \approx \frac{  y_i(t+h)-y_i(t) }{h} .
\end{align}
Substituting into the objective function (\ref{numericalobjective0}), we arrive at a fully discrete problem, which is to minimize the sum of the kinetic energy, the potential energy, and the nonlocal energy,
 \begin{align} \label{numericalobjective}
& \min_{ { y}_{i,j} } \left\{ {\rm KE}( {y}_{i,j})+{\rm PE}(y_{i,j}) + {\rm NE}( {y}_{i,j}) : y_{i,1} = y_{i,0} \ \forall i = 1, \dots, N \right\} , \end{align}
for
\begin{align} 
{\rm KE}({y}_{i,j})  &:= \frac{M-1}{ N}\sum_{i=1}^N\left|\texttt{diff}( { y}_{i,j})\right|^2   \label{numericalKE} \\
{\rm PE}( {y}_{i,j})  &:= \frac{1}{N(M-1)} \sum_{i=1}^N \sum_{j=1}^M L_N(y_{i,j})  \label{numericalPE} \\
{\rm NE}( {y}_{i,j}) &:=  \frac{1}{\epsilon} \left[\frac{1}{N^2}\sum_{i,k=1}^N  K_\delta \left( {y}_{i,M}- {y}_{k,M}\right) - \frac{2}{N}\sum_{i=1}^N (K_\delta*m_1)( {y}_{i,M})  \right] , \label{numericalNE}
\end{align}
where $i = 1, \dots, N$ and $ j = 2, \dots, M$, with
\begin{equation*}
    \texttt{diff}( { y}_{i,j}) := \left[ {y}_{i,2}-y_{i,1}, {y}_{i,3}- {y}_{i,2}, \ldots,  {y}_{i,M} - {y}_{i,M-1} ,0 \right].
\end{equation*}
The resulting optimization problem is unconstrained. 
  In what follows, we commit a mild abuse of notation and let $y_i(t)$ denote the linear interpolation of $\{ y_{i,j}\}_{j=1}^M$ in time. In the majority of the simulations that follow, we consider the case without obstacles $L_N \equiv 0$; section \ref{obstaclesec} considers the optimal trajectories in the presence of obstacles.
  
  In the present simulations, we choose $k_\delta, K_\delta: \Rd \to \R$ to be Gaussian mollifiers,
\begin{equation}
    K_\delta(\theta) = (2\pi \delta^2)^{-d/2}\exp{(-\|\theta\|^2/2\delta^2)}, \quad k_\delta = K_{\delta/\sqrt{2}} .\label{eq:gaussian_mollifier}
\end{equation}

\begin{remark}[target measures $m_1 \not \in   L^2(\Rd)$] \label{notL2}
While our theoretical convergence results require $m_1 \in L^2(\Rd)$ (see Assumption \ref{asmp1}), our numerical approach extends naturally to  empirical target measures, $m_1 = \frac{1}{N}\sum_{i=1}^N\delta_{w_i}$. 
In this case,   NE becomes
\begin{align} \label{numericalNEempirical}
 {\rm NE}( {y}_{i,j})    =   \frac{1}{\epsilon N^2}\sum_{i,k=1}^N     \left[K_\delta \left( {y}_{i,M}- {y}_{k,M}\right) -  2  K_\delta( {y}_{i,M} - w_k)    \right]  .
\end{align}

\end{remark}

Finally, once we have arrived at the fully discrete minimization problem (\ref{numericalobjective}), we compute an approximate minimizer by gradient descent. With this approach, updating the trajectory of the $i$th particle only requires local information at the terminal time regarding the proximity to other particles and the value of the target $m_1$. Derivatives are calculated automatically by PyTorch \cite{paszke2019pytorch}, with a prescribed learning rate $\alpha$ and maximum number of descent steps $n$. In addition, we use standard learning rate reduction and early stopping mechanisms.
If the objective function value does not decrease for two steps, the learning rate $\alpha$ is reduced to $0.2\alpha$, as long as the reduction is larger than $10^{-8}$, as implemented in PyTorch. If the objective function value does not decrease for 5 steps, the algorithm is terminated early.

As we illustrate in our simulations (see section \ref{losslandscapesection}), we do not expect the loss landscape to be convex, so we do not have guarantees that our gradient descent approach will converge to a global minimizer of (\ref{numericalobjective}). However, as our numerical results attest, even our simple gradient descent approach leads to reasonable results. We leave the important question of developing more accurate methods for computing an approximate minimizer to future work, as the main focus of the present paper is to analyze the effects of discretizing dynamic optimal transport by a regularized particle method.

Due to the nonconvexity of the loss landscape, we expect that the results of gradient descent will depend strongly on the initialization. In the following simulations, we initialize the trajectories to be straight lines terminating at the center of mass of the target distribution. In other words, we choose  $ {y}_{i,j}$ so that 
\begin{align} \label{typicalinitialization}
& {y}_{i,M} = \int y m_1(y) dy \text{ for all } i = 1, \dots, N , \\
&[ {y}_{i,j}]_{j=2}^{M-1} \text{ linearly interpolates between }  y_{i,0} \text{ and } {y}_{i,M} \text{ for all }i = 1, \dots, N. \label{typicalinitialization2}
\end{align}

Python code for all  experiments is available  at \url{https://github.com/HarlinLee/BlobOT}.

\subsection{Parameter selection} \label{parameterselection}
In practice, we expect the gradient descent dynamics for computing an approximate minimizer of (\ref{numericalobjective}) to depend on the choice of the parameters $N, \delta, \epsilon, M, \alpha$ and $n$. For the continuum formulation of the problem, there is no downside to choosing $\epsilon$ to be arbitrarily small, independent of $N$ and $\delta$, since it simply enforces the terminal constraint more exactly. For this reason, in our simulations, we typically choose $\epsilon$ to be a fixed, very small parameter. 

On the other hand, we cannot hope for good results by choosing $\delta$ arbitrarily small without regard to $N$. When $\delta$ is very small, the ``sensing radius'' of the mollifier $K_\delta$ becomes very small, preventing particles from detecting the correct target distribution, unless they were coincidentally initialized very close to target. For example, when the target distribution is an empirical measure, particles would not be able to detect desired target locations $w_i$ too far from their current terminal location $y_{i,M}$ because the small radius of concentration of $K_\delta$. In a similar way, if one chooses $N$ arbitrarily large without regard to $\delta$, assuming that the empirical measure $\frac{1}{N} \sum_{i=1}^N \delta_{y_{i,M}}$ narrowly converges to a limiting probability measure $\mu_1$ the nonlocal energy approximates
\begin{align*}
 {\rm NE}(y_{i,j})   \approx F_{\epsilon,\delta}(\mu_1) =  \frac{1}{\epsilon} \| k_\delta*\mu_1 - k_\delta*m_1\|_{L^2}^2.
 \end{align*}
 It is a classical result   that there exist $\mu_1$, $m_1$ arbitrarily far apart for which the $L^2$ differences of their regularizations is arbitrarily small, so in this way, when $\delta >0$ is fixed too large with respect to $N$,  ${\rm NE}(y_{i,j})$ is not as accurate for imposing the terminal constraint.

For these reasons, we allow $\delta \to 0$ and $N \to +\infty$ simultaneously. Inspired by quantitative error estimates available for  classical vortex blob methods for the Euler and Navier-Stokes equations \cite{BealeMajda1,BealeMajda2} and blob methods for the aggregation equation \cite{CraigBertozzi}, we let $\delta$ scale with $N$,
\begin{align} \label{deltaNscaling}
\delta = N^{-k/d} , \quad \text{ for } 0 \ll k < 1 .
\end{align}
In the present manuscript, we take $k = 0.99$.
As a consequence, two particles distance $O(N^{-1/d})$ apart in $d$-dimensions will be able to ``sense'' one another via the mollifier $K_\delta$. In particular, when $K_\delta$ is a Gaussian mollifier, equation (\ref{deltaNscaling}) ensures its standard deviation is slightly larger than the interparticle distance on a regular grid.

With regard to the number of time steps $M$, since we expect good regularity in time, we anticipate higher accuracy when $M$ is large, at the expense of increasing the  dimension of the optimization problem (\ref{numericalobjective}). 
Due to the fact that the optimal trajectories of the dynamic 2-Wasserstein problem will always follow straight lines, in the case   $L_N \equiv 0$, it suffices to take $M=2$. However, we allow for flexibility in choice of $M$, to accommodate more general formulations of the problem, including obstacles and acceleration constraints. In the numerical experiments that follow, we choose $M>2$, even though this is an unnecessarily computational expense in the $L_N \equiv 0$ case, in order to illustrate that the kinetic energy ${\rm KE}$ on its own is able to effectively straighten the trajectories. 
 
Finally, we consider the choice of parameters  $\alpha$ and $n$ for our gradient descent dynamics.  As shown in section \ref{losslandscapesection}, when $\epsilon$ is large, the energy landscape flattens, but we can still achieve good results with a larger learning rate $\alpha$. For this reason, we typically take $\alpha$ to scale with $\epsilon$. As a consequence,  when $\epsilon$ is small, the learning rate $\alpha$ is small, and it may take more steps $n$ to reach an approximate global minimizer. In this way, while from the perspective of the continuum problem, there is no downside to taking $\epsilon$ arbitrarily small, from the perspective of the fully discrete problem, we do see that when $\epsilon$ is very small, it can require more iterations of gradient descent to reach an approximate global minimizer; see section \ref{erroralonggraddescent}.

\subsection{Comparison with Python Optimal Transport}
\begin{figure}[h!]
    \centering
    \includegraphics[width=.82\linewidth]{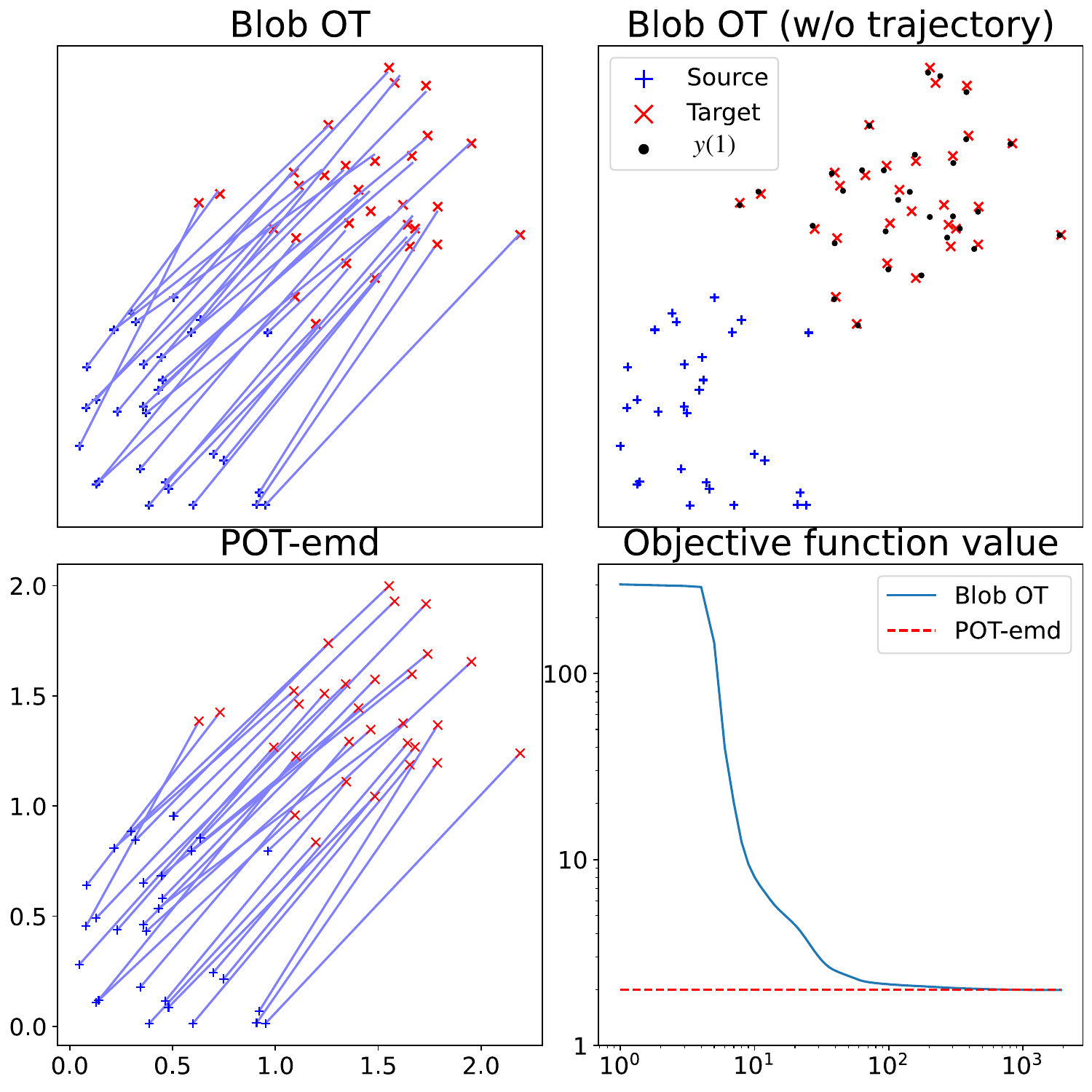}
    \caption{Comparison between our particle method for computing the dynamic formulation of the 2-Wasserstein distance and classical approaches to computing optimal transport, e.g., via the Kantorovich formulation of the problem, as implemented in Python Optimal Transport emd (POT-emd).} 
    \label{fig:uniform_to_uniform}
\end{figure}

In Figure \ref{fig:uniform_to_uniform}, we begin by comparing the results of Blob OT to classical approaches for solving the 2-Wasserstein optimal transport problem. In particular, we compare our method to the Python Optimal Transport emd (POT-emd) function, which computes optimal transport between two empirical measures via solving the Kantorovich formulation of the problem using the Hungarian algorithm \cite{flamary2021pot}. As illustrated in the figure, our source and target measures are given by  $N=30$ equally weighted Dirac masses in $d=2$ spatial dimensions. We choose $M = 3$ time steps,  $\delta = N^{-0.99/d}$, $\epsilon = 0.01$,  learning rate $\alpha=0.01$, and $n=2000$ gradient descent steps. (See section \ref{numericalerrorestimates} for a more detailed discussion of the relationship between the parameters $\delta$, $N$, and $\epsilon$ and section \ref{erroralonggraddescent} for a discussion of how these relate to the gradient descent parameters $\alpha$ and $n$.)

The top left panel of Figure \ref{fig:uniform_to_uniform} shows the trajectories $\{y_{i,j} \}$  computed by our method, linearly interpolating between the $M=3$ time steps to better illustrate the paths of the particles. The top right panel compares the terminal values of the trajectories $\{ y_{i,M} \}$ to the locations of particles in the target $\{w_i\}$. While our method uses a soft constraint to match source to target particles, we observe overall good agreement. The bottom left panel shows the optimal transport matching computed via POT-emd, which is qualitatively similar to our solution, though not identical. However, in the bottom right panel, we see by comparing the Blob OT and POT-emd solutions in terms of the value of the objective function (\ref{numericalobjective}), they are extremely close in terms of the degree to which they match source to target measure with the smallest possible kinetic energy. While the gradient descent of Blob-OT is initialized far from optimum, we see that, over the course of the gradient descent, it  converges to the nearly optimal objective function value achieved by POT-emd. After $n=2,000$ steps, Blob OT has objective function value of 1.9884 =  1.9779 (kinetic energy, \ref{numericalKE}) +  0.0105 (nonlocal energy, \ref{numericalNE}), while POT-emd has objective function value of 1.9908 = 1.9908 (kinetic energy, \ref{numericalKE}) +  0 (nonlocal energy, \ref{numericalNE}), rounded to 4 decimal points. Note that while the kinetic energy for the Blob OT is smaller than POT-emd, this is due to the fact that Blob OT has a soft terminal constraint, not because it has indeed found a superior matching between source and target.

\subsection{Continuous target Gaussian distribution} \label{gaussiansection}

\begin{figure}[h]
\centering
\includegraphics[width=.82\linewidth]{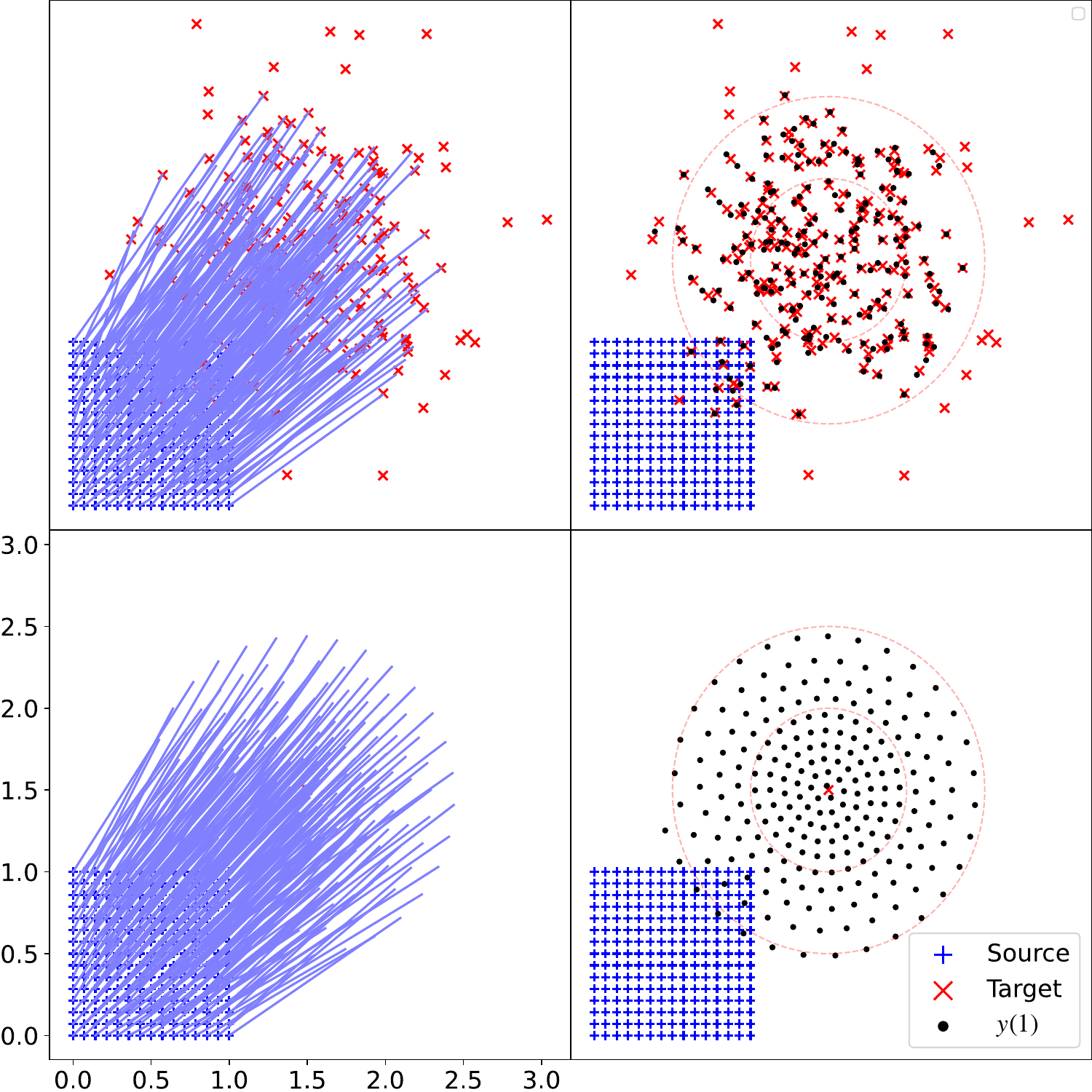}
\caption{Comparison between our particle method for computing the dynamic formulation of the 2-Wasserstein distance between a uniform grid and iid samples from a Gaussian (top row) and the continuum Gaussian itself (bottom row).}
    \label{fig:grid_to_gaussian}
\end{figure}

In Figure \ref{fig:grid_to_gaussian}, we contrast the behavior of our method for computing the 2-Wasserstein optimal transport when the target distribution is given by a continuum Gaussian distribution $\mathcal{N}(\mu, \sigma^2)$, with   $\mu = (1.5, 1.5), \sigma=0.5$, versus when the target is given by iid samples from the same Gaussian, in $d=2$ dimensions. In both cases, we take our source distribution to be $N=225$ equally weighted Dirac masses arranged on a grid, as an approximation of  the 
uniform measure on the square $m_0 = 1_{[0,1]\times[0,1]}$. We take $M=5$ time steps, $\epsilon=0.01$, $\delta = N^{-0.99/d}$, $n = 5\cdot 10^6$ gradient descent steps, and learning rate  $\alpha = 0.01\epsilon$.

In the top row of Figure \ref{fig:grid_to_gaussian}, we use our method to compute trajectories from the source measure to iid samples of the Gaussian, drawn using the numpy function randn. In the bottom row, we compute trajectories from the source measure to the continuum gaussian, taking $m_1(y) = \varphi_{0.5}(y - (1.5,1.5))$ and using the approximation $K_\delta*m_1(y) = \varphi_{\sqrt{ 0.5 + \delta^2}}(y - (1.5, 1.5)) \approx m_1(y)$ in the definition of the nonlocal energy, equation (\ref{numericalNE}).   The left column shows the trajectories $\{y_{i,j}\}$ computed by our method, linearly interpolating between time steps to show the path of each particle. The right column shows the terminal locations of the trajectories, which can be interpreted as samples of $m_1$, and compares them to the target measure. The dotted circles around the Gaussian target measure illustrate its standard deviations, $\sigma$ and $2\sigma$.

In both the top and bottom row, we  see that the particles primarily end up within two standard deviations of the Gaussian mean, largely ignoring outliers/tails. We believe this is due to the finite particle nature of our approximation, and we anticipate that, as $N \to +\infty$, $\delta \to 0$, and $\epsilon \to 0$, minimizers of (\ref{numericalobjective}) will indeed match to more outliers/tails. While our method succeeds at capturing the irregular samples in the top row, we find it especially interesting that, when the source measure is taken to be a continuum Gaussian in the second row, the final particle locations self-organize with much more regular structure. For this reason, we believe that methods based on our approach may show promise within the context of sampling, especially when the samples are used to discretize integrals of smooth functions.

\begin{figure}[h]
\centering
\includegraphics[width=.82\linewidth]{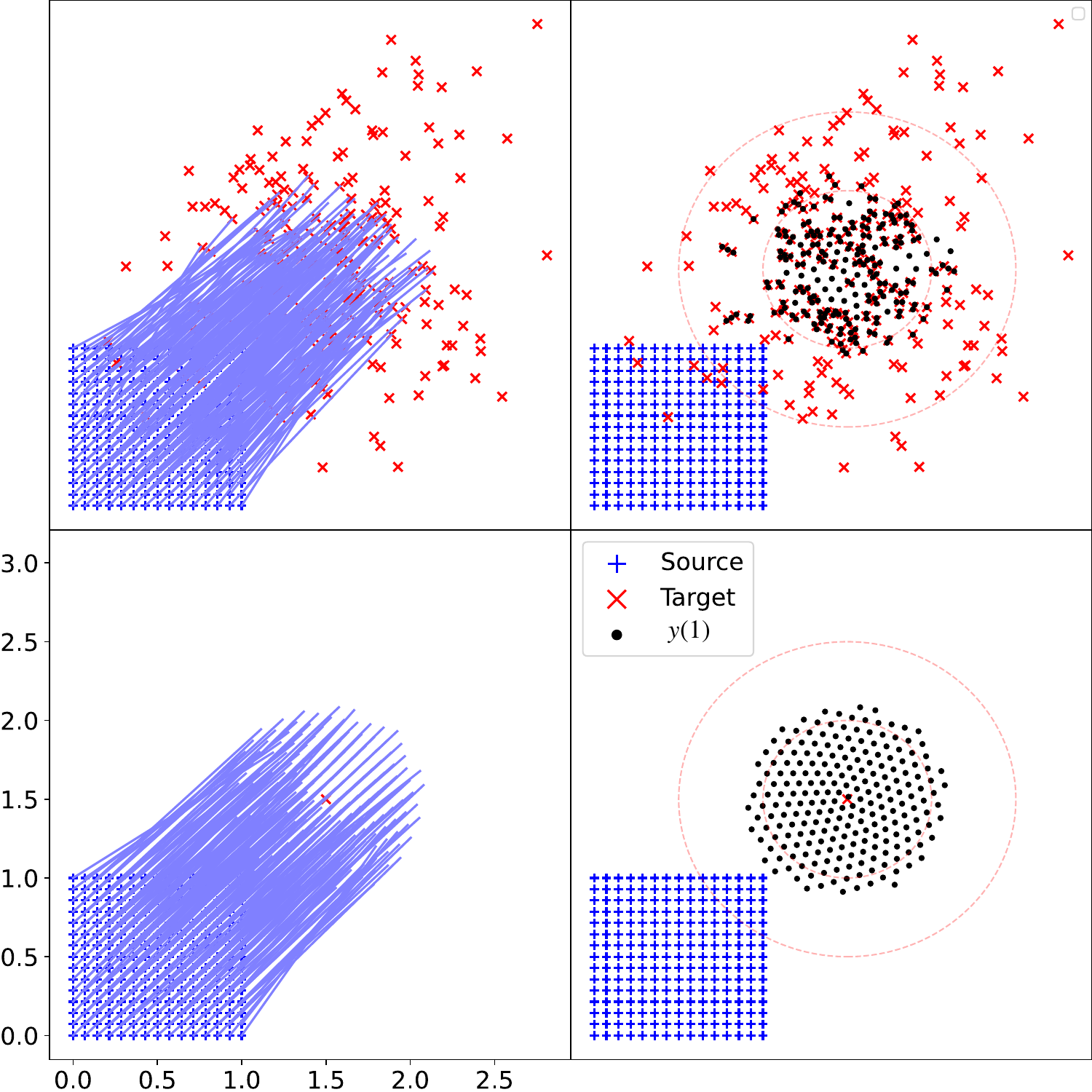}
\caption{Comparison between computing optimal transport from a uniform grid to iid samples from a Gaussian (top row) and the continuum Gaussian itself (bottom row) when the regularization parameter $\delta$ is intentionally chosen too small.}
    \label{fig:grid_to_gaussian_smalldelta}
\end{figure}

In order to illustrate the importance of the choice of parameters on behavior of our method, we illustrate how the behavior of the previous figure changes when $\delta$ is chosen too small. In Figure \ref{fig:grid_to_gaussian_smalldelta}, we consider the same numerical experiment as in the previous Figure \ref{fig:grid_to_gaussian}, except that instead of choosing $\delta = N^{-0.99/d}$, we choose $\delta=\sigma d^{-0.5} N^{-0.99/d}$. For this smaller value of $\delta$, we observe that the particles are less able to match the outliers/tails of the target Gaussian. The reason for this is similar in both the iid sample case (top row) and the continuum Gaussian case (bottom row). In the case of iid samples,  this is due to the fact that, in our objective function (\ref{numericalobjective}), the terminal locations of our particles $\{y_{i,M}\}$ are less able to sense the samples $\{w_i\}$  when $\delta$ is too small, since $y \mapsto \varphi_\delta(y - w_i)$ decays quickly away from $w_i$. Interestingly, in this case, when the terminal locations of our particles are unable to sense a sufficiently near sample, they organize themselves to be roughly evenly spaced in the empty regions.

In the case of a continuum Gaussian, $y \mapsto m_1(y)$ becomes very small outside two standard deviations of the mean. In both cases, the particle interactions mediated by the first term in the nonlocal energy $K_\delta(y_{i,M} - y_{j,M})$ become weak when $\delta$ is too much smaller than the distance between particles, preventing the particles from repelling each other and exploring more regions of parameter space.

\subsection{Obstacles} \label{obstaclesec}

\begin{figure}[h]
\centering
\begin{subfigure}[b]{0.45\textwidth}
        \includegraphics[width=\textwidth]{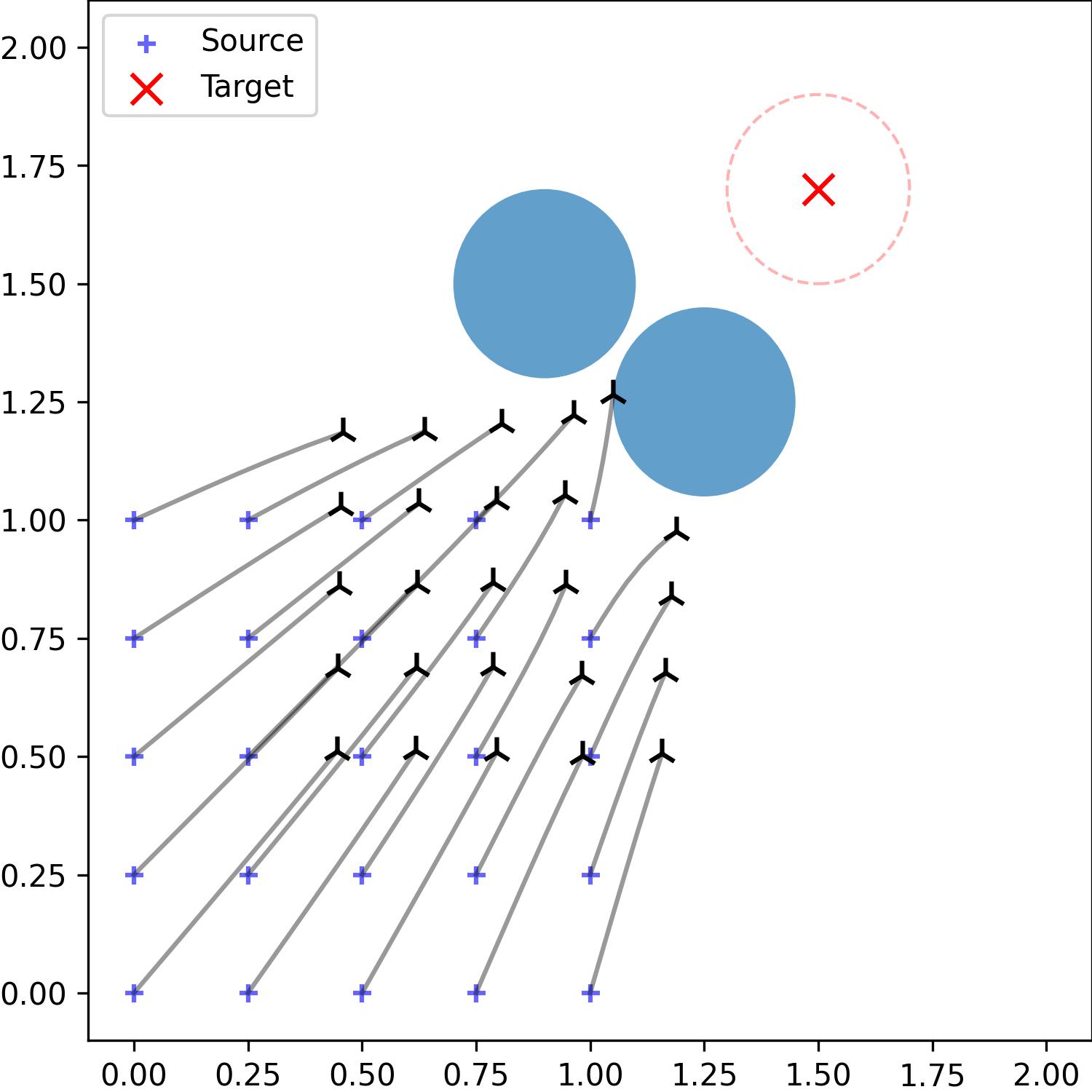}
        \caption{$t=0.3$}
\end{subfigure}%
\begin{subfigure}[b]{0.45\textwidth}
        \includegraphics[width=\textwidth]{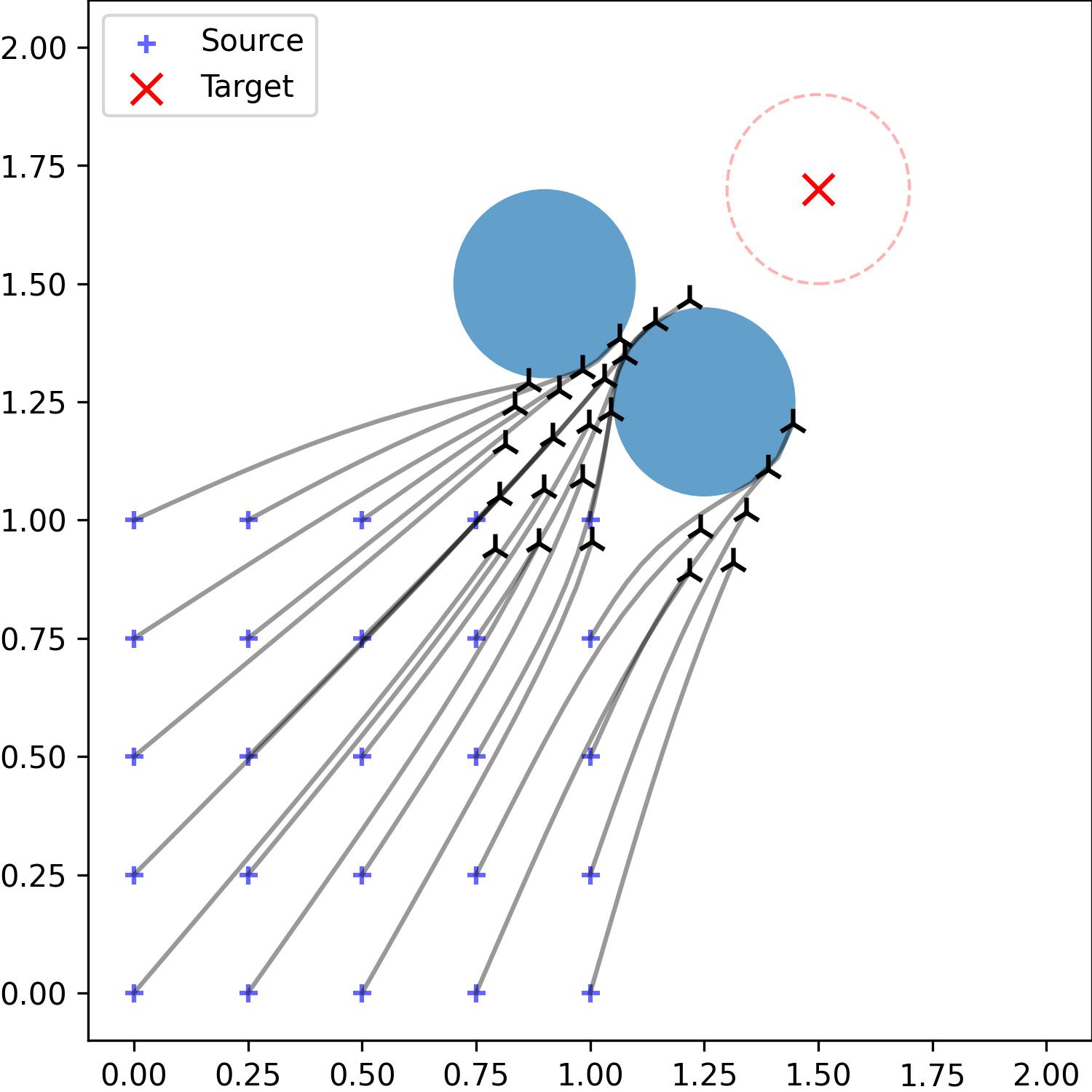}
        \caption{$t=0.55$}
\end{subfigure}
\begin{subfigure}[b]{0.45\textwidth}
        \includegraphics[width=\textwidth]{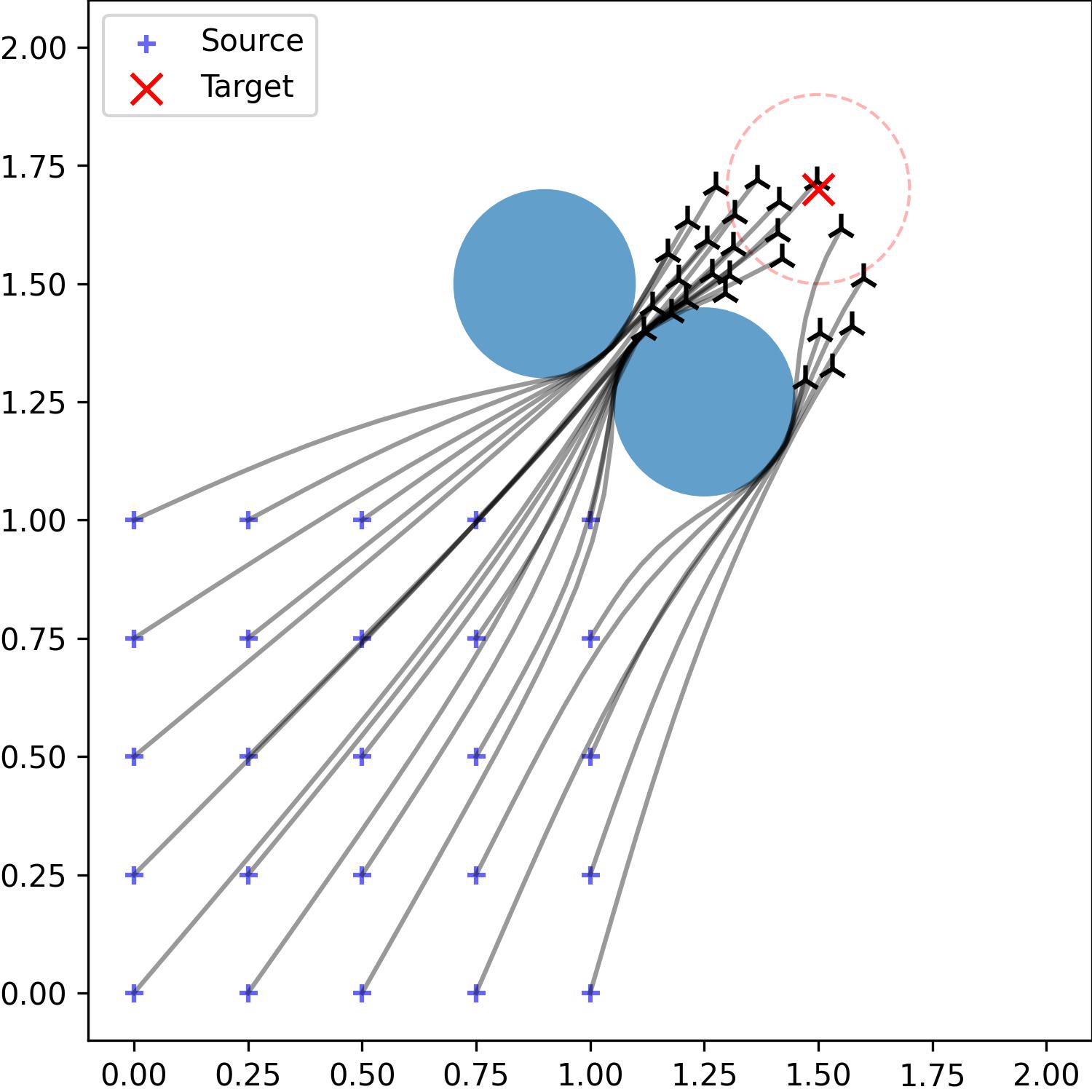}
        \caption{$t=0.85$}
\end{subfigure}%
\begin{subfigure}[b]{0.45\textwidth}
        \includegraphics[width=\textwidth]{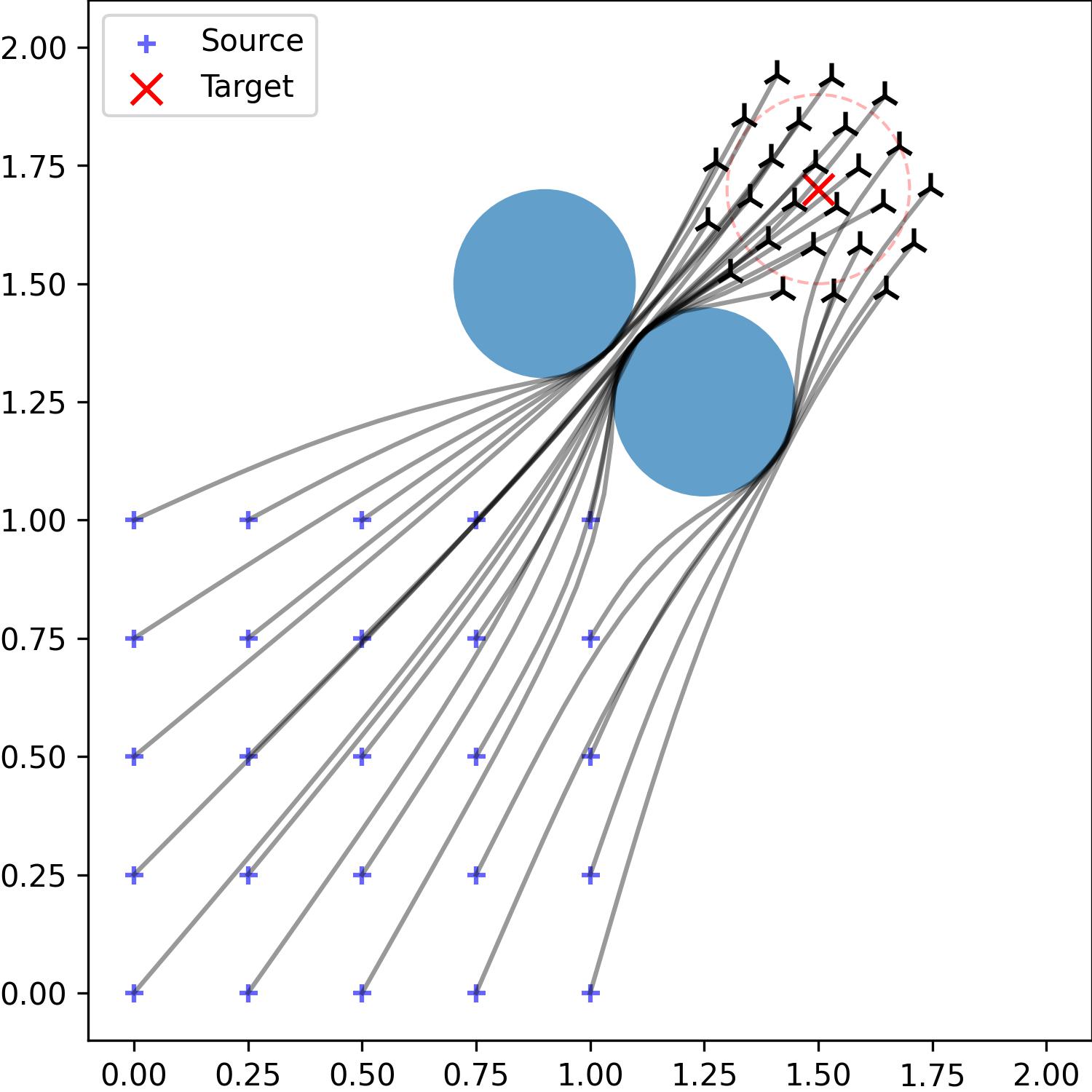}
        \caption{$t=1$}
\end{subfigure}
\caption{2-Wasserstein optimal transport from an empirical measure on a uniform grid to a continuum Gaussian distribution, in the presence of two circular obstacles.} 
\label{fig:obstacles}
\end{figure}

Building on the preceding section, we now consider optimal transport in the presence of an obstacle $\Omega$, which is given by the union of the interior of two circles of  radius $r_k = 0.2$ centered at $c_1 = (1.0,1.5)$ and $c_2 = (1.25,1.25)$ . We represent the obstacle in our optimization problem via the state and measure cost
\begin{equation}
    L_N(y) =c_\Omega(N)  \sum_{k=1}^2 \max\{r_k^2 - (y-c_k)^2, 0\} .
\end{equation}
for  $c_\Omega(N)$ satisfying $\lim_{N \to +\infty} c_\Omega(N) = +\infty$. By definition, we see that $L_N \nearrow  L $, where 
\[ L(y) = \begin{cases} 0 &\text{ if } y \in \Omega^c, \\ +\infty &\text{ if } y \in \Omega. \end{cases} \]


The source measure is given by $N=25$ equally weighted Dirac masses arranged in a uniform grid on the unit square, $m_0 = 1_{[0,1]\times[0,1]}$. For the continuum gaussian target measure, we take $m_1(y) = \varphi_{0.2}(y - (1.5,1.7))$    and use the approximation $K_\delta*m_1(y)   \approx m_1(y)$ in the definition of the nonlocal energy, equation (\ref{numericalNE}).  We take $M=21$ time steps, $\epsilon=0.1$, $\delta = N^{-0.99/d}$, $n = 2\cdot 10^5$ gradient descent steps, and learning rate  $\alpha = 0.001\epsilon$. The obstacle constant is  $c_\Omega(N)=(h\epsilon)^{-1}$.
 
 Figure \ref{fig:obstacles} shows the evolution of the linear interpolations of the trajectories $y_i(t)$, $i = 1, \dots N$ at times $t=0.3, 0.55, 0.85$, and $1.0$. For the time continuous problem in the limit $\delta \to 0, \epsilon \to 0$, the particles should follow straight lines, bending only to follow the boundaries of the obstacles and straightening again as they leave the obstacle and approach their terminal points. Overall, we observe good agreement between our numerical approximation and the continuum solution, with only mild bending close to the obstacle.

%
%
%

\subsection{Acceleration control}\label{accsection}

We now consider the performance of our method in the case of measure transport subject to acceleration controls. We begin by describing how to discretize the time continuous approach, in analogy with our approach for the velocity control problem, described at the beginning of section \ref{numericalimplementation}. 

In order to solve (\ref{MFCa}), we apply the nonlocal terminal constraint and particle discretization (\ref{MFCN}), to arrive at the continuous time formulation: 
			\begin{align} \label{acccts1}
			&\inf_{x_i, v_i, a_i} \frac{1}{N}\sum_{i=1}^N\int_0^1|a_i(t)|^2 dt \\
			&   + \frac{1}{\epsilon} \left[ \frac{1}{N^2}\sum_{i,k=1}^N  K_\delta((x_i(1),v_i(1))- (x_k(1),v_k(1))) - \frac{2}{N}\sum_{i=1}^N(K_\delta*m_1)(x_i(1),v_i(1)) \right] \nonumber
			\end{align}
			subject to the constraints
			\begin{eqnarray}
			\begin{cases} \dot{x_i}(t) &= v_i(t)  \\ x_i(0)&=x_{i,0} \end{cases} \qquad \qquad \label{acccts3}
			\begin{cases} \dot{v}_i(t) &= a_i(t) \\  v_i(0)&=v_{i,0} \end{cases} \label{acccts2}
			\end{eqnarray}
As before, we discretize time $t \in [0,1]$ on a uniform grid with $M$ grid points and time step $h = 1/(M-1)$, approximating $x_i(t)$ by the vector
   \begin{align}
    [x_{i,j}]_{j=1}^M &\approx \left[x_{i,0}, x_{i,0}+h v_{i,0}  ,x_i(2h),  \ldots, x_i(1-h), x_i(1) \right],
\end{align}
which, by definition, incorporates initial condition  constraints  (\ref{acccts3}). As before, we approximate the velocity and the acceleration by first order finite differences,
\begin{align*}
v_{i,j} := \frac{x_{i,j} - x_{i,j-1}}{h} , \quad a_{i,j} := \frac{v_{i,j} - v_{i,j-1}}{h} .
\end{align*}
Substituting into the objective function (\ref{acccts1}), we arrive at the fully discrete problem, which is to minimize the sum of the control cost and the nonlocal position/velocity energy 
\begin{align*}
  \min_{ { x}_{i,j}} ~& \left\{ {\rm CC}(x_{i,j}) + {\rm NPVE}(x_{i,j}) : x_{i,1} = x_{i,0}, \  x_{i,2} = x_{i,0} + h v_{i,0}  \ \forall i = 1, \dots, N\right\}
\end{align*}
for
\begin{align*} \ 
{\rm CC}(x_{i,j})  : = &\frac{(M-1)^3}{  N}\sum_{i=1}^N\left|\texttt{ddiff}( { x}_{i,j})\right|^2  \\
{\rm NPVE}(x_{i,j}) := & \\
\frac{1}{\epsilon} &\left[ \frac{1}{N^2}\sum_{i,k=1}^N  K_\delta((x_{i,M},v_{i,M})- (x_{k,M},v_{k,M}))  - \frac{2}{N}\sum_{i=1}^N (K_\delta* m_1)(x_{i,M},v_{i,M})  \right]
\end{align*}
where $i = 1, \dots, N$ and $ j = 2, \dots, M$, with
\begin{eqnarray*}
    \texttt{ddiff}( { x}_{i,j}) := \left[x_{i,1}-2x_{i,2}+x_{i,3},x_{i,2}-2x_{i,3} +x_{i,4} ,  \ldots, x_{i,M-2} -2x_{i,M-1} +x_{i,M}, 0 ,0\right] .
\end{eqnarray*}

As in the velocity control case, while the continuum problem motivating the above problem is  well-posed only when $m_1 \in L^2(\Rd \times \Rd)$,  one particular case of interest is the case when the target measure  $m_1$   is given by a  sum of Dirac masses, $\frac{1}{N}\sum_{i=1}^N\delta_{(w_i^x, w_i^v)}$. This can be handled within our framework analogously to the previously described velocity control case.

\begin{figure}[h]
\centering
\includegraphics[width=\linewidth]{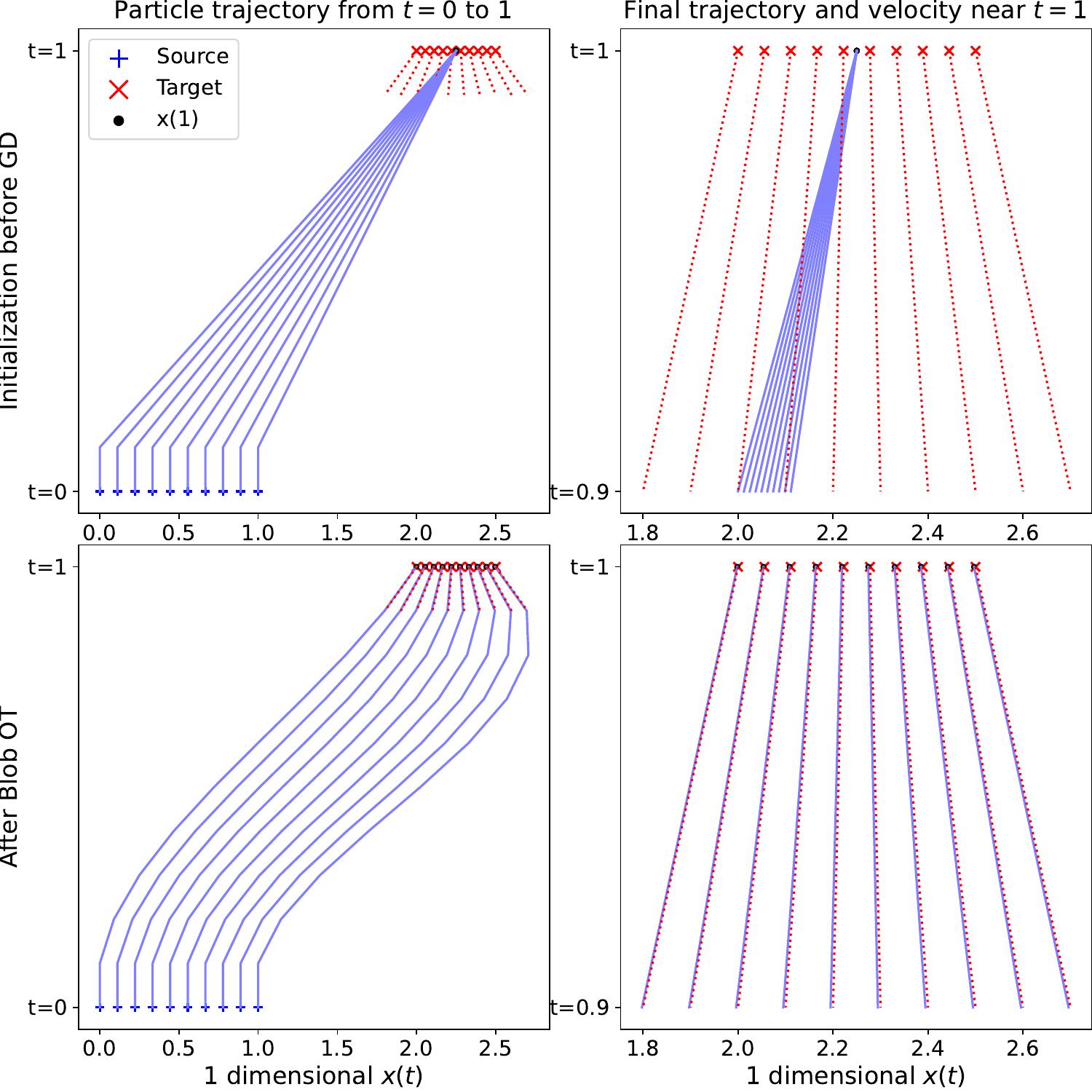} 
  \caption{Measure transport subject to acceleration control in one spatial dimension, between two well-ordered empirical distributions. In the top row, we show the initialization of our gradient descent, and in the bottom row, we show the approximate optimizer.  }
      \label{fig:acc_grid}
\end{figure}

In Figure \ref{fig:acc_grid}, we plot the results of an acceleration control optimization problem in one spatial dimension, in the case that the source distribution is an empirical measure with $N =10$ particles evenly spaced on $[0,1]$, with initial velocities zero, and the target distribution is an empirical measure with $N =10$ particles evenly spaced on $ [2, 2.5] $, with target velocities evenly spaced on $[-2,2]$. Note that this example satisfies the condition of our main theorem that $\supp m_0 \subseteq N(A)$. We consider $M = 11$ time steps, $ \epsilon=10^{-4}$, $\delta = N^{-0.99/2}$, $n = 10^6$ iterations, and $\alpha = 4 \epsilon / 10^3$ learning rate. (Note that the dimension in our choice of $\delta$ is two dimensional, since our mollifier $\varphi_\delta(x,v)$ is a function on $\R \times \R$.)
 
 The top left panel of Figure \ref{fig:acc_grid} shows our initialization of the gradient descent $x_{i,j}$, linearly interpolating in time. The top right panel zooms in on this initialization and compares it to the desired target distribution. The bottom left panel shows the approximate minimizer $x_{i,j}$ obtained via gradient descent, again linearly interpolating in time, and the bottom right panel compares this to the desired target distribution. The approximate minimizer computed by our method exhibits trajectories with low acceleration and that agree with the desired target distribution in position and velocity.

%
%
%

\subsection{Illustrating the Loss Landscape} \label{losslandscapesection}
\begin{figure}[h!]
    \includegraphics[width=\linewidth]{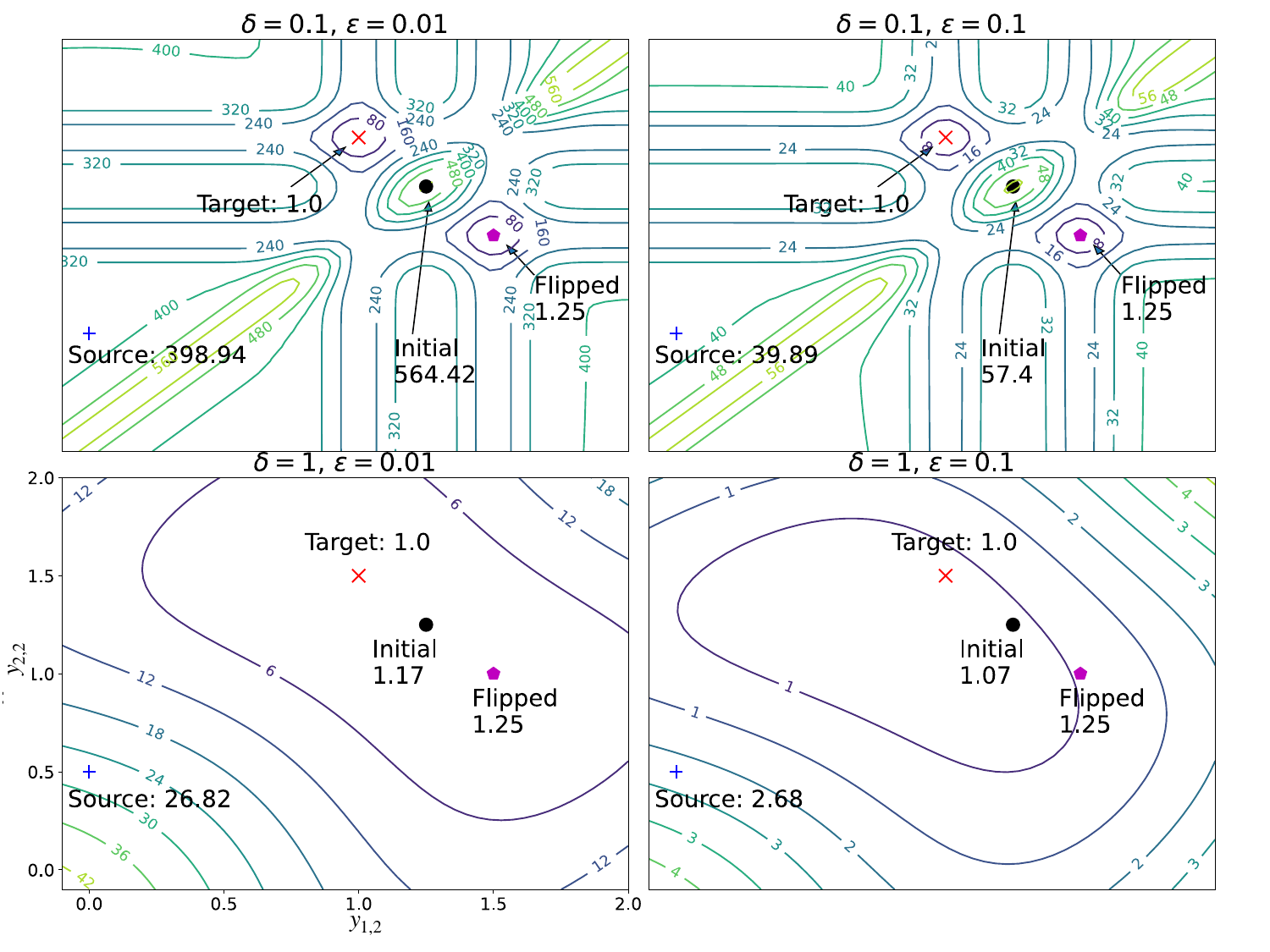}
\caption{Contour plots of nonconvex loss landscape for our optimization problem (\ref{numericalobjective}) between two Dirac masses, as a function of the terminal points $y_{1,2}$ and $y_{2,2}$,  for different choices of  $\epsilon$ and $\delta$.}
\label{fig:loss_contours}
\end{figure}

We now turn to  an example illustrating basic   properties of the loss landscape for our particle discretization of the optimal transport problem (\ref{numericalobjective}). For simplicity of visualization, we consider $L_N=0$, $N = 2$ and $M=2$, with one dimensional source distribution $ \frac{1}{2}(\delta_0 +\delta_{0.5})$ and target distribution $\frac{1}{2}(\delta_1 +\delta_{1.5})$. In Figure \ref{fig:loss_contours}, we plot the loss landscape as a function of $y_{1,2}$ and $y_{2,2}$, for different choices of $\delta$ and $\epsilon$.  

For reference, we have plotted the following values on the loss landscape:
\begin{itemize}
\item \textbf{Source distribution}, {\color{RoyalBlue}$+$}: the value of the objective function if the particles stay at their initial locations $(y_{1,2}, y_{2,2}) = (y_{1,1}, y_{2,1}) = (0, 0.5)$;  we observe that the objective function is large at the source distribution, reflecting the fact that the particles are far from the desired target distribution.
\item \textbf{Target distribution}, $\color{red} \times$: the value of the objective function if the particles are optimally transported from source to target, with the leftmost particle in the source getting mapped to the leftmost particle in the target $(y_{1,2}, y_{2,2})  = (1.0, 1.5)$;  we observe that the objective function is minimized for this configuration.
\item \textbf{Flipped distribution}, $\color{Purple} \pentagonblack$: the value of the objective function if the particles in the source are exactly mapped to the particles in the target, but in a nonoptimal order, with the leftmost source particle mapped to the rightmost target particle $(y_{1,2}, y_{2,2})  = (1.5, 1.0)$; we observe that the value of the objective function is small, since the terminal configuration agrees exactly matches the target distribution, but not minimal, since the kinetic energy of achieving this configuration is not as small as possible.
\item \textbf{Initial distribution}, $\bullet$: the value of the objective function at our initialization of gradient descent, as described in equations (\ref{typicalinitialization}-\ref{typicalinitialization2}).
\end{itemize}
Note that the value of the objective function is the same at the ``Target'' and ``Flipped'' distributions in all four plots, due to the fact that, for both of these configurations, the value of the nonlocal energy term in the objective function (\ref{numericalobjective}) is zero for any $\delta$ and $\epsilon$.

As anticipated, for all values of $\delta$ and $\epsilon$,   the loss landscape is   nonconvex.  Comparing the first row and second row of Figure \ref{fig:loss_contours}, we observe that smaller values of $\delta$ (top row) better distinguish between small changes in $(y_{1,2}, y_{2,2})$, due to the fact that, when $\delta$ is small, there is less smoothing in our approximation of the source and target measures.

Comparing the first column and second column of Figure \ref{fig:loss_contours}, we observe that smaller values of $\epsilon$ (left column) place more weight on the nonlocal energy term on the objective function (\ref{numericalobjective}),  rewarding proximity to the target measure. We observe that the plots in the left column are highly symmetric across $y_{1,2} = y_{2,2}$, due to the fact that the objective function primarily considers the final locations of the particles, rather than whether the particles were optimally transported to those locations. On the other hand, larger values of $\epsilon$ (right column) place more weight on the kinetic energy term. In this case, we observe less symmetry across   $y_{1,2} = y_{2,2}$, due to the fact that  the kinetic energy term prioritizes moving particles a shorter distance from source to target.

Finally, we observe that larger values of $\delta$ (bottom row) and $\epsilon$ (right column) flatten the energy landscape. In this way, when using a gradient-based optimization method to compute minimizers of the objective function (\ref{numericalobjective}), an optimal choice of learning rate will depend on the choices of $\delta$ and $\epsilon$. Since, in the present simulations, we typically choose $\delta$ very small, the size of $\epsilon$ is the main factor in the flatness of the energy landscape. For this reason, we choose our learning rate $\alpha$ to scale with $\epsilon$, taking bigger steps on the flatter energy landscape and smaller steps on the steeper landscape.

\subsection{Estimating the error from the optimal transport map} \label{erroralonggraddescent}

\begin{figure}[h]
\centering
    \includegraphics[width=.85\linewidth]{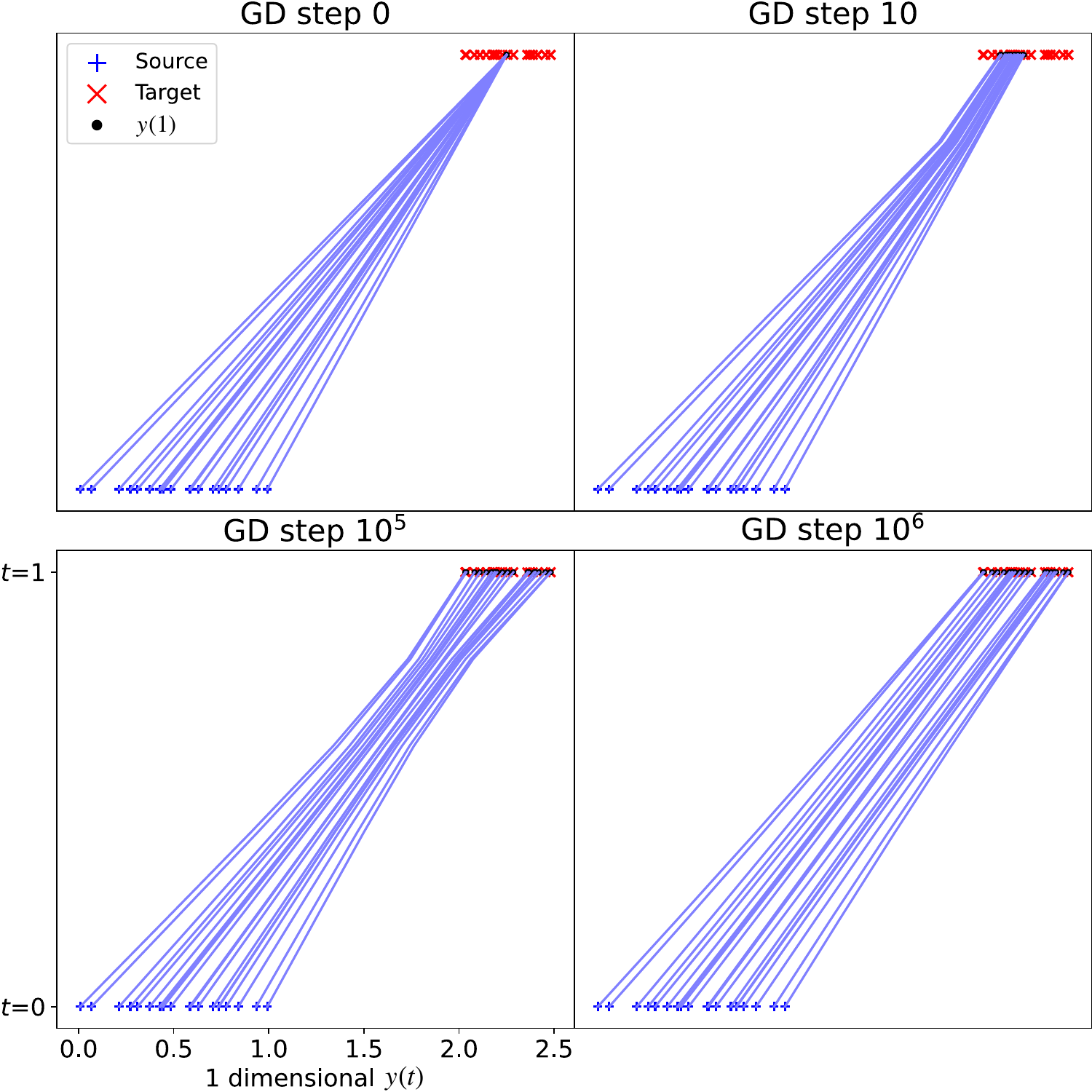}
     \caption{Illustration of how trajectories change along gradient descent to compute an approximate minimizer of (\ref{numericalobjective}),  when computing 2-Wasserstein optimal transport problem between two empirical measures.  }
     \label{fig:small_eps_traj}
\end{figure}

%
%
%
%

We now analyze the error between the approximate solution of the 2-Wasserstein optimal transport problem  computed by our method and the exact optimal transport map. In particular, we examine how this error behaves along the gradient descent which computes the approximate optimizer. On one hand, we expect higher accuracy when $\epsilon$ is small, given that this more strongly imposes the terminal constraint via the nonlocal energy (\ref{numericalobjective}). On the other hand, due to the fact that we let the learning rate $\alpha$ depend  on $\epsilon$ (see sections \ref{parameterselection} and \ref{losslandscapesection}), when $\epsilon$ is small, we also have to run gradient descent for more iterations to compute our approximate minimizer. Consequently, in practice, one seeks a value of $\epsilon$ that is small enough to give accurate results but large enough to be computationally efficient. 

In the following figures,   we consider optimal transport in one dimension, where the source measure is given by  $m_0 = 1_{[0,1]}$, discretized as $N=20$ particles on a uniform grid on $[0, 1]$. The target measure is given by  $m_1 =2 \cdot 1_{[2,2.5]}$. In both cases, we take nonlocal regularization $\delta=N^{-0.99}$, $M= 5$ time steps, and learning rate  $\alpha=\max(0.001\epsilon,10^{-5})$ .

At the continuum level, the 2-Wasserstein geodesic from $m_0$ to $m_1$ at time  $t$ is given by
\begin{align} T(y, t):= (1-t)y + t(0.5y+2) .
\end{align}
Likewise, $T(y,1)$ is the optimal transport map from $m_0$ to $m_1$.
We consider two measures of error: the average error across all time points in our discretization,
\begin{equation}
\text{[Error at all time points]} = \sqrt{\frac{1}{N(M-1)}\sum_{i=1}^N\sum_{j=2}^{M}\left|y_{i,j}  - T \left(y_{i,1}, \frac{j-1}{M-1} \right)\right|^2}, \label{eq:error_ot_all}
\end{equation}
and the error at terminal time $t=1$,
\begin{equation}
\text{[Error at }t=1] = \sqrt{\frac{1}{N}\sum_{i=1}^N\left|y_{i,M}  - T(y_{i,1}, 1)\right|^2}. \label{eq:error_ot_end}
\end{equation}
We also analyze the decay of both the kinetic energy term (\ref{numericalKE}) and the nonlocal energy term (\ref{numericalNEempirical}), as well as the total loss, given by their sum (\ref{numericalobjective}). 

Figure \ref{fig:small_eps_traj} shows the gradient descent dynamics for  a very small choice of $\epsilon = 0.001$ and $n = 10^6$ gradient descent steps. The top left panel shows the initialization of the gradient descent $y_{i,j}$, the top right panel shows $y_{i,j}$ after $10$ gradient descent steps, the bottom left shows the behavior after $10^5$ steps, and the bottom right panel shows the behavior after $10^6$ steps.  
 \begin{figure}[h]
     \centering
    \includegraphics[width=\linewidth]{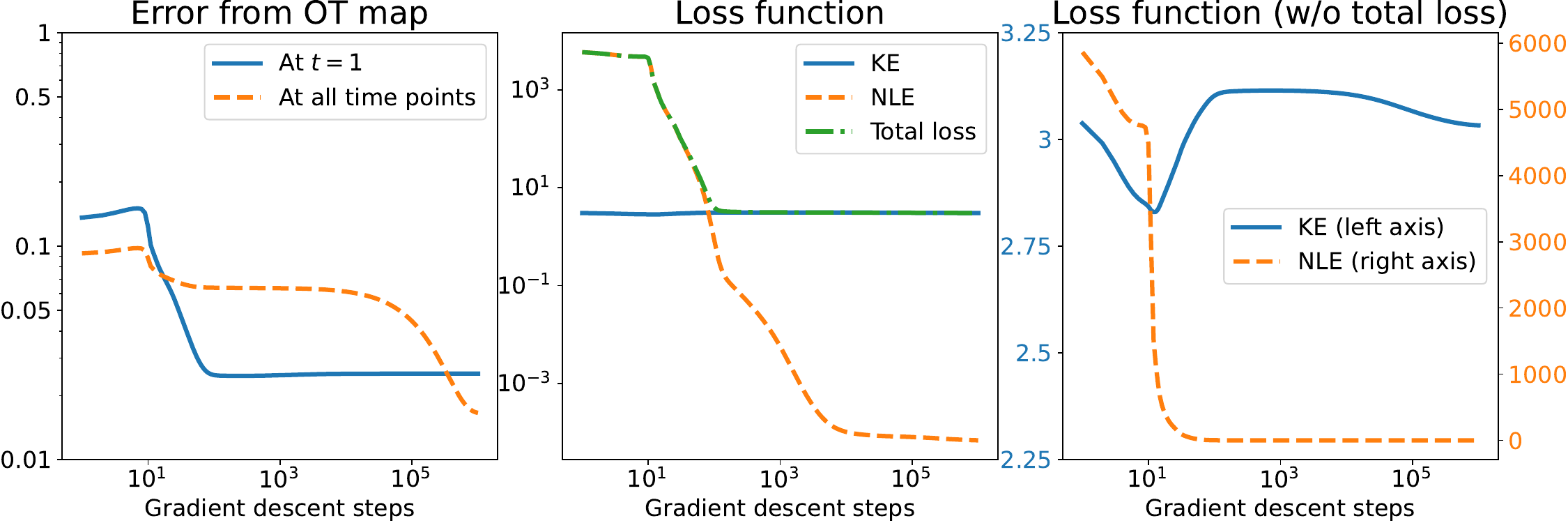}
    \caption{Error  analysis between the approximate solution of an  optimal transport problem,  computed by our method, and the exact  solution of the optimal transport problem. Left: While the error at $t=1$ quickly decays, many more gradient descent steps are required for the error at all time points to decay. Middle: Illustration of how the kinetic energy and nonlocal energy contribute to the total loss. Right: While the nonlocal energy quickly decreases to its optimal value, many more gradient descent steps are required for the kinetic energy to approach its optimal value. Note the dual axis.}
    \label{fig:small_eps_obj}
\end{figure}

A similar phenomena can be observed in Figure \ref{fig:small_eps_obj}, which shows the behavior of the error and loss function along the gradient descent iterations. In the left plot, we see that both the error at terminal time and the error at all time points decay to zero, though the error at all time points requires a very large number of iterations to decay. This reflects the fact that the very small value of $\epsilon$ forces the trajectories to match the terminal points more strongly than it enforces the kinetic energy, which straightens the trajectories. In the middle plot, we observe the decay of both the kinetic energy and nonlocal energy along iterations. In the right plot, plotting the kinetic energy and nonlocal energy on different axes, we see again that the small value of $\epsilon$ causes the nonlocal energy to decay first, while the kinetic energy takes many more iterations to decay.

\subsection{Numerical analysis of rate of convergence to continuum} \label{numericalerrorestimates}

\begin{figure}[h]
    \centering
\includegraphics[width=.8\linewidth]{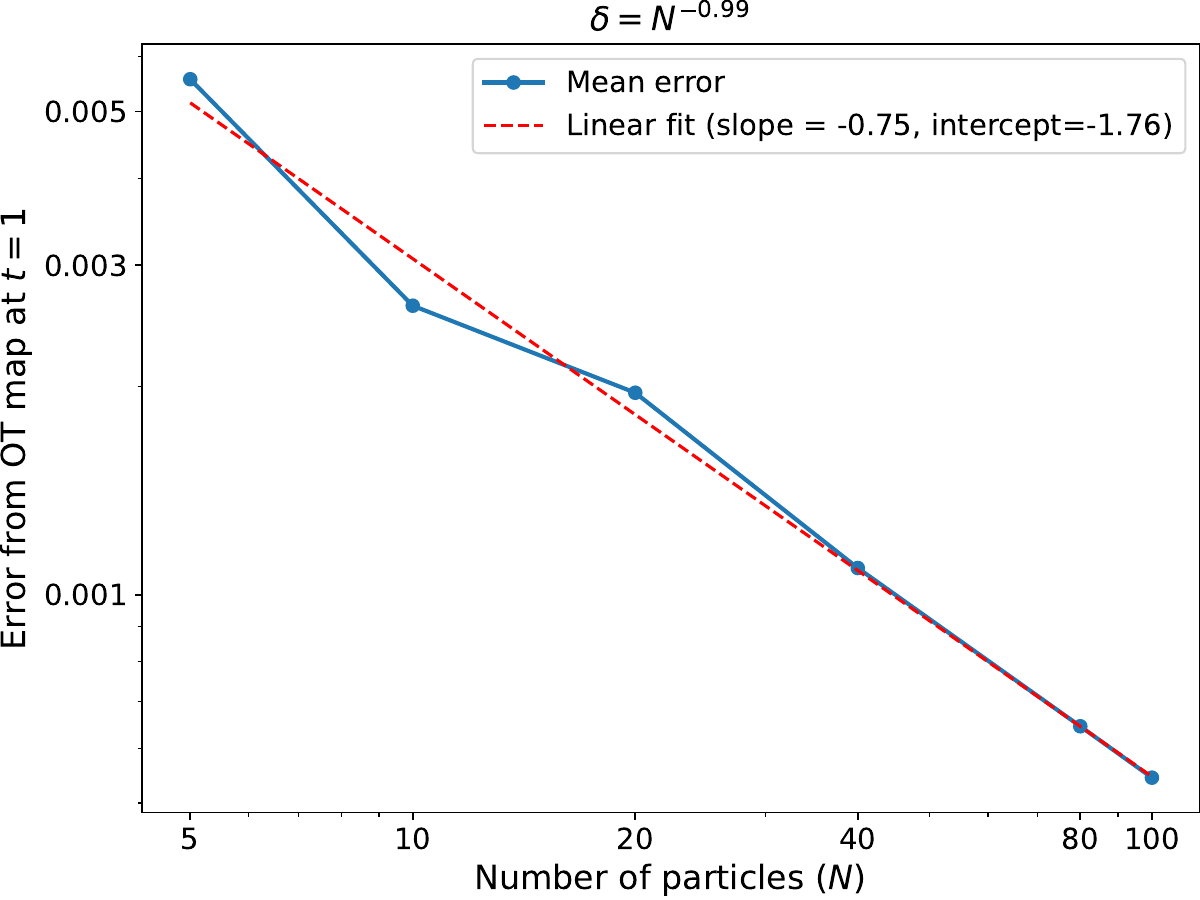}
    \caption{Analysis of the rate of convergence of our particle approximation of 2-Wasserstein optimal transport to the exact solution   in one dimension, as $N \to +\infty$ and $\delta = N^{-0.99} \to 0$.  We observe slightly slower than first order convergence in $N$.
}
    \label{fig:delta_scaling_grid}
\end{figure}

We conclude by analyzing the rate of convergence to the continuum formulation of the optimal transport problem as $\epsilon \to 0$, $\delta \to 0$,  and $N \to +\infty$. We consider one dimensional source and target distributions given by $m_0 = 1_{[0,1]}$ and $m_1 = 2 \cdot 1_{[2,2.5]}$, discretized on a uniform grid with $N$ particles, where we allow $N$ to vary. We consider $\epsilon = 0.01$, $M =5$ time steps, $n = 2 \cdot 10^6$ gradient descent steps, and learning rate $\alpha = 0.003\epsilon$. 
As explained in section (\ref{parameterselection}), we allow $\delta \to 0$ and $N \to +\infty$ simultaneously according to the rate in equation (\ref{deltaNscaling}), with $k = 0.99$. 

In order to obtain numerical estimates on the rate of convergence as $N \to +\infty$ and $\delta \to 0$, we consider the error at the terminal time $t=1$, see equation (\ref{eq:error_ot_end}), using the fact that, in the present example the optimal transport map for the continuum problem is given by $T(y,1) = y/2+2$.

The results of our numerical study are shown in Figure \ref{fig:delta_scaling_grid}. Here, we plot how the error from the optimal transport map at time $t=1$ varies as the number of particles $N$  increases from 5, 10, 20, 40, 80, to 100. Comparing the decay of the error on a log-log scale to the line of best fit (calculated by minimizing least squares deviation as implemented in NumPy \cite{harris2020array}), we observe slightly slower than first order convergence as $N \to +\infty$.


\appendix

\section{Basic Properties of Mean Field Control Problems} \label{basicpropertiesappendix}
In this section, we collect several results and proofs regarding basic properties of mean field control problems. 
We begin by proving Lemma \ref{fancydoubling}, on elementary properties of our control cost $\psi$ and admissible function $\phi$.
 \begin{proof}[Proof of Lemma \ref{fancydoubling}]
Part (\ref{psisuperlinearity}) is an immediate consequence of the fact that $\phi$ is superlinear at $+\infty$ and inequality (\ref{phipsirelation}).
 
Now, we show part (ii). To begin, we show that  there exists  $K'' >1$ so that
\begin{align} \label{pepsi}
\phi(2r) + 2r \leq K''(\phi(r)+ r) , \quad \text{ for all } r \in [0,+\infty).
\end{align}
For any $\phi$ satisfying Assumption \ref{asmp1prime}, there exists $r_0$ so that $\phi(r_0)=1$. Then, for all $r \geq r_0$, the fact that $\phi$ is increasing ensures
\begin{align*}
\phi(2r) \leq K(1+\phi(r)) \leq 2K \phi(r)  \implies \phi(2r) + 2r \leq 2K \phi(r) +2r .
\end{align*}
On the other hand, by the convexity of $\phi$, $\phi(r) + \phi'(r)(s-r) \leq \phi(s)$ for all $r, s \in [0,+\infty)$. Thus, for all $r \leq r_0$,
\begin{align*}
\phi(2r) \leq \phi(r) + r\phi'(2r) \leq \phi(r) + r \phi'(2r_0)  \implies \phi(2r) + 2r \leq \phi(r) +(2 + \phi'(2r_0))r .
\end{align*}
Therefore, letting $K'' = \max\{ 2K, 2+ \phi'(2r_0) \}$ gives inequality (\ref{pepsi}).

To prove part (ii), note that, for any $D\geq 1$, choosing $n = \lceil \log_2(D) \rceil$,
\begin{align*}
\phi(Dr) +Dr &\leq \phi(2^n r) + 2^n r \leq (K'')^n (\phi(r) + r) \leq (K'')^{\log_2(D)+1} (\phi(r) + r) \\
&= K'' D^{\log_2(K'')} (\phi(r) + r).
\end{align*}
\end{proof}

Now, we prove Lemma \ref{momequiv} on the equivalence of the original formulation of the mean field control problem (\ref{star}) and the formulation in momentum coordinates (\ref{mainproblem}).

\begin{proof}[Proof of Lemma \ref{momequiv}]
By definition of distributional solutions to the respective PDE constraints, as in equations (\ref{firstdistributionalsolution}) and (\ref{seconddistributionalsolution}), if $(\mu, \bu) \in \mathcal{A}(m_0,m_1)$ and the value of the (\ref{star}) objective function is finite, then  $(\mu, \bu \ d \mu_t \otimes dt) \in \mathcal{C}(m_0)$, $\F(\mu_1)  = 0$, and $\E(\mu,\bu \ d \mu_t \otimes dt) <+\infty$.  Conversely, if $(\mu,\bnu) \in \mathcal{C}(m_0)$ and $\E(\mu,\bnu)<+\infty$, then $\Psi(\bnu|\mu)<+\infty$ and $\F(\mu_1)<+\infty$. The fact that $\Psi(\bnu|\mu)<+\infty$ implies that $d\bnu(y,t) = \bu(y,t)d\mu_t(y)   dt$ for some $\bu \in L^1_{d \mu_t \otimes dt}(\Rd \times [0,1];U)$, and the fact   $\F(\mu_1)<+\infty$ implies $\mu_1 = m_1$. Thus $(\mu, \bu) \in \mathcal{A}(m_0,m_1)$ and at this point, the value of the objective function in (\ref{star}) is finite.

In this way, there is a one to one correspondence between feasible points for (\ref{star}) and (\ref{mainproblem}). Furthermore, at any feasible point of either problem, the value of the objective functions coincide. This gives the result.
\end{proof}

Now, we prove Lemma \ref{lem:epcont}, which shows the energy $\G$ that appears in (\ref{mainproblem}) is lower semicontinuous.

\begin{proof}[Proof of Lemma \ref{lem:epcont}]
First, note that convergence in $C([0,1],\mathcal{P}_1(\mathbb{R}^d)) $ implies narrow convergence in $\P(\Rd \times[0,1])$.  The lower semicontinuity of the functional $\Psi$ follows from \cite[Lemma 5.4.4]{ambrosio2005gradient}.

It remains to consider lower semicontinuity of the second term in $\G$. If $L$ satisfies assumption (\ref{Llsc}), it is is independent of the measure  and lower semicontinuous, and this is an immediate consequence of \cite[Lemma 5.1.7]{ambrosio2005gradient}. Now suppose $L$ satisfies assumpetion (\ref{Lcts}), so that it depends on the measure but is uniformly continuous and real-valued. Suppose $(\mu^n,\bnu^n)_{n=1}^{\infty}$ is a sequence in $C([0,1],\mathcal{P}_1(\mathbb{R}^d)) \times \mathcal{M}(0,T \times \mathbb{R}^d, U)$ converging to a limit $(\mu,\bnu) $. Then,
\begin{align*}
&   \big |\int_{\Rd}L(\cdot ,\mu_t)d\mu_t  - \int_{\Rd} L(\cdot,\mu_t^n)d\mu_t^n  \big | \\
& \quad  \leq 
   \big |\int_{\Rd}L(\cdot,\mu_t)d\mu_t - \int_{\Rd} L(\cdot,\mu_t)d\mu^n_t \big |  +  \int_{\Rd}\big |L(\cdot,\mu_t) - L(\cdot,\mu^n_t) \big |d\mu^n_t 
\end{align*}
As $n \to +\infty$, both terms tend to zero uniformly in  $t \in [0,1]$, since  $L$ is jointly uniformly continuous. Fatou's lemma then ensures that
\begin{equation}
\int_0^T \int_\Rd L(x,\mu_t) d \mu_t(x)  \leq  \liminf_{n \rightarrow \infty} \int_0^T \int_\Rd L(x,\mu^n_t) d \mu_t(x)  .\end{equation}
\end{proof}

Next, we prove the following  estimate on the time regularity of feasible measures.
\begin{proposition} \label{pdeconstraintprop}
Under the hypotheses of Assumption \ref{asmp1}, suppose
\begin{align*}
&(\mu,\bnu) \in C([0,1]; \P_1(\Rd)) \times \mathcal{M}( \Rd \times [0,1]; U) ,
\end{align*}
with $\mu_0 = m_0 $ and $
   \partial_t \mu_t  + \nabla \cdot \left( (\bF(\cdot, \mu_t) \mu_t  + \bnu_t \right) = 0$, in the sense of distributions. Furthermore, assume   $\Psi(\bnu|\mu)<+\infty$. 
 Then, there exists $C'>0$ depending on $C_F, C_F', \psi$, $\Psi(\bnu |\mu)$, and $M_1(m_0)$ so that
 \begin{align} \label{equictssecond2}
 W_1(\mu_s,\mu_t) &\leq C'(t-s) + \int_s^t \int_\Rd |\bu(\cdot, r)| d \mu_{r} dr .
 \end{align}
 \end{proposition}
 
 \begin{proof}
 Since $\Psi(\bnu|\mu)<+\infty$, we have $d\bnu(y,t) = \bu(y,t)d\mu_t(y)   dt$ where
 \[  \int_0^1 \int_\Rd \psi(\bu(y,t))d{\mu}_t(y)dt  <+\infty . \]
Thus, recalling our notions of distributional solution from equations (\ref{firstdistributionalsolution}) and (\ref{seconddistributionalsolution}), we see that  $\mu$ is a distributional solution of the continuity equation with velocity
  \[ {\mathbf{v}}(y,t):=\bF(y, \mu_{t})+ \bu(y,t). \] 
  Furthermore, by Lemma \ref{fancydoubling}, we have
\begin{align} \label{linearphibound2}
  \int_0^1 \int_\Rd |\bu(y,t)| d \mu_{t}(y) dt \leq   \int_0^1 \int_\Rd \psi(\bu(y,t)) d \mu_{t}(y) dt + R_\phi +1   <+\infty .
\end{align}
 
Combining the above inequalities with our assumptions on $\bF$, we have
 \begin{align} \label{L1velocitybound}
 \left\|\bF(\cdot, \mu_{t})+ \bu(\cdot,t) \right\|_{L^1(\mu_{t})}   &\leq   \int_\Rd ( C_F+ C_F'|y| + C_F' M_1(\mu_t) ) d \mu_{t}(y) + \int_\Rd \left| \bu(\cdot,t) \right| d \mu_{ t} \nonumber \\
 &\leq C_F + 2C_F'M_1(\mu_{t}) + \int_\Rd \left| \bu(\cdot,t) \right| d \mu_{ t} \nonumber \\
 &= C_F + 2C_F' W_1(\mu_{t} ,\delta_0) + \int_\Rd \left| \bu(\cdot,t) \right| d \mu_{t} .   \end{align}
 In particular, since $\mu \in C([0,1];\P_1(\Rd))$, we have $t \mapsto W_1(\mu_{t} ,\delta_0)$ is uniformly bounded on $[0,1]$. Combining this with (\ref{linearphibound2}), we see that the right hand side  of equation (\ref{L1velocitybound}) is integrable in time. 
 
 Thus, by \cite[Theorem 3.4]{ambrosio2014continuity}, there exists $\bEta \in \mathcal{M}(\Rd \times C([0,1]; \Rd))$ so that $\bEta$ is concentrated on sets of pairs so that $\gamma$ is an absolutely continuous integral solution of 
 \[ \dot{\gamma}(t) = {\mathbf{v}}(\gamma(t),t) , \ \gamma(0) = x \]
 and $\mu_t = e_t   \# \bEta$, where    $e_t: \Rd \times C([0,1];\Rd) \to \Rd : (x,\gamma) \to \gamma(t)$. Therefore, $(e_t,e_s)\#\bEta$ is a transport plan from $\mu_t$ to $\mu_s$, so applying the definition of the 1-Wasserstein metric, Jensen's inequality,  Tonelli, and inequality (\ref{L1velocitybound}), we obtain
 \begin{align} \label{equicontinuityfirst} 
W_1(\mu_s,\mu_t) &\leq \int_{\Rd \times \Rd} |\pi^1-\pi^2| d  (e_s,e_t)\#\bEta \\
  &=   \int_{\Rd \times C([0,1];\Rd)} |e_s-e_t| d \bEta \nonumber \\
    &=   \int_{\Rd \times C([0,1];\Rd)} |\gamma(s)-\gamma(t)| d \bEta(x,\gamma) \nonumber \\
    &\leq   \int_{\Rd \times C([0,1];\Rd)} \int_s^t | \dot{\gamma}(r)|dr \ d \bEta(x,\gamma)  \nonumber \\
        &=    \int_{\Rd \times C([0,1];\Rd)}  \int_s^t | \bF(\gamma(r), \mu_{r})+ \bu(\gamma(r),r)| dr \  d \bEta(x,\gamma)  \nonumber \\
 &= \int_s^t \left\|\bF(\cdot, \mu_{r})+ \bu(\cdot,r) \right\|_{L^1(\mu_{r})} dr \nonumber  \\
 &\leq \int_s^t \left( C_F + 2C_F'M_1(\mu_{r}) + \int_\Rd \left| \bu(\cdot,r) \right| d \mu_{r}  \right) dr. \nonumber \end{align}

 In particular, for all $t \in [0,1]$, inequality (\ref{linearphibound2}) ensures that, for all $t \in [0,1]$,
 \begin{align*}
 M_1(\mu_t) &= W_1(\mu_t, \delta_0)  \leq W_1(\mu_0, \mu_t) + W_1(\mu_0, \delta_0)  \\
  &\leq (C_F + R_\phi + 1) + 2C_F' \int_0^t M_1(\mu_r)dr + \Psi(\bnu |\mu) + M_1(m_0) .
 \end{align*}
 Thus, by Gronwall's inequality, there exists $C>0$ depending on $C_F, R_\phi$, $\Psi(\bnu |\mu)$, and $M_1(m_0)$ so that
 \begin{align} \label{firstmomentbound} 
 M_1(\mu_t) \leq C( 1+ 2C_F't e^{2C_F't}) \leq C(1+2C_F'e^{2C_F'}) := C' \ , \ \forall t \in [0,1] , \ \forall n \in \mathbb{N}. 
 \end{align}
 Substituting this fact into inequality (\ref{equicontinuityfirst}), we obtain the result.
 \end{proof}

Next, we prove  Lemma \ref{compactnesseps}, which shows that sublevels of $\G$ in the constraint set $\C(m_0)$ are sequentially compact.
\begin{proof}[Proof of Lemma \ref{compactnesseps}]
We begin by observing that, since $\sup_{n} \G(\mu_n,\bnu_n) <+\infty$, we also have 
\begin{align} \label{Psibound} \sup_{n  } \Psi(\bnu_n|\mu_n) <+\infty 
\end{align}
and, for all $n \in \mathbb{N}$, there exists $\bu_n$ so that $d\bnu_n(y,t) = \bu_n(y,t)d\mu_{n,t}(y)  dt$ and
\begin{align} \label{Psiformula} \Psi(\bnu_n| \mu_n) = \int_0^1 \int_\Rd \psi(\bu_n(y,t)) d \mu_{n, t}(y) dt . \end{align}
Furthermore, by inequality (\ref{phipsirelation}) and Lemma \ref{fancydoubling},   we have
\begin{align} \label{linearphibound}
\sup_n  \int_0^1 \int_\Rd |\bu_n(y,t)| d \mu_{n, t}(y) dt \leq \sup_n  \int_0^1 \int_\Rd \phi(|\bu_n(y,t)|) d \mu_{n, t}(y) dt + R_\phi <+\infty.
\end{align}
 
First, we apply Arzel\'a-Ascoli to show the convergence of $\mu_n$, up to a subsequence. We begin by showing $\{\mu_{n, t} \}_{n \in \mathbb{N}, t \in[0,1]}$ is relatively sequentially compact with respect to $W_1$.
Since $(\mu_n, \bnu_n) \in \C(m_0)$ and $\sup_{n } \G(\mu_n,\bnu_n) <+\infty$, by \cite[Proposition 5.3]{fornasier2019mean}, there exists $C>0$ and a nonnegative, continuously differntiable, convex, and superlinear at $+\infty$ function $\theta$, depending on the choice of $\psi, \bF$ and $m_0$, as in Assumption \ref{asmp1}, and the value of $\sup_{n} \G(\mu_n,\bnu_n)$ so that 
\begin{align*}
\sup_{t \in [0,1], n \in \mathbb{N}} \int_\Rd \theta(|x|) d \mu_{n,t}(x) &\leq C \left( 1+ \int_\Rd \theta(|x|) m_0(x) dx \right) <+ \infty.
\end{align*}
Since $\theta$ is superlinear at $+\infty$, for all $\epsilon >0$, there exists $R$ so that $r \geq R$ ensures $\theta(r)/r \geq 1/\epsilon$ and 
\[ \int_{\Rd \setminus B_0(R)} |x| d \mu_{n, t}(x)  \leq \epsilon \int_{\Rd \setminus B_0(R)} \theta( |x|) d \mu_{n, t}(x) \leq \epsilon C \left( 1+ \int_\Rd \theta(|x|) m_0(x) dx \right) .\]
Thus, $\{\mu_{n, t} \}_{n \in \mathbb{N}, t \in[0,1]}$ have uniformly integrable first moments, so by  \cite[Proposition 7.1.5]{ambrosio2005gradient}, $\{\mu_{n, t} \}_{n \in \mathbb{N}, t \in[0,1]}$ is relatively sequentially compact with respect to $W_1$ convergence.
 
Next, we show equicontinuity of $\mu_{n,t}$.  By Proposition \ref{pdeconstraintprop}, there exists $C' >0$ so that, for all $n \in \mathbb{N}$, $s,t \in [0,1]$,
\begin{align} \label{equictssecond}
 W_1(\mu_{n,s},\mu_{n,t}) &\leq   C'(t-s) + \int_s^t \int_\Rd |\bu_n(\cdot, r)| d \mu_{n,r} dr .
 \end{align}
 
Since  $\phi$ is strictly increasing,  $\phi^{-1}:[0,+\infty) \to (0, +\infty)$ is well defined and strictly increasing. Since $\phi$ is superlinear at $+\infty$, we must have $\lim_{s \to +\infty} \phi^{-1}(s)/s = 0$. Therefore, we may use Jensen's inequality to  estimate the second term above by
\begin{align} \label{equictsthird}
\int_s^t \int_\Rd |\bu_n(\cdot, r)| d \mu_{n,r} dr &= (t-s) \phi^{-1} \circ \phi \left( \frac{1}{t-s} \int_s^t \int_\Rd |\bu_n(\cdot, r)| d \mu_{n,r} dr \right)  \\
&\leq (t-s) \phi^{-1} \left( \frac{1}{t-s} \int_s^t \int_\Rd \phi(|\bu_n(\cdot,r)|) d \mu_{n, r} dr \right) \nonumber \\
&\leq (t-s) \phi^{-1} \left(\frac{C''}{t-s} \right) \nonumber
\end{align}
where the last  inequality follow from (\ref{linearphibound}), for some $C'' >0$ independent of $n \in \mathbb{R}^d$ and $s, t \in [0,1]$. Combining (\ref{equictssecond}) and (\ref{equictsthird}) shows that $\{\mu_{n,t}\}_{n \in \mathbb{N}}$ are equicontinuous.
 
 Thus, by Arzel\'a-Ascoli, there exists $\mu \in C([0,1]; \P_1(\Rd))$ so that, up to a subsequence, 
 \begin{align} \label{muconvergence} \mu_n \rightarrow \mu \text{ in } C([0,1]; \P_1(\Rd)) . \end{align}
 In particular, this implies $\mu_0 = m_0$ and  $d\mu_{n,t} \otimes dt \to d \mu_t \otimes dt$ narrowly in $\P(\Rd \times [0,1])$.

 Now, we turn to the convergence of $\bnu_n$. By inequality (\ref{Psibound}),  equation (\ref{Psiformula}), and equation (\ref{muconvergence}), we may argue in a similar way to  \cite[Theorem 5.4.4]{ambrosio2005gradient} to obtain that there exists $\bu: \Rd \times [0,1] \to U$ so that, up to a subsequence, for all $\xi \in C^\infty_c(\Rd \times [0,1]; \Rd) $,
 \begin{align*}
 \lim_{n \to +\infty} \int_0^1 \int_\Rd \xi(y,t) \bu_n(y,t) d \mu_{n,t}(y) dt = \int_0^1 \int_\Rd \xi(y,t) \bu(y,t) d \mu_t(y) ,
 \end{align*}
 and 
 \begin{align} \label{Psilimitfinite}
 \int_0^1 \int_\Rd |\bu(y,t)| d \mu_t(y) dt \leq  \int_0^1 \int_\Rd \psi(\bu(y,t)) d \mu_t(y) dt +R_\phi+1 < +\infty .
  \end{align}
 Thus, defining $d \bnu(y,t) := \bu(y,t) d \mu_t(y) dt$, we see that, up to a subsequence,  $\bnu_n \to \bnu$ narrowly in $  \M( \Rd \times [0,1];U)$. 

 It remains to show that $(\mu,\bnu) \in \C(m_0)$. We have already shown that $\mu_0 = m_0$. Since $(\mu_n, \bnu_n)$ is a distributional solution of the continuity equation, in the sense of equation (\ref{seconddistributionalsolution}),  for all $n \in \mathbb{N}$ and $\varphi \in C^\infty_c(\Rd \times[0,1])$, we have
 \begin{align*}
\int_0^1 \int_{\Rd} \partial_t\varphi d\mu_{n,t}dt + \int_0^1 \int_{\Rd}\nabla \varphi \cdot \bF(\cdot,  \mu_{n,t})d\mu_{n,t}dt  + \int_0^1 \int_{\Rd}\nabla \varphi \cdot d\bnu_{n} = 0 .
\end{align*}
By the convergence in equation (\ref{convergencecompactness}), it is clear that we may pass to the limit in the first and third terms. For the second term, note that, since  $\bF$ is  uniformly continuous we can conclude that 
\[\int_0^1 \int_{\Rd}\nabla \phi \cdot [\bF(\cdot,  \mu_{n,t}) - \bF(\cdot,  \mu_t) ]d\mu_{n,t}dt \rightarrow 0. \]
Thus, we may likewise pass to the limit in the second term, since $(y,t) \mapsto \bF(y, \mu_t)$ is continuous and $d\mu_{n,t} \otimes dt \to d \mu_t \otimes dt$ narrowly in $\P(\Rd \times [0,1])$. This shows that $(\mu,\bnu)$ is a distributional solution of the continuity equation.

It remains to show $\mu \in AC([0,1];\P_1(\Rd))$.  By Proposition \ref{pdeconstraintprop}, there exists $C' >0$ so that, for all $n \in \mathbb{N}$, $s,t \in [0,1]$,
 \begin{align} \label{equictssecond3}
 W_1(\mu_s,\mu_t) &\leq  C'(t-s) + \int_s^t \int_\Rd |\bu(\cdot, r)| d \mu_{r} dr = \int_s^t \omega(r) dr
 \end{align}
for
 \[ \omega(r):=  C' + \int_\Rd |\bu(\cdot, r)| d \mu_{r} .\]
 By inequality (\ref{Psilimitfinite}), we have $\omega \in L^1(\Rd)$. Therefore $\mu \in AC([0,1];\P_1(\Rd))$.
 
 Finally, the fact that $\G(\mu,\bnu) \leq \sup_n \G(\mu_n,\bnu_n)$ is an immediate consequence of the lower semicontinuity proved in  Lemma \ref{lem:epcont}.
\end{proof}

 \section{Convergence as $\epsilon \to 0, \delta \to 0, N \to +\infty$} \label{epstozeroappendix}
 We begin with our proof of Proposition \ref{propgammaeptozero}, which shows $\Gamma$-convergence of the objective functionals from (\ref{eq:main_optimrrbb}) to (\ref{mainproblem}) as $\epsilon \to 0$.
   \begin{proof}[Proof of Proposition \ref{propgammaeptozero}]
First we prove part (\ref{firstpartep}). Without loss of generality, we may pass to a subsequence so that  
\[ \lim_{\epsilon \rightarrow 0}\E_{\epsilon}(\mu_\epsilon,\bnu_\epsilon)  = \liminf_{\epsilon \rightarrow 0}\E_{\epsilon}(\mu_\epsilon,\bnu_\epsilon)  < + \infty . \]
Thus $\F_\epsilon(\mu_{\epsilon,1})$, as defined in equation (\ref{Fepsdef}), must be bounded uniformly in $\epsilon$, which shows  that $\lim_{\epsilon \rightarrow 0} \|\mu_{\epsilon,1} - m_1\|_2 =0$. Since, $\mu_\epsilon \rightarrow \mu$ in $ C([0,1]; \P_1(\Rd))$, we have $\mu_{\epsilon,1}$ converges to $\mu_1$ in $\P_1(\Rd)$. Thus, by uniqueness of limits, $\mu_1 = m_1$, so $\F(\mu_1) = 0$.
Due to the lower semicontinuity of the functional $\G$ 
from Lemma \ref{lem:epcont}, we can conclude that
\[\E(\mu,\bnu) = \G(\mu,\bnu) \leq \liminf_{\epsilon \rightarrow 0}\E_{\epsilon}(\mu_\epsilon,\bnu_\epsilon). \]

Next, we show part  (\ref{secondpartep}). We may assume without loss of generality that $\E(\mu,\bnu)< +\infty$, so $\F(\mu_1) =0$ and $\mu_1 = m_1 \in L^2(\Rd)$. The result then follows, since $\F_\ep(\mu_1) \equiv 0$.
\end{proof}

We now apply this $\Gamma$-convergence result to prove Proposition \ref{thm:convep}, which shows that, as $\epsilon \to 0$, minimizers of (\ref{eq:main_optimrrbb}) converge to a minimizer of (\ref{mainproblem}), up to a subsequence.

\begin{proof}[Proof of Proposition \ref{thm:convep}]
First, note that by the feasibility of (\ref{mainproblem}), there exist $(\mu',\bu') \in \mathcal{C}(m_0)$ so that $\E(\mu',\bu')<+\infty$, so $\mu'_1 = m_1$. Combining this with the fact that $(\mu_\epsilon, \bnu_\epsilon)$ are minimizers of  (\ref{eq:main_optimrrbb}), we have
\[ \G(\mu_\epsilon, \bnu_\epsilon) \leq  \E_\epsilon(\mu_\epsilon,\bnu_\epsilon) \leq \E_\epsilon( {\mu}', {\bu}')  = \G( {\mu}', {\bu}') + \F_\epsilon( {\mu}_1') = \G( {\mu}', {\bu}') + 0 , \ \forall \epsilon >0.\]
Therefore, $\sup_{\epsilon} \G(  \mu_\epsilon, \bnu_\epsilon) < +\infty$. By Lemma \ref{compactnesseps}, there exists $(\mu,\bnu) \in \mathcal{C}(m_0)$ so that, up to a subsequence, (\ref{desiredconvergencepstozero}) holds. 

Furthermore, for any  $(\mu'',\bnu'') \in C(m_0)$, Proposition \ref{propgammaeptozero} and the fact that $(\mu_\ep,\bnu_\ep)$ are minimizers ensure that  
\[\E({\mu}'',{\bnu}'') \geq \limsup_{\ep \rightarrow 0} \E_\ep({\mu}'',{\bnu}'') \geq \liminf_{\ep \rightarrow 0} \E_\ep(\mu_{\ep},\bnu_{\ep}) \geq \E(\mu,\bnu) \]
Since $({\mu}'',{\bnu}'') \in \C(m_0)$ was arbitrary, this gives the result.
\end{proof}

Next we prove Theorem \ref{deltatozerothm}, on the convergence of minimizers of (\ref{eq:main_optimsrrbb2}) to minimizers of (\ref{eq:main_optimrrbb}) as $\delta \to 0$.

\begin{proof}[{Proof of Theorem \ref{deltatozerothm}}]
First, note that, by the assumption that (\ref{eq:main_optimrrbb}) is feasible, there exist $(\mu', \bnu') \in \mathcal{C}(m_0)$ so that $\E_\ep(\mu',\bnu') <+\infty$. Combining this with the fact that $(\mu_\delta, \bnu_\delta)$ are minimizers, we have
\[ \G(\mu_\delta, \bnu_\delta) \leq  \E_{\epsilon,\delta}(\mu_\delta,\bnu_\delta) \leq \E_{\epsilon, \delta}( {\mu}', {\bnu}')  =  \G( {\mu}', {\bnu}') +  \ep^{-1} \|k_\delta*\mu_1' - k_\delta*m_1 \|_{L^2(\Rd)}^2 . \]
Since $\E_\ep(\mu',\bnu')<+\infty$, $\mu_1' \in L^2(\Rd)$, and the right hand side is bounded uniformly in $\delta$. Thus, by Lemma \ref{compactnesseps}, there exists $(\mu,\bnu) \in \mathcal{C}(m_0)$ so that, up to a subsequence, (\ref{desiredconvergencepstozeroagain}) holds. 

    Furthermore, for any  $(\mu'',\bnu'') \in C(m_0)$, Proposition \ref{gammadelta} and the fact that $(\mu_\delta,\bnu_\delta)$ are minimizers ensure that  
\[\E_\ep({\mu}'',{\bnu}'') \geq \limsup_{\delta \rightarrow 0} \E_{\ep,\delta}({\mu}'',{\bnu}'') \geq \liminf_{\delta \rightarrow 0} \E_{\ep,\delta}(\mu_{\delta},\bnu_{\delta}) \geq \E(\mu,\bnu) \]
Since $({\mu}'',{\bnu}'') \in \C(m_0)$ was arbitrary, this gives the result.

\end{proof}

 Now, we prove Proposition \ref{discretesolutionsexist}, which develops sufficient conditions to ensure that a solution of (\ref{MFCN})   exists. Our proof is a mild adaptation of \cite[Proposition 4.2]{fornasier2019mean}, extending to the case when $L_N$ satisfies Assumption  \ref{MFCNas}(\ref{LNlsc}).
 \begin{proof}[Proof of Proposition \ref{discretesolutionsexist}]
 First, we consider feasibility of (\ref{MFCN}).
As observed in the paragraph following \cite[equation (3.8)]{fornasier2019mean}, classical results on existence of solutions to ordinary differential equations ensure that there exists $(\by,0) \in \mathcal{A}_N(\by_0)$ is nonempty. If   $L_N$ satisfies Assumption \ref{MFCNas}(\ref{LNcts}), then for  $(\by,0) \in \mathcal{A}_N(\by_0)$, $\E_{\ep,\delta,N}(\by,\bu)<+\infty$. 

On the other hand, suppose  $L_N$ satisfies Assumption  \ref{MFCNas}(\ref{LNlsc}), the control is unconstrained $U=\Rd$, and the initial particle locations are contained in $\{L_N<+\infty\}$. If we consider the curve $y_i(t) \equiv y_{i,0}$ and the velocity  $u_i(t) = -\bF_N(y_i(t),\by(t)) = -\bF_N(y_{i,0}, \by_0)$, then $(\by,\bu) \in \mathcal{A}_N(\by_0)$, and since $y_i(t) \in \{L_N <+\infty\}$ for all $t \in [0,1]$, we have $\E_{\epsilon,\delta,N}(\by,\bu)<+\infty$.

Now, suppose  (\ref{MFCN}) is feasible, and we will show that a minimizer exists. Take a minimizing sequence $(\by_k,\bu_k) \in \mathcal{A}_N(\by_0)$. Since 
\begin{align*}
\sup_k \frac{1}{N}\sum_{i=1}^N  \int_0^1 \psi(u_{i,k}(t))dt < +\infty ,
\end{align*}
as in \cite[Proposition 4.2]{fornasier2019mean}, we obtain that there exist $(\by,\bu) \in \mathcal{A}_N(\by_0)$ so that, up to a subsequence, $\by_k \to \by$ in $C([0,1];(\Rd)^N)$ and $\bu_k \to \bu$ in $L^1([0,1];U^N)$. Furthermore, the proof of \cite[Proposition 4.2]{fornasier2019mean} shows that
\begin{align*} &\liminf_{k \to +\infty} \frac{1}{N}\sum_{i=1}^N  \int_0^1 \psi(u_{i,k}(t))dt + \frac{1}{N} \sum_{i=1}^N \int_0^1  L_N(y_{i,k}(t),\mathbf{y}_k(t)) dt \\
&\qquad  \qquad \geq \frac{1}{N}\sum_{i=1}^N  \int_0^1 \psi(u_i(t))dt + \frac{1}{N} \sum_{i=1}^N \int_0^1  L_N(y_i(t),\mathbf{y}(t)) dt  . \end{align*}
By continuity of $K_\delta$ and $K_\delta*m_1$, we may likewise pass to the limit in $k$ in the third term in $\E_{\epsilon,\delta,N}$. This shows $ \liminf_{k \to +\infty} \E_{\epsilon,\delta,N}(\by_k,\bu_k) \geq \E_{\epsilon,\delta,N}(\by,\bu) . $
Since $(\by_k,\bu_k)$ was a minimizing sequence, this shows that $(\by,\bu)$ is a minimizer of (\ref{MFCN}).
 \end{proof}
 
Now, we prove Theorem \ref{Ntoinftyminimizers}, which shows that, for fixed $\epsilon, \delta >0$, as $N \to +\infty$, minimizers of the spatially discrete problem (\ref{MFCN}) converge to a solution of (\ref{eq:main_optimsrrbb2}), up to a subsequence. 

\begin{proof}[Proof of Theorem \ref{Ntoinftyminimizers}]
First, we will show that  there exists such a sequence $\by_{N,0}$. Let
\begin{align} \label{em0def}
e_{m_0} := \inf_{(\mu, \bnu) \in \C(m_0)} \E_{\epsilon,\delta}(\mu,\bnu) , \quad \E_{\epsilon,\delta}(\mu,\bnu):= \G(\mu,\bnu)+ \F_{\epsilon,\delta}(\mu_1) .
\end{align}
By our hypothesis that (\ref{eq:main_optimsrrbb2}) is feasible, we have that $e_{m_0}<+\infty$. Thus, for all $N \in \mathbb{N}$, there exists $(\mu_N,\bnu_N) \in C(m_0)$ so that $\E_{\epsilon,\delta}(\mu_N,\bnu_N) < e_{m_0} + 1/N$. By Proposition \ref{gammathm}(\ref{fulllimsup}), there exists $\by_{N,0}$ so that (\ref{whynot}-\ref{whynot1}) hold and there exist $(\tilde{\by}_N,\tilde{\bu}_N) \in \mathcal{A}_N(\by_{N,0})$ so that $\E_{\epsilon,\delta, N}(\tilde{\by}_N,\tilde{\bu}_N) \leq e_{m_0} + 2/N$. Now, suppose $(\by_N,\bu_N) \in \mathcal{A}_N(\by_{N,0})$ is an optimizer of (\ref{MFCN}). Then, we have 
\begin{align} \label{whyyes}
\limsup_{N\to +\infty} \E_{\epsilon,\delta, N}(\by_N,\bu_N)  \leq \limsup_{N\to +\infty} \E_{\epsilon,\delta, N}(\tilde{\by}_N,\tilde{\bu}_N) \leq e_{m_0} 
\end{align}
and
\begin{align*}
 \sup_{N} \frac{1}{N}\sum_{i=1}^N  \int_0^1 \psi(u_{i,N}(t))dt \leq  \sup_{N} \E_{\ep,\delta,N}( \by_N,\bu_N) \leq  \sup_{N} \E_{\ep,\delta,N}( \tilde{\by}_N, \tilde{\bu}_N) \leq e_{m_0}+1 <+\infty.
\end{align*}
Then \cite[Theorem 3.1]{fornasier2019mean} ensures that there exists $(\mu, \bnu) \in \mathcal{C}(m_0)$ so that (\ref{conv1})-(\ref{conv2}) hold. Furthermore, for any such limit point $(\mu,\bnu)$, by Proposition \ref{gammathm}(\ref{fullliminf}),
\begin{align} \label{whynotineq2}
e_{m_0} \leq    \E_{\epsilon, \delta}(\mu ,\bnu) \leq \liminf_{N \to +\infty} \E_{\ep,\delta,N}(\by_N, \bu_N) .
\end{align}
Combining (\ref{whyyes}) and (\ref{whynotineq2}), we obtain (\ref{goodinitialdata})
and that  $(\mu,\bnu)$ is a minimizer of (\ref{eq:main_optimsrrbb2}). 
\end{proof}

Now, we turn to the proof of Corollary \ref{existenceofsolutions}, which provides sufficient conditions to ensure existence of minimizers to the continuum optimization problems we consider.
\begin{proof}[Proof of Corollary \ref{existenceofsolutions}]
Since (\ref{mainproblem}) is feasible, we have that (\ref{eq:main_optimrrbb}) and (\ref{eq:main_optimsrrbb2}) are feasible for all $\epsilon, \delta>0$.  Likewise, Proposition \ref{discretesolutionsexist} ensures that, for any $\by_0 \in (\Rd)^{N}$ with 
\begin{align} \label{goodgood}
y_{i,0} \in \{L_N <+\infty\} \text{ for all } i = 1, \dots, N,
\end{align} minimizers of (\ref{MFCN}) exist.

By Theorem \ref{Ntoinftyminimizers}, for all $N \in \mathbb{N}$, there exists $\by_{N,0} \in (\Rd)^N$ satisfying (\ref{goodgood}) so that, for any sequence of minimizers $(\by_N, \bu_N) \in \mathcal{A}(\by_{N,0})$  of (\ref{MFCN}), up to a subsequence, $(\by_N, \bu_N)$ converges to a minimizer of (\ref{eq:main_optimsrrbb2}). Thus, minimizers of  (\ref{eq:main_optimsrrbb2}) exist.

By Theorem \ref{deltatozerothm}, any sequence of minimizers of (\ref{eq:main_optimsrrbb2}) converges, up to a subsequence, to a minimizer of (\ref{eq:main_optimrrbb}), so minimizers of (\ref{eq:main_optimrrbb}) exist.

Finally, by Proposition \ref{thm:convep}, any sequence of minimizers of (\ref{eq:main_optimrrbb}) converges, up to a subsequence, to a minimizer of (\ref{mainproblem}), so minimizers of (\ref{mainproblem}) exist.

\end{proof}
We conclude with the proof of our main theorem, Theorem \ref{maincoffeetheorem1}, on the convergence of minimizers of (\ref{MFCN}) to a minimizer of (\ref{mainproblem}) as $\epsilon \to 0, \delta \to 0$, and $N \to +\infty$.

\begin{proof}[Proof of Theorem \ref{maincoffeetheorem1}]
Recall that, as a consequence of our assumption that (\ref{mainproblem}) is feasible, we also have that (\ref{eq:main_optimrrbb})  and (\ref{eq:main_optimsrrbb2}) are feasible.

By Theorem \ref{NtoinftyminimizersRdLTI}, for any $\epsilon,\delta >0$,
 if $ (\by_{N}^{\epsilon,\delta}, \bu_{N}^{\epsilon,\delta}) \in \mathcal{A}_N(\by_{N,0}) $ minimizes  (\ref{MFCN}),
 there is a  subsequence $N_l^{\epsilon,\delta}$ of $N$  and   $(\mu_{\epsilon,\delta},\bnu_{\epsilon,\delta}) \in \C(m_0)$ that minimizes (\ref{eq:main_optimsrrbb2}) for which 
\begin{align} \label{conv1superfinal}
&\frac{1}{N_l^{\epsilon,\delta}}\sum_{i=1}^{N_l^{\epsilon,\delta}} \delta_{y^{\epsilon,\delta}_{i,N_l^{\epsilon,\delta}}} \xrightarrow{l\to +\infty} \mu_{\epsilon,\delta} \text{ in }C([0,1]; \P_1(\Rd))\\ 
&\frac{1}{N_l^{\epsilon,\delta}} \sum_{i=1}^{N_l^{\epsilon,\delta}} u^{\epsilon,\delta}_{i,N_l^{\epsilon,\delta}}(t) \delta_{y^{\epsilon,\delta}_{i,N_l^{\epsilon,\delta}}(t)} dt \xrightarrow{l\to +\infty} \bnu_{\epsilon,\delta} \text{ in } \M(\Rd \times [0,1]; U) . \label{conv2superfinal}
\end{align}
 
 Likewise, for any $\epsilon$, Theorem \ref{deltatozerothm} ensures that there exists a subsequence $\delta_k^\epsilon$ of $\delta$ so that $(\mu_{\epsilon,\delta_k^\epsilon}, \bnu_{\epsilon,\delta_k^\epsilon})$ converges to some $(\mu_{\epsilon}, \bnu_{\epsilon}) \in \mathcal{C}(m_0)$, where $(\mu_{\epsilon}, \bnu_{\epsilon})$ minimizes (\ref{eq:main_optimrrbb}).
Finally, by Theorem \ref{thm:convep}, up to a subsequence, $(\mu_{\epsilon}, \bnu_{\epsilon})$ converges to some $(\mu, \bnu)$, where $(\mu,\bnu)$ is a minimizer of (\ref{mainproblem}).

To obtain the result, we now apply a standard diagonal argument. For simplicity of notation, recall that convergence in the topologies of equations (\ref{coffee1}-\ref{coffee2}) is metrizable, and let $d$ denote such a metric. 
For $m \in \mathbb{N}$, we may choose a subsequence $(\mu_{\epsilon_{m}}, \bnu_{\epsilon_{m}})$ so that
\[ d\left((\mu_{\epsilon_{m}}, \bnu_{\epsilon_{m}}), (\mu,\bnu) \right) < \frac{1}{3m}. \]
Likewise, we may choose a subsequence $\left(\mu_{\epsilon_{m},\delta_{m}^{\epsilon_m}}, \bnu_{\epsilon_{m},\delta_{m}^{\epsilon_m}} \right)$ so that 
\[ d\left( \left(\mu_{\epsilon_{m},\delta_{m}^{\epsilon_m}}, \bnu_{\epsilon_{m},\delta_{m}^{\epsilon_m}} \right), (\mu_{\epsilon_{m}}, \bnu_{\epsilon_{m}}) \right) < \frac{1}{3m}. \]
Finally, we may choose a subsequence  $\left(\by^{\epsilon_m,\delta_m}_{N_{m}^{j_m,k_m}}, \bu^{\epsilon_m,\delta_m}_{N_m^{j_m,k_m}} \right)$ so that, defining $\mu_m$ and $\bnu_m$ as on the left hand side of   (\ref{cookies1}-\ref{cookies2}),
\[ d\left((\mu_{m},\bnu_m), \left(\mu_{\epsilon_{m},\delta_{m}^{\epsilon_m}}, \bnu_{\epsilon_{m},\delta_{m}^{\epsilon_m}} \right) \right) < \frac{1}{3m}. \]
The result then follows by the triangle inequality.
\end{proof}

\textbf{Acknowledgements:} K. Craig would like to thank Amir Sagiv for an interesting discussion on the case of optimal transport around obstacles and the connection to manifolds with holes.

\bibliographystyle{siamplain} 
\bibliography{blobOTref}
\end{document}